\documentclass{amsart}
\usepackage{graphicx}
\usepackage{subfigure}
\usepackage{amssymb,amsthm, amsmath}
\usepackage{mathrsfs}
\usepackage{bm}
\usepackage[monochrome]{color}
\usepackage{stmaryrd}
\usepackage{epic,overpic}
\usepackage{float}
\usepackage{tikz}
\usetikzlibrary{arrows}
\usepackage[hidelinks]{hyperref}

\def\N{\mathbb N}
\def\Z{\mathbb Z}
\def\Q{\mathbb Q}
\def\R{\mathbb R}

\def\P{\mathscr P}

\def\A{\mathcal{A}} 
\def\AA{\mathscr{A}} 
\def\T{\bm{T}}
\def\TT{\mathscr{T}}
\def\eps{\varepsilon}
\def\supp{\mathrm{supp}}
\def\Per{\mathrm{Per}}
\def\H{\mathscr H}
\def\Card{\mathrm{Card}}

\newcommand{\SL}{\mathfrak{L}}
\DeclareMathOperator{\Nut}{Nut}

\newcommand{\SLT}{\mathrm{SLT}}
\newcommand{\SLH}{\mathrm{SLH}}
\DeclareMathOperator{\st}{st}
\DeclareMathOperator{\VC}{VC}
\DeclareMathOperator{\Isom}{Isom}

\newcommand{\MH}[1]{\mathrm{MH}_{#1}}
\newcommand{\AT}{\AA_{\TT}}
\newcommand{\TA}{\TT_{\AA}}
\newcommand{\TTH}{\TT_{\mathrm{B}}}
\newcommand{\BDsim}{\stackrel{\mathrm{BD}}{\sim}}
\newcommand{\sfteq}{\stackrel{\mathrm{sft}}{=}}

\newcommand{\str}[1]{.\bm{#1}\phantom{.}}
\newcommand{\Triplet}[3]{\left\langle{#1}, {#2}, {#3}\right\rangle}
\newcommand{\abs}[1]{\left|{#1}\right|}
\newcommand{\floor}[1]{\left\lfloor{#1}\right\rfloor}
\newcommand{\ceil}[1]{\left\lceil{#1}\right\rceil}
\newcommand{\round}[1]
{\left\llbracket{#1}\right\rrbracket}

\renewcommand{\setminus}{\smallsetminus}
\definecolor{HMDcomment}{rgb}{0,0,0.5}

\newtheorem{thm}{Theorem}
\newtheorem{prop}{Proposition}
\newtheorem{cor}{Corollary}

\newtheorem{lem}{Lemma}
\newtheorem{rem}{Remark}
\newtheorem{ex}{Example}

\theoremstyle{definition}
\newtheorem{definition}{Definition}

\begin{document}
\raggedbottom
\title{Sturmian lattices and Aperiodic tile sets}
\address{Institute of Mathematics, University of Tsukuba}
\author{Shigeki Akiyama}
\email{akiyama@math.tsukuba.ac.jp}
\author{Tadahisa Hamada}
\email{taratta@jca.or.jp}
\author{Katsuki Ito}
\email{katsuki@math.tsukuba.ac.jp}
\date{}

\begin{abstract}
We give an explicit algorithm to construct
aperiodic tile sets based on Sturmian words of quadratic slopes. 
The method works for any quadratic irrational slope, and 
we can produce infinitely many aperiodic tile sets 
whose underlying scaling constant is a unit of 
any real quadratic field.
There are two key ingredients in our construction. 
The first one is ``Sturmian lattices''; 
an interesting grid structure generated by Sturmian words that 
emerged in an aperiodic monotile called Smith Turtle found by Smith, Myers, Kaplan and Goodman-Strauss.
We shall give a classification of Sturmian lattices.
The second is the bounded displacement equivalence of Delone sets, 
which plays a central role in this construction. 
\end{abstract}
\maketitle

\tableofcontents

\section{Introduction}
\label{sec:Intro}

Smith Turtle found in \cite{SMKGS:23_1} is an aperiodic monotile, i.e., it tiles the plane but only in non-periodic ways. 
There are several proofs of this fact and our starting point is a proof
in \cite{Akiyama-Araki:23}.
We draw the Ammann bars on the tile as in Figure~\ref{fig:Turtle}.

\begin{figure}[htb]\centering
\includegraphics[pagebox = artbox, 
width = .5\linewidth, page = 9]
{Example3.pdf}
\caption{Ammann bars on Smith Turtle}
\label{fig:Turtle}
\end{figure}

\begin{figure}[htb]\centering
\includegraphics[
width = \linewidth, page = 1]
{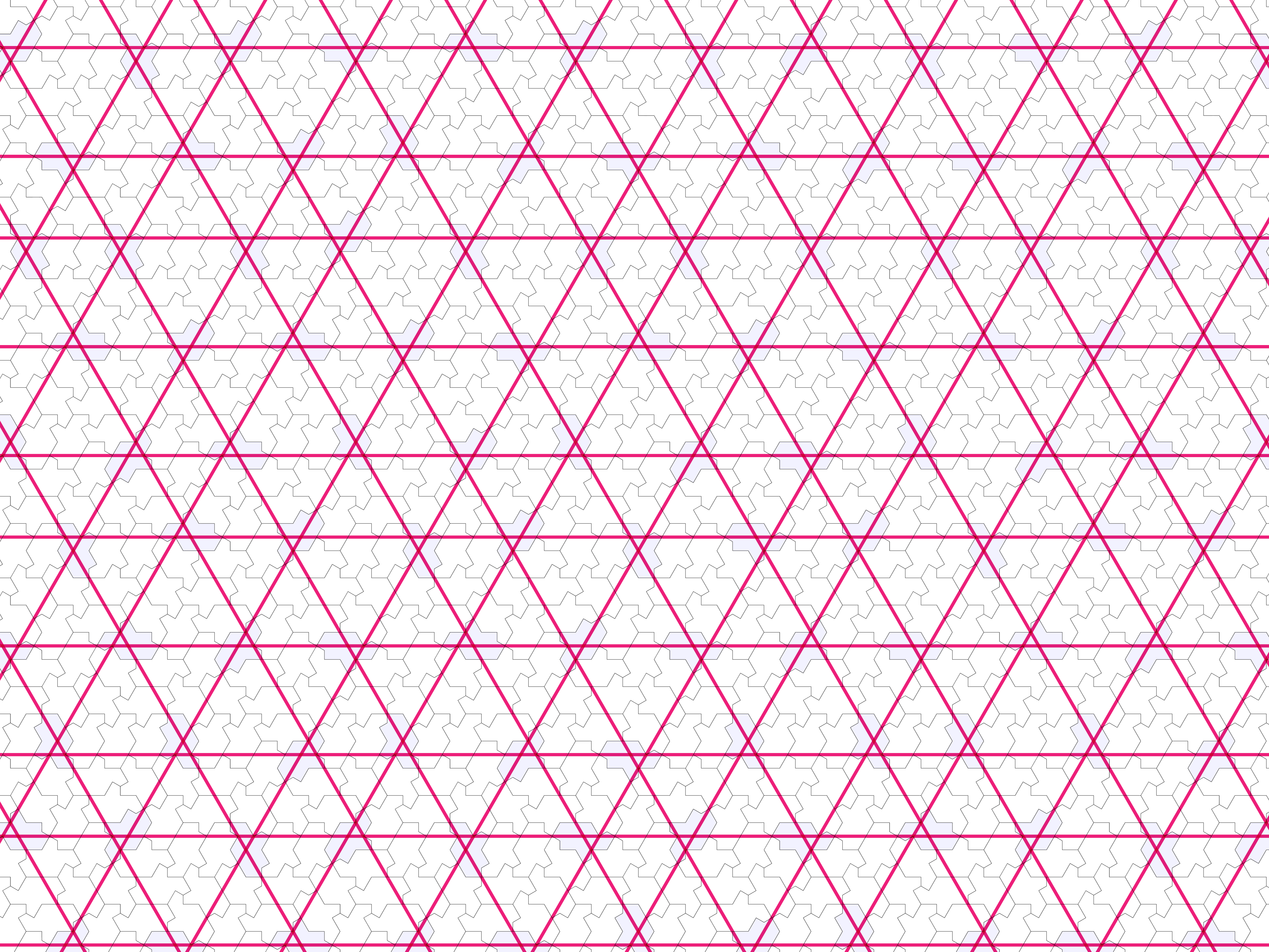}
\caption{Sturmian lattice on a Turtle tiling}
\label{fig:intro}
\end{figure}

Assuming patches extend to a tiling, 
the bars must continue across the edge and form a
straight line, enforcing a lattice-like structure as in Figure~\ref{fig:intro}. This structure may look simple, but it is not. 

Whenever two lines meet, a third line of the remaining slope passes close by forming a tiny regular triangle.
Fixing one line, 
the tiny triangles appear above or below this line
without an obvious rule.
In this paper, we define the {\it Sturmian lattice} by 
generalizing this configuration and give a classification in \S \ref{sec:main}
in Theorems~\ref{thm:irrational}--\ref{thm:others}.
 It turned out that the distribution of widths between adjacent parallel lines forms balanced sequences, i.e., these Sturmian lattices already enforce the
structure of the Sturmian words but include rational slopes.
We also show in \S \ref{sec:SSL}
that Sturmian lattices have natural hierarchical structures 
which are intimately related to the natural extension of the continued fraction algorithm.

The proof in \cite{Akiyama-Araki:23} 
relies on two facts:
(i) the tile and its flipped version carry Ammann bars of length $4:1$, and
(ii) flipped tiles are in one-to-one correspondence with the tiny triangles that
appear at the crossings.
Indeed, 
this already implies that the slope corresponding to the
Sturmian lattice must be irrational, implying that the resulting tiling cannot be periodic.
Further, we see that adjacent gaps of 
the parallel Ammann bars correspond to the Sturmian sequence of 
slope $(-1 + \sqrt{5})/2$ or $(3 - \sqrt{5})/2$ in Corollary~\ref{cor:Turtle-slope} in \S \ref{sec:symb}.
Using the above classification
of Sturmian lattices, 
we wish to generalize this idea 
 to an arbitrary
Sturmian lattice. 
Our main Theorem~\ref{thm:Voronoi-aperiodic} in \S \ref{sec:Voronoi} reads
\begin{quote}\itshape
For the Sturmian lattice of any quadratic irrational slope, 
we can construct an aperiodic tile set; a finite set of tiles
that can tile the plane but only in non-periodic ways.
\end{quote}
A Sturmian word of a quadratic slope $\alpha$
is generated by a certain 
substitution whose expansion constant 
is a unit of a quadratic field $\Q(\alpha)$, see \S\ref{sec:exp-const}.
Therefore the corresponding 
Sturmian lattice has a self-similar structure to this quadratic
unit expansion.

Theorem~\ref{thm:Voronoi-aperiodic} provides infinitely many essentially 
different aperiodic tile sets, since their 
underlying expansion constants are the units 
in arbitrarily chosen real quadratic fields.

A common historical way to show the
aperiodicity of a given tile set 
is to show that it possesses the hierarchical unique composition property.
S.~Mozes gave an idea to construct 
an aperiodic Wang tile set in \cite{Mozes:89}
from tilings generated by a substitution rule. 
There are several attempts along this line, e.g., 
\cite{Goodman-Strauss:98},
\cite{Gloannec-Ollinger:12},
\cite{Bedaride-Fernique:15}. 
However 
the proof often becomes intricate and requires computer assistance. 
Here we propose a quite different method. 
We compare frequencies of 
local patterns which appear in the Sturmian lattice and
enforce some algebraic constraints among them.
The proof of aperiodicity is mathematically pretty simple.
Everything boils down to a simple statement:
\begin{quote}\itshape
A quadratic curve and a line intersect in at most two points.
\end{quote}
It is noteworthy that 
we do not use the self-similarity of Sturmian lattices in the 
proof of Theorem~\ref{thm:Voronoi-aperiodic}. 
In fact, the tile sets we obtain need not enforce self-similarity of the
tilings they admit: a substitution tiling space associated with a
primitive invertible tiling substitution has zero topological entropy with
respect to the translation action
\cite[Theorem~7.16]{Robinson:04},
whereas we construct a tile set whose associated tiling space has positive
topological entropy; see \S\ref{sec:sqrt2-2}.

Sturmian lattices can be naturally viewed as tilings.
Our key tool is bounded displacement equivalence among Delone sets
to find an appropriate
 correspondence among tiles emerging in the tiling by a Sturmian lattice.
The proof of Theorem~\ref{thm:Voronoi-aperiodic} depends on the
result of Laczkovich \cite{Laczkovich:90,Laczkovich:92}. 

To make this proof into an algorithm, we introduce the \emph{cross-BD pair} which was inspired by the idea of Duneau~and~Oguey~\cite{DO:91} to reduce the problem to 
a one-dimensional one.
Our method is simpler than the one in \cite{DO:91} and gives a BD correspondence between $2$-dimensional lattices, 
where the ratio of the areas of their fundamental domains is an integer
We fix a particular way of construction of aperiodic tile sets using this BD equivalence
and show 
in Theorem~\ref{thm:lambda} in \S \ref{sec:ex1}
that:
%
\begin{quote}\itshape
The number of tiles of the aperiodic tile set generated by 
a Sturmian lattice of quadratic slope
is bounded by a constant multiple of its expansion constant.
\end{quote}
There are many possible ways of constructing BD equivalences. 
We do not know the best way to
minimize the number of tiles with respect to the choices of slopes.
Also by our construction, the tiles may not be topological disks; even worse,
they could be disconnected in general. 
It is an intriguing remaining problem to minimize the number of tiles
having nice topological properties. 
Our method might be applicable to finding a new 
aperiodic monotile of a different type.
Apart from the above rough estimate, we tried to minimize the number of tiles, 
using the self-similar structure of the Sturmian lattices in \S \ref{sec:SSL}.
We find that there is an aperiodic tile set 
consisting of 3 tiles; all of them are homeomorphic to a ball, and
one of its resulting tilings is self-similar 
with an expansion constant $\sqrt{2}-1$.
We also discuss how Smith Turtle can be retrieved by our method in \S \ref{sec:exTurtle}.
The construction of aperiodic tile sets using self-similarity seems handy, 
but we are not sure if this way works or not 
for the general quadratic slopes. 

The logical dependencies among the sections are summarized in
Figure~\ref{fig:roadmap}.
Our main theorem, Theorem~\ref{thm:Voronoi-aperiodic} in \S\ref{sec:Voronoi},
rests on the tiling framework of \S\ref{sec:aperiodic} and the BD
correspondences of \S\ref{sec:BD}, together with the classification of
\S\ref{sec:main}.
For a shortest route to the proof, the reader may follow
\S\ref{sec:bal}, \S\ref{sec:SL}, \S\ref{sec:main}, \S\ref{sec:aperiodic},
\S\ref{sec:BD}, and \S\ref{sec:Voronoi}, in this order.

\begin{figure}[H]
\centering
\hyphenpenalty=10000\exhyphenpenalty=10000
\begin{tikzpicture}[>=stealth,
  sbox/.style={draw=blue, rounded corners, align=center, text width=3cm,
    font=\footnotesize, inner sep=3.5pt, minimum height=0.85cm},
  mbox/.style={draw=blue, very thick, rounded corners, align=center, text width=3.2cm,
    font=\footnotesize, inner sep=3.5pt, minimum height=0.85cm},
  flow/.style={->, draw=blue, thick}
]
\node[sbox] (s3)  at (0,0)     {\S3 \textbf{Sturmian lattices}\\{\scriptsize geometric axiom $\Leftrightarrow$ symbolic coding}};
\node[sbox] (s4)  at (0,-1.7)  {\S4 \textbf{Classification}\\{\scriptsize Thms \ref{thm:irrational}--\ref{thm:others}}};
\node[mbox] (s8)  at (0,-4.1)  {\S8 \textbf{Nuts and Bolts}\\{\scriptsize \textbf{Main Thm \ref{thm:Voronoi-aperiodic}}}\\{\scriptsize SL $\Leftrightarrow$ Nut/Bolt tilings}};
\node[sbox] (s9)  at (0,-6.3)  {\S9 \textbf{Examples I}\\{\scriptsize Thm \ref{thm:lambda}}};
\node[sbox] (s10) at (0,-8.0)  {\S10 \textbf{Examples II}\\{\scriptsize $3$ tiles; Smith Turtle}};
\node[sbox] (s2) at (-4.4,0)     {\S2 \textbf{Balanced words}\\{\scriptsize Markoff--Reutenauer}};
\node[sbox] (s6) at (-4.4,-2.4)  {\S6 \textbf{Aperiodic tile sets}\\{\scriptsize tiles; Ammann bars}};
\node[sbox] (s7) at (-4.4,-4.1)  {\S7 \textbf{BD correspondences}\\{\scriptsize cross-BD pair (Laczkovich)}};
\node[sbox] (s5) at (4.4,-1.7) {\S5 \textbf{Super Sturmian lattices}\\{\scriptsize self-similarity; fundamental SL}};
\draw[flow] (s2) -- (s3);
\draw[flow] (s2.south) to[out=-90,in=170] (s4.west);
\draw[flow] (s3) -- (s4);
\draw[flow] (s3) to[out=0,in=170] (s5.west);
\draw[flow] (s4) -- (s8);
\draw[flow] (s6) -- (s8);
\draw[flow] (s7) -- (s8);
\draw[flow] (s8) -- (s9);
\draw[flow] (s8.west) to[out=-160,in=150] (s10.west);
\draw[flow] (s5.south) to[out=-80,in=0] (s10.east);
\end{tikzpicture}
\caption{Logical roadmap of the paper. Arrows point from prerequisites to the sections that depend on them.}
\label{fig:roadmap}
\end{figure}

\section{Balanced words}
\label{sec:bal}

Let $\A:= \{0, 1\}$ be the set of letters.
Denote by $\A^{*}$ the monoid generated by $\A$ by concatenation with the identity $\lambda$, the empty word.
Let $\A^{\Z}$ be the set of bi-infinite words generated by $\A$.
We usually use a notation $w = (w_{n})_{n\in \Z} = \cdots w_{-1}w_{0}w_{1}\cdots$ for a bi-infinite word, and  if necessary we also write $w = \cdots w_{-1}.w_{0}w_{1}\cdots$ with a decimal point.
For word $w\in \A^{\Z}$, we denote by
\begin{align*}
w_{[a,b]}&:= w_{a}w_{a+1}\cdots w_{b},&
w_{[a,\infty)}&:= w_{a}w_{a+1}\cdots,&
w_{(-\infty,b]}&:= \cdots w_{b-1}w_{b}.
\end{align*}
We often write $w_{[a,b)}:= w_{[a,b-1]}$.
This notation is compatible with the additivity: $w_{[a,b)}w_{[b,c)} = w_{[a,c)}$.
The \emph{mirror word} $\tilde{w}$ of a word $w$ is
\[
\tilde{w}:= \begin{cases}
w_{b}w_{b-1}\cdots w_{a} & \text{if }w = w_{[a,b]}\\
\cdots w_{a+1}w_{a} & \text{if }w = w_{[a,\infty)}\\
w_{b}w_{b-1}\cdots & \text{if }w = w_{(-\infty,b]}\\
\cdots w_{1}w_{0}.w_{-1}\cdots & \text{if }w = \cdots w_{-1}.w_{0}w_{1}\cdots
\end{cases}.
\]
For $w = w_{1}w_{2}\cdots w_{\ell}\in \A^{*}$, denote by $|w|:= \ell$ the \emph{length} of $w$, and $|w|_{1}$ be the number of occurrences of $1$ in $w$, which is called the \emph{height} of $w$.
A word $u\in \A^{*}$ is a \emph{factor} of $w\in \A^{*}\cup \A^{\Z}$ if $u$ appears as a subword of $w$, and we denote by $u\prec w$ (or $w\succ u$).
A word $w\in \A^{*}\cup \A^{\Z}$ is said to be \emph{$C$-balanced} if
\[
-C\le |u|_{1} - |v|_{1}\le C
\]
for all $u, v\prec w$ with $|u| = |v|$, and $w$ is said to be \emph{strictly $C$-balanced} if $w$ is $C$-balanced but not $(C-1)$-balanced.
If $w\in \A^{\Z}$ is $C$-balanced, then it is plain to see that there exists a limit
\[
\lim_{N\to \infty}
\frac{1}{2N}\abs{w_{[-N,N)}}_{1},
\]
called the \emph{slope} of $w$.

A word $w\in \A^{\Z}$ is \emph{positive periodic}, if there exist $u\in \A^{*}$ and $\ell\in \Z$ that
\[
w = w_{(-\infty,\ell]}uuu\cdots =: w_{(-\infty,\ell]}u^{\infty}.
\]
Similarly it is \emph{negative periodic} if
\[
w = \cdots uuuw_{[\ell,\infty)} =:{}^{\infty}uw_{[\ell,\infty)}.
\]
A word $w\in \A^{\Z}$ is \emph{periodic} if there exists $u\in \A^{*}$ and
\[
w = \cdots uuu\cdots =: {}^{\infty}u^{\infty} = u^{\Z}.
\]
By \cite[Lemma~6.2.2]{Fogg}, for a $1$-balanced word $w\in \A^{\Z}$, positive periodicity implies negative periodicity and vice versa.

A word $w\in \A^{\Z}$ is \emph{lower mechanical} if there exists $\alpha, \rho\in \R$ that
\[
w_{n} = \floor{(n+1)\alpha + \rho} - \floor{n\alpha + \rho}
\]
for all $n\in \Z$.
Similarly it is \emph{upper mechanical} if
\[
w_{n} = \ceil{(n+1)\alpha + \rho} - \ceil{n\alpha + \rho}
\]
for all $n\in \Z$.
A word $w$ is \emph{mechanical} if it is lower mechanical or upper mechanical, and $w$ is \emph{Sturmian} if it is a mechanical word with an irrational $\alpha$.
Note that necessarily $\alpha\in [0, 1]$, and then the slope of a mechanical word $w$ coincides with $\alpha$.

The theorem of Morse--Hedlund reads that a word $w\in \A^{\Z}$ is Sturmian if and only if it is a $1$-balanced and not positive periodic, see \cite[Chapter~6]{Fogg}.
Indeed, the classification of $1$-balanced words in $\A^{\Z}$ dates back to Markoff~\cite{Markoff:79}.
We use the classification detailed in \cite{Reutenauer:06}.
If the slope is irrational, then it can be written by a mechanical word.
It is further classified into two cases:
$(\MH2)$ when $\rho\notin \Z + \alpha\Z$ and $(\MH3)$ when $\rho\in \Z + \alpha\Z$.
In the case $(\MH2)$, the upper and lower mechanical words are identical.
If the slope is rational, then they are classified into periodic case $(\MH1)$ and skew case $(\MH4)$.
The case $(\MH4)$ is not periodic but both positive and negative periodic.
We shall recall the precise shape in Lemma~\ref{lem:Markoff}.
Reutenauer~\cite{Reutenauer:06} made clear that this fourfold classification measures the maximal symmetry of the balanced words with respect to the word $01$ or $10$.

A word $c\in \A^{*}$ is \emph{central} if the words $0c0$, $0c1$, $1c0$, and $1c1$ are all $1$-balanced, see \cite[Proposition~2.2.34]{Lothaire:02}.
Note that a central $c$ is a palindrome (i.e., $c = \tilde{c}$) and the word $0c1$ (resp.\ $1c0$) is called the \emph{lower} (resp.\ \emph{upper}) \emph{Christoffel word}.
For an irreducible fraction $\alpha = p/q\in (0, 1)\cap \Q$, the lower Christoffel word is explicitly written as
\[
\bigl(\floor{(n+1)\alpha} - \floor{n\alpha}\bigr)_{n = 0, 1, \dots, q-1}.
\]
The upper Christoffel word is obtained in the same manner by replacing $\floor{\cdots}$ with $\ceil{\cdots}$.
The bi-infinite words $(0c1)^{\Z}$ and $(1c0)^{\Z}$ are $1$-balanced of a rational slope $\alpha = p/q$.
By abuse of notation, we call this $\alpha$ the \emph{slope} of the central word.

\begin{lem}[Markoff]
\label{lem:Markoff}
Let $w\in \A^{\Z}$ be a $1$-balanced word of a rational slope $\alpha$ and $c$ the central word of slope $\alpha$.
Then $w$ must be one of the following three forms:
\begin{enumerate}
\item
$(0c1)^{\Z}$,
\item
${}^{\infty}(0c1)0c0(1c0)^{\infty}$,
\item
${}^{\infty}(1c0)1c1(0c1)^{\infty}$.
\end{enumerate}
\end{lem}

\begin{proof}
In the notation of \cite{Reutenauer:06}, the cases of rational slope are classified into the periodic case $(\MH1)$, and the skew case $(\MH4)$.
Here (2) and (3) are skew cases.
For the convenience of the reader, we give a brief proof due to \cite[Proposition~5.5]{Akiyama-Kaneko:21}.
The key observation is that the map $\phi$ which decreases run lengths of the fixed letter by one, i.e.,
\[
\phi(\cdots 10^{u(-1)}10^{u(0)}10^{u(1)}1\cdots) =
\cdots 10^{u(-1)-1}10^{u(0)-1}10^{u(1)-1}1\cdots
\]
or
\[
\phi(\cdots 01^{u(-1)}01^{u(0)}01^{u(1)}0\cdots) =
\cdots 01^{u(-1)-1}01^{u(0)-1}01^{u(1)-1}0\cdots
\]
preserves the balancedness, provided it is well-defined.
By repeated applications of $\phi$, the $1$-balanced word of rational slope eventually falls into
\[
\{0^{\Z}, 1^{\Z}, {}^{\infty}010^{\infty}, {}^{\infty}101^{\infty}\},
\]
i.e., the one of slope of $0$ or $1$.
This is nothing but a slow continued fraction algorithm (Farey map) applied on rational numbers.
Taking successive inverse images of $\phi$, we obtain this result.
Here we used the fact that $(c01)^{\Z}$ and $(c10)^{\Z}$ are shift equivalent.
One can deduce this fact by using the decomposition of central words $c = p01q$ by two central words $p$, $q$, see \cite[Corollary~2.2.9]{Lothaire:02} and
\cite[Appendix]{Akiyama-Kaneko:21}.
\end{proof}

\section{Sturmian lattices}
\label{sec:SL}

In this section we study the lattice-like structures called ``Sturmian lattices'', which appear in the Turtle tilings shown in Figure~\ref{fig:intro}.
In \S\ref{sec:geom} we introduce them through a geometric axiom and derive their fundamental properties.
In \S\ref{sec:symb} we recast Sturmian lattices in combinatorial terms and study their balancedness.

\subsection{Geometric Sturmian lattices}
\label{sec:geom}

\begin{figure}[htb]\centering
\includegraphics[width = .7\linewidth]
{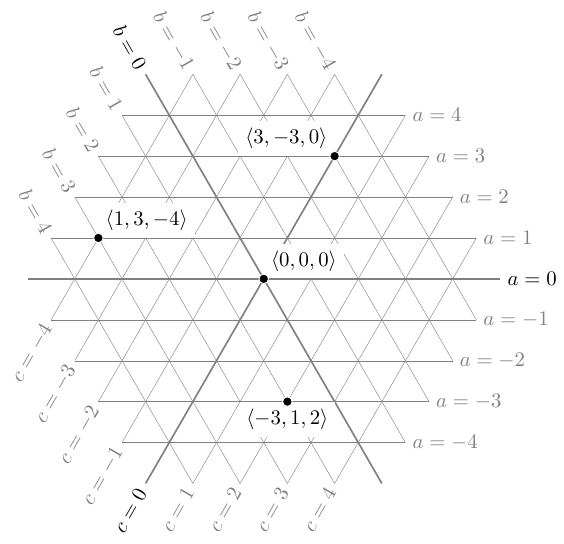}
\caption{Triangular coordinates}
\label{fig:triplet}
\end{figure}

We consider a collection of infinitely many lines in the Euclidean plane $\R^{2}$; there are three directions.
We introduce three axes labeled by $a$, $b$, and $c$ and encode the position of each line as a one-dimensional coordinate.
In the Cartesian $xy$-coordinate system, we define straight lines
\begin{align*}
\{a = \eta_{0}\}&:=
\{y = \eta_{0}\},\\
\{b = \eta_{1}\}&:=
\{y = -\sqrt{3}x - 2\eta_{1}\},\\
\{c = \eta_{2}\}&:=
\{y = \sqrt{3}x - 2\eta_{2}\}
\end{align*}
for $\eta_{0}, \eta_{1}, \eta_{2}\in \R$.
Note that three lines $a = \eta_{0}$, $b = \eta_{1}$, and $c = \eta_{2}$ form a regular triangle of height $\bigl|\eta_{0} + \eta_{1} + \eta_{2}\bigr|$.
Moreover, if $\eta_{0} + \eta_{1} + \eta_{2} > 0$ then the corresponding triangle is downward ($\nabla$-like), and if negative then it is upward ($\Delta$-like).
Of course $\eta_{0} + \eta_{1} + \eta_{2} = 0$ means that the three lines intersect at a single point, denoted by $\Triplet{\eta_{0}}{\eta_{1}}{\eta_{2}}\in \R^{2}$.
Figure~\ref{fig:triplet} helps us understand the setting.

From a triple $(a(i), b(j), c(k))$ of bi-infinite real sequences $(a(i))_{i\in \Z}$, $(b(j))_{j\in \Z}$, $(c(k))_{k\in \Z}$, we obtain a family of lines in the three directions,
\begin{equation}
\label{eq:geom}
\mathfrak{A}\cup
\mathfrak{B}\cup
\mathfrak{C}:=
\{a = a(i)\mid i\in \Z\}\cup
\{b = b(j)\mid j\in \Z\}\cup
\{c = c(k)\mid k\in \Z\}.
\end{equation}
The structure observed in Figure~\ref{fig:intro} is formalized as follows.

\begin{definition}
\label{def:geom}
A family of lines $\mathfrak{A}\cup \mathfrak{B}\cup \mathfrak{C}$ is called a \emph{geometric Sturmian lattice} if there exist a triple of real sequences $(a(i), b(j), c(k))$ that defines it in the sense of \eqref{eq:geom} and a constant $\eps > 0$ such that
\begin{align}
\label{eq:def:geom}
a(i+1) - a(i)&\ge 2\eps,&
b(j+1) - b(j)&\ge 2\eps,&
c(k+1) - c(k)&\ge 2\eps
\end{align}
and
\begin{equation}
\label{eq:def:geom-axiom}
i + j + k = 0\implies
\abs{a(i) + b(j) + c(k)} = \eps.
\end{equation}
\end{definition}

Condition \eqref{eq:def:geom-axiom} in Definition~\ref{def:geom} states that every three lines $a=a(i)$, $b=b(j)$, $c=c(k)$ with $i+j+k=0$ bound a triangle of one fixed size, and that these triangles are arranged in a lattice pattern.
We shall call such a triangle \emph{tiny}.
Condition \eqref{eq:def:geom} then asserts that these tiny triangles never overlap one another.
Since, for any $r > 0$, the scalar multiple
\[
r\cdot (a(i), b(j), c(k)):=
(r\cdot a(i), r\cdot b(j), r\cdot c(k))
\]
of a geometric Sturmian lattice $(a(i), b(j), c(k))$ is again a geometric Sturmian lattice, we fix $\eps = 1/2$ throughout, except in \S\ref{sec:SSL}.
We write
\[
(a(i), b(j), c(k))\propto
(a'(i), b'(j), c'(k))
\]
if there exists $r > 0$ that $a'(i) = r\cdot a(i)$, $b'(j) = r\cdot b(j)$, and $c'(k) = r\cdot c(k)$.

\begin{figure}[htb]\centering
\subfigure[Example~\ref{ex:Kagome} (Kagome lattice)]{%
\includegraphics[width = .48\linewidth]{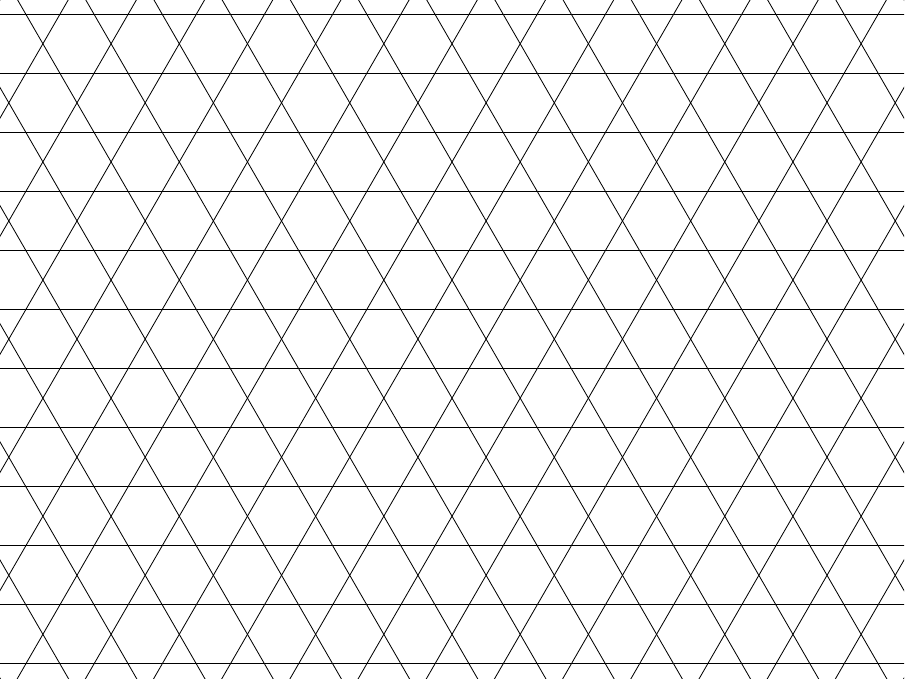}
}\hfill%
\subfigure[Example~\ref{ex:periodic}]{%
\includegraphics[width = .48\linewidth]{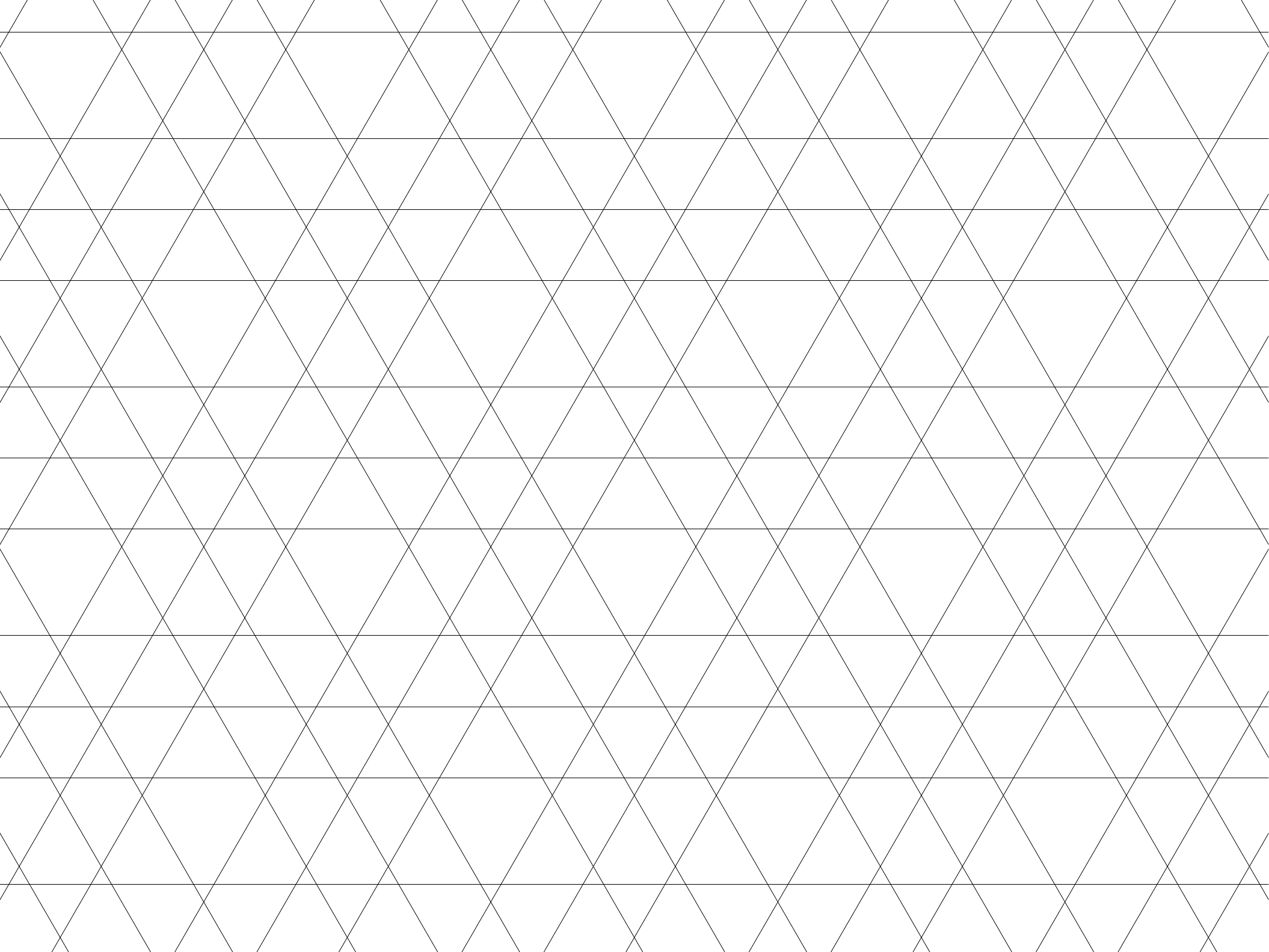}
}
\caption{Examples of periodic Sturmian lattices}
\label{fig:example}
\end{figure}

By abuse of terminology, we also call the defining triple $(a(i), b(j), c(k))$ itself a geometric Sturmian lattice.
We give a few examples of geometric Sturmian lattices:

\begin{ex}
\label{ex:Kagome}
\begin{align*}
a(i)&= \frac{1}{2} + i,&
b(j)&= j,&
c(k)&= k
\end{align*}
yields the Kagome lattice in Figure~\ref{fig:example}(a).
\end{ex}

\begin{ex}
\label{ex:periodic}
For $r\in \Z$ and $i, j, k \in \{0, 1, 2\}$, setting
\begin{align*}
a(3r - i)&= \frac{1}{2} + 7r - 2i,&
b(3r + j)&= 7r + 2j,&
c(3r + k)&= 7r + 2k
\end{align*}
yields a geometric Sturmian lattice in Figure~\ref{fig:example}(b), which is periodic in each of the three directions.
\end{ex}

The dihedral group
\[
D_{12} = \langle
\sigma, \tau\mid
\sigma^{6} = \tau^{2} = 1,~
\tau\sigma\tau = \sigma^{-1}
\rangle
\]
and the translations $\R^{2}$ act naturally on geometric Sturmian lattices.
We describe these actions in terms of transformations of the defining sequences.

\begin{rem}
\label{rem:rotate-geom}
For a geometric Sturmian lattice 
$(a(i), b(j), c(k))$,
\[
\sigma(a(i), b(j), c(k)) =
(-b(-i), -c(-j), -a(-k))
\]
is the image of $(a(i), b(j), c(k))$ rotated $\pi/3$.
\[
\tau(a(i), b(j), c(k)) =
(-a(-i), -c(-j), -b(-k))
\]
is the mirror image with respect to the horizontal line $a = 0$.
Any translation can be described by
\begin{align*}
a(i)&\mapsto a(i) + \eta_{0},&
b(j)&\mapsto b(j) + \eta_{1},&
c(k)&\mapsto c(k) + \eta_{2}
\end{align*}
for some $\bm{\eta} = (\eta_{0}, \eta_{1}, \eta_{2})\in \R^{3}$ with $\eta_{0} + \eta_{1} + \eta_{2} = 0$.
This moves the origin $\Triplet{0}{0}{0}$ to $\Triplet{\eta_{0}}{\eta_{1}}{\eta_{2}}\in \R^{2}$.
\end{rem}

The following notion of shift-equivalence is also useful when we classify Sturmian lattices.

\begin{rem}
\label{rem:shift}
We identify two Sturmian lattices
\[
(a(i), b(j), c(k))\sfteq
(a(i_{0} + i), b(j_{0} + j), c(k_{0} + k))
\]
if $i_{0} + j_{0} + k_{0} = 0$.
They are said to be \emph{shift-equivalent}.
Two shift-equivalent Sturmian lattices give the same pattern $\mathfrak{A}\cup \mathfrak{B}\cup \mathfrak{C}$.
\end{rem}

For a geometric Sturmian lattice $(a(i), b(j), c(k))$, we call the unbounded strip
\[
\{a(i)\le a\le a(i+1)\}
\]
the \emph{corridor} of \emph{width} $a(i+1) - a(i)$, and similarly for the other two directions.
For a real number $\eta\ge 0$ and a triple $(a(i), b(j), c(k))$ defining a geometric Sturmian lattice, the change of corridor widths
\begin{equation}
\label{eq:comb-equiv1}
(a'(i), b'(j), c'(k)) =
(a(i) + i\eta, b(j) + j\eta, c(k) + k\eta)
\end{equation}
again defines a geometric Sturmian lattice. In particular, it satisfies
\begin{equation}
\label{eq:comb-equiv2}
i + j + k = 0\implies
a(i) + b(j) + c(k) = a'(i) + b'(j) + c'(k).
\end{equation}
We call the equivalence relation on geometric Sturmian lattices determined by \eqref{eq:comb-equiv2} \emph{combinatorial equivalence}. Conversely:
\begin{prop}
\label{prop:comb-equiv}
Suppose that two geometric Sturmian lattices $(a(i), b(j), c(k))$ and $(a'(i), b'(j), c'(k))$ are combinatorially equivalent, that is, satisfy \eqref{eq:comb-equiv2}. Then there exists $\eta\in \R$ such that the two lattices satisfy \eqref{eq:comb-equiv1} up to translation.
\end{prop}

\begin{proof}
Without loss of generality with respect to translation, we may assume $b(0) = b'(0)$ and $c(0) = c'(0)$. Combinatorial equivalence then forces $a(0) = a'(0)$. Put
\[
\eta:= a'(1) - a(1)\in \R.
\]
By combinatorial equivalence,
\[
a(1) + b(-1) + c(0) =
a'(1) + b'(-1) + c(0),
\]
so that $b'(-1) = b(-1) - \eta$.
Likewise $c'(-1) = c(-1) - \eta$.
Applying combinatorial equivalence once more,
\[
a(2) + b(-1) + c(-1) =
a'(2) + b'(-1) + c'(-1),
\]
that is, $a'(2) = a(2) + 2\eta$.
Repeating this argument yields the claim.
\end{proof}

\subsection{Axiomatic Sturmian lattices}
\label{sec:axiom}

Under combinatorial equivalence, condition \eqref{eq:def:geom} of Definition~\ref{def:geom} no longer carries essential information. This leads us to the following axiomatic notion of a Sturmian lattice.

\begin{definition}
\label{def:SL}
A triple $(a(i), b(j), c(k))$ is a \emph{Sturmian lattice} if, for every $(i, j, k)\in \Z^{3}$,
\begin{equation}
\label{eq:def:axiom}
i + j + k = 0\implies
\abs{a(i) + b(j) + c(k)} = \frac{1}{2}.
\end{equation}
\end{definition}

Our immediate goal is to determine and classify the global structure of a Sturmian lattice from this local information. 
Let us preview the arguments that follow in the simplest case, that of a triangular lattice.

\begin{prop}
\label{prop:trigonal}
Suppose a triple $(a(i), b(j), c(k))$ satisfies
\[
i + j + k = 0\implies
\abs{a(i) + b(j) + c(k)} = 0.
\]
Then there exist $\kappa\in \R$ and $(\eta_{0}, \eta_{1}, \eta_{2})\in \R^{3}$ with $\eta_{0} + \eta_{1} + \eta_{2} = 0$ such that
\begin{align*}
a(i)&= i\kappa + \eta_{0},&
b(j)&= j\kappa + \eta_{1},&
c(k)&= k\kappa + \eta_{2}.
\end{align*}
\end{prop}

This asserts that the local information ``every three concurrent lines meet at a single point'' already forces the global structure of a triangular lattice. It is noteworthy that each sequence then has a single increment,
\[
a(i+1) - a(i) =
b(j+1) - b(j) =
c(k+1) - c(k) = \kappa.
\]

\begin{proof}
Set
\begin{align*}
a(0)&= \eta_{0},&
b(0)&= \eta_{1},&
c(0)&= \eta_{2},
\end{align*}
which satisfy $\eta_{0} + \eta_{1} + \eta_{2} = 0$ by assumption.
For any $i, j\in \Z$, choosing $k\in \Z$ with $i + j + k + 1 = 0$ gives
\[
a(i+1) + b(j) + c(k) = 0 =
a(i) + b(j+1) + c(k),
\]
hence
\[
a(i+1) - a(i) = b(j+1) - b(j).
\]
Therefore there exists $\kappa\in \R$ such that, for all $i, j\in \Z$,
\begin{align*}
a(i+1) - a(i)&= \kappa,&
b(j+1) - b(j)&= \kappa.
\end{align*}
Similarly $c(k+1) - c(k) = \kappa$, which yields the claim.
\end{proof}

For $C\in \mathbb{N}\cup \{\infty\}$, a sequence $(a(i))$ is said to be \emph{$C$-color} if
\[
\Card(\{a(i+1) - a(i)\in \R \mid i\in \Z\})\le C.
\]
We also use these terms for a triple; $(a(i), b(j), c(k))$ is \emph{$C$-color} if
\[
\Card(\{a(n+1) - a(n), b(n+1) - b(n), c(n+1) - c(n)\in \R \mid n\in \Z\})\le C.
\]
A sequence $(a(i))$ (or a triple $(a(i), b(j), c(k))$) is said to be \emph{strictly $C$-color} if it is $C$-color but not $(C-1)$-color.

A Sturmian lattice need not be equidistant, but its corridor widths (so called even though they may take negative values) take only finitely many values.
This fact suggests that Sturmian lattices work well with symbolic sequences as we will discuss in \S\ref{sec:symb}.

\begin{prop}
\label{prop:at-most-3}
Any Sturmian lattice $(a(i), b(j), c(k))$ has at most $3$ colors.
\end{prop}

To prove the proposition, we use four lemmas.
Lemmas~\ref{lem:mutually-bal}~and~\ref{lem:2-bal} lead some ``balanced'' properties of the corridor colors.
Lemmas~\ref{lem:3-col-case}~and~\ref{lem:2-col-case} form the main part of the proof of Proposition~\ref{prop:at-most-3}.

\begin{lem}
\label{lem:mutually-bal}
Let $(a(i), b(j), c(k))$ be a Sturmian lattice.
Then
\begin{align*}
\bigl|(a(i+1) - a(i)) - (b(j+1) - b(j))\bigr|&\in \{0, \pm1\},\\
\bigl|(b(j+1) - b(j)) - (c(k+1) - c(k))\bigr|&\in \{0, \pm1\},\\
\bigl|(c(k+1) - c(k)) - (a(i+1) - a(i))\bigr|&\in \{0, \pm1\}
\end{align*}
hold for any $i, j, k\in \Z$.
\end{lem}

\begin{proof}
We only prove for $(a(i))$ and $(b(j))$.
Take $k\in \Z$ with $i + j + k + 1 = 0$.
By applying \eqref{eq:def:axiom} to $(i+1, j, k), (i, j+1, k)\in \Z^{3}$, we have
\begin{align*}
a(i+1) + b(j) + c(k)&\in \left\{\frac{1}{2}, -\frac{1}{2}\right\},&
a(i) + b(j+1) + c(k)&\in \left\{\frac{1}{2}, -\frac{1}{2}\right\}.
\end{align*}
Subtraction yields the desired result.
The other cases follow similarly.
\end{proof}

\begin{lem}
\label{lem:2-bal}
Let $(a(i), b(j), c(k))$ be a Sturmian lattice.
Then
\begin{align*}
\bigl|(a(i+1) - a(i)) - (a(i'+1) - a(i'))\bigr|&\in \{0, \pm1, \pm2\},\\
\bigl|(b(j+1) - b(j)) - (b(j'+1) - b(j'))\bigr|&\in \{0, \pm1, \pm2\},\\
\bigl|(c(k+1) - c(k)) - (c(k'+1) - c(k'))\bigr|&\in \{0, \pm1, \pm2\}
\end{align*}
hold for any $i, i', j, j', k, k'\in \Z$.
\end{lem}

\begin{proof}
The lemma immediately follows from Lemma~\ref{lem:mutually-bal}.
Indeed, for $(a(i))$, fix $j\in\Z$ and use the triangle inequality for
\begin{align*}
\abs{(a(i+1) - a(i)) - (b(j+1) - b(j))}\le 1,\\
\abs{(a(i'+1) - a(i')) - (b(j+1) - b(j))}\le 1,
\end{align*}
to obtain $\abs{(a(i+1) - a(i)) - (a(i'+1) - a(i'))}\le 2$.
Since the left-hand side is an integer, it belongs to $\{0,\pm1,\pm2\}$.
The other cases follow similarly.
\end{proof}

\begin{lem}
\label{lem:3-col-case}
Let $(a(i), b(j), c(k))$ be a Sturmian lattice and suppose that there are some $i', i''\in \Z$ with
\[
\bigl(a(i'+1) - a(i')\bigr) - \bigl(a(i''+1) - a(i'')\bigr) = 2.
\]
Then there exists $\kappa\in \R$ such that
\begin{align*}
a(i+1) - a(i)&\in \{\kappa - 1, \kappa, \kappa + 1\},&
b(j+1) - b(j)&= \kappa,&
c(k+1) - c(k)&= \kappa
\end{align*}
hold for any $i, j, k\in \Z$.
In particular, $(a(i), b(j), c(k))$ is strictly $3$-color under the assumption.
\end{lem}

\begin{proof}
Under our assumption, take $\kappa\in \R$ with
\begin{align*}
a(i'+1) - a(i')&= \kappa + 1,&
a(i''+1) - a(i'')&= \kappa - 1.
\end{align*}
By Lemma~\ref{lem:mutually-bal} and the former equation, we have
\[
\bigl(b(j+1) - b(j)\bigr) - \bigl(a(i'+1) - a(i')\bigr)\in \{-1, 0, 1\}
\]
for any $j\in \Z$ so that
\[
b(j+1) - b(j)\in \{\kappa, \kappa + 1, \kappa + 2\}.
\]
Similarly, the latter $a(i''+1) - a(i'') = \kappa - 1$ implies that
\[
b(j+1) - b(j)\in \{\kappa - 2, \kappa - 1, \kappa\}.
\]
Thus $b(j+1) - b(j) = \kappa$ for $j\in \Z$.
Similarly, $c(k+1) - c(k) = \kappa$ for $k\in \Z$.
\end{proof}

\begin{figure}[htb]
\includegraphics[width = \linewidth]
{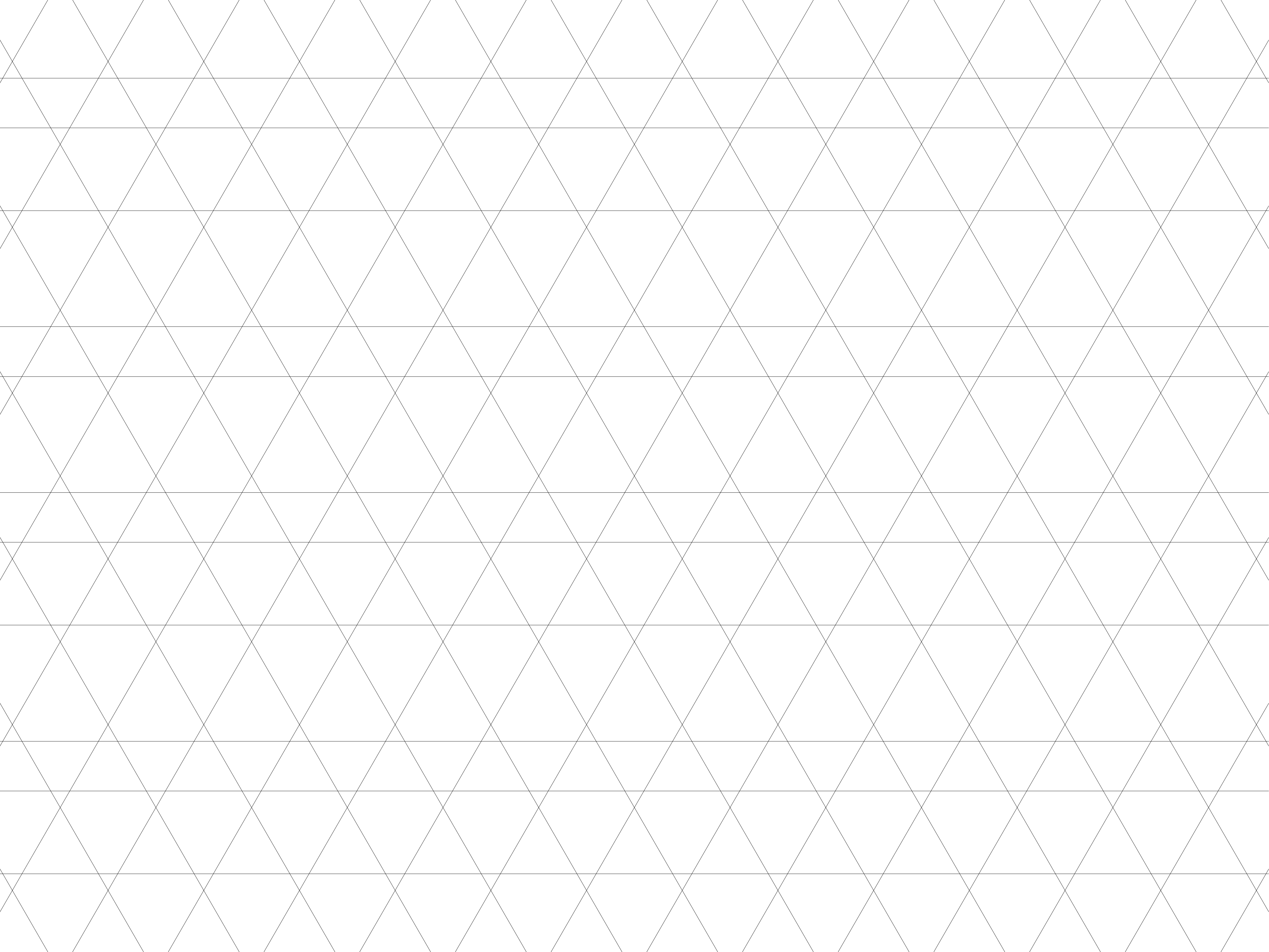}
\caption{Example of a $3$-color Sturmian lattice}
\label{fig:misfit}
\end{figure}

Figure~\ref{fig:misfit} shows an example of a strictly $3$-color geometric Sturmian lattice.
Here, the widths of the oblique corridors are all equal to $\kappa$, while the horizontal corridor widths take the values $\kappa$, $\kappa + 1$, and $\kappa - 1$.

\begin{lem}
\label{lem:2-col-case}
Let $(a(i), b(j), c(k))$ be a Sturmian lattice and suppose that
\begin{align*}
\abs{(a(i+1) - a(i)) - (a(i'+1) - a(i'))}&\le 1,\\
\abs{(b(j+1) - b(j)) - (b(j'+1) - b(j'))}&\le 1,\\
\abs{(c(k+1) - c(k)) - (c(k'+1) - c(k'))}&\le 1
\end{align*}
hold for any $i, i', j, j', k, k'\in \Z$.
Then there exists $\kappa\in \R$ such that
\[
a(i+1) - a(i), b(j+1) - b(j), c(k+1) - c(k)\in \{\kappa, \kappa + 1\}
\]
hold for any $i, j, k\in \Z$.
In particular, $(a(i), b(j), c(k))$ is $2$-color.
\end{lem}

\begin{proof}
Let
\[
\kappa = \inf_{i,j,k\in \Z}\bigl\{a(i+1) - a(i), b(j+1) - b(j), c(k+1) - c(k)\bigr\}.
\]
Since every width may take a value in discrete set $\kappa + \Z_{\ge0}$ by Lemmas~\ref{lem:mutually-bal}~and~\ref{lem:2-bal}, we may take $i\in \Z$ with $a(i+1) - a(i) = \kappa$ without loss of generality.
If there exists $i'\in \Z$ that $a(i'+1) - a(i')\ge \kappa + 2$, then it contradicts our assumption.
If $b(j+1) - b(j)\ge \kappa + 2$ or $c(k+1) - c(k)\ge \kappa + 2$ happen, then it violates Lemma~\ref{lem:mutually-bal}.
\end{proof}

\begin{proof}[Proof of Proposition~\ref{prop:at-most-3}]
The statement immediately follows from Lemmas~\ref{lem:3-col-case}~and~\ref{lem:2-col-case}.
If $\bigl(a(i'+1) - a(i')\bigr) - \bigl(a(i''+1) - a(i'')\bigr) = 2$ happens (or in other direction), then $(a(i), b(j), c(k))$ is strictly $3$-color by Lemma~\ref{lem:3-col-case}; otherwise, we may apply Lemma~\ref{lem:2-col-case} so that it is $2$-color.
\end{proof}

\subsection{Symbolic Sturmian lattices}
\label{sec:symb}

By Proposition~\ref{prop:at-most-3}, any Sturmian lattice is $3$-color.
In this subsection, we only treat $2$-color Sturmian lattices.
For a $2$-color Sturmian lattice $(a(i), b(j), c(k))$, we define
\[
\kappa:= \min_{i,j,k\in \Z}\bigl\{a(i+1) - a(i), b(j+1) - b(j), c(k+1) - c(k)\bigr\}\in \R,
\]
called the \emph{passage} of $(a(i), b(j), c(k))$.
We wish to forget the geometric information of Sturmian lattices, and study their combinatorial properties.
We extend the notation of heights of factors: for $w = (w_{n})\in \A^{\Z}$ and integers $n, n'\in \Z$, we define
\[
\abs{w_{[n,n')}}_{1}:= \begin{cases}
+\abs{w_{[n,n')}}_{1} & \text{if }n < n'\\
0 & \text{if }n = n'\\
-\abs{w_{[n',n)}}_{1} & \text{if }n > n'
\end{cases}.
\]
Note that the extended height still satisfies the additivity:
\[
\abs{w_{[n,n')}}_{1} +
\abs{w_{[n',n'')}}_{1} =
\abs{w_{[n,n'')}}_{1}
\]
for any $n, n', n''\in \Z$.

\begin{definition}
\label{def:symb}
A tuple $(a_{i}, b_{j}, c_{k})_{i,j,k}$ of bi-infinite words $(a_{i}), (b_{j}), (c_{k})\in \A^{\Z}$ is a \emph{symbolic Sturmian lattice} if
\begin{equation}
\label{eq:prop:axiom}
\abs{a_{[i,i')}}_{1} + \abs{b_{[j,j')}}_{1} + \abs{c_{[k,k')}}_{1}\in \{-1, 0, 1\}
\end{equation}
holds for any $(i, j, k), (i', j', k')\in \Z^{3}$ with $i + j + k = i' + j' + k' = 0$.
\end{definition}

Here we identify the two monochromatic Sturmian lattices $(0^{\Z}, 0^{\Z}, 0^{\Z})$ and $(1^{\Z}, 1^{\Z}, 1^{\Z})$.
For a $2$-color 
Sturmian lattice $(a(i), b(j), c(k))$ with passage $\kappa$, we define $\pi(a(i), b(j), c(k)) = (a_{i}, b_{j}, c_{k})$ by
\begin{equation}
\label{eq:proj}
\begin{aligned}
a_{i}&:= \begin{cases}
0 & \text{if }a(i+1) - a(i) = \kappa\\
1 & \text{if }a(i+1) - a(i) = \kappa + 1
\end{cases},\\
b_{j}&:= \begin{cases}
0 & \text{if }b(j+1) - b(j) = \kappa\\
1 & \text{if }b(j+1) - b(j) = \kappa + 1
\end{cases},\\
c_{k}&:= \begin{cases}
0 & \text{if }c(k+1) - c(k) = \kappa\\
1 & \text{if }c(k+1) - c(k) = \kappa + 1
\end{cases}.
\end{aligned}
\end{equation}
Then $(a_{i}, b_{j}, c_{k})$ is a symbolic Sturmian lattice.
In Proposition~\ref{prop:proj}, we shall show that the $2$-color Sturmian lattices and the symbolic Sturmian lattices are in one-to-one correspondence, apart from the trivial exception identified above.
Note that there is a relationship between the height of a factor in one direction and the widths of the corresponding corridors:
\begin{equation}
\label{eq:height}
\abs{a_{[i,i+n)}}_{1} = \bigl(a(i+n) - a(i)\bigr) - n\kappa,
\end{equation}
and similarly for the other directions.

Two words $x, y\in \A^{*}\cup \A^{\Z}$ are said to be \emph{mutually balanced} if
\[
-1\le |u|_{1} - |v|_{1}\le 1
\]
holds for any $u\prec x$, $v\prec y$ with $|u| = |v|$.
We see that if two $C$-balanced words $x, y\in \A^{\Z}$ are mutually balanced, then each slope of $x$ and $y$ coincide.

Propositions~\ref{prop:mutually-bal}~and~\ref{prop:2-bal} are ``balanced'' properties rephrased in terms of symbolic sequences.

\begin{prop}
\label{prop:mutually-bal}
Let $(a_{i}, b_{j}, c_{k})$ be a symbolic Sturmian lattice.
Then any two of the bi-infinite words $(a_{i}), (b_{j}), (c_{k})\in \A^{\Z}$ are mutually balanced.
\end{prop}

\begin{proof}
The proof is essentially the same as that of Lemma~\ref{lem:mutually-bal}.
For example, to see the mutually balancedness of $(a_{i})$ and $(b_{j})$, we may apply \eqref{eq:prop:axiom} to $(i, j+n, k), (i+n, j, k)\in V$ so that
\[
\abs{a_{[i,i+n)}}_{1} + \abs{b_{[j+n,j)}}_{1} + \abs{c_{[k,k)}}_{1} = \abs{a_{[i,i+n)}}_{1} - \abs{b_{[j,j+n)}}_{1}
\]
lies in $\{-1, 0, 1\}$.
Since $i$, $j$, and $n$ are arbitrary, we finished the proof.
\end{proof}

\begin{prop}
\label{prop:2-bal}
Let $(a_{i}, b_{j}, c_{k})$ be a symbolic Sturmian lattice.
Then each of the bi-infinite words $(a_{i}), (b_{j}), (c_{k})\in \A^{\Z}$ is $2$-balanced.
\end{prop}

\begin{proof}
We may use Proposition~\ref{prop:mutually-bal} and follow the discussion for Lemma~\ref{lem:2-bal}.
\end{proof}

We introduce the slope of a Sturmian lattice with $2$ colors.

\begin{prop}
\label{prop:slope}
Let $(a_{i}, b_{j}, c_{k})$ be a symbolic Sturmian lattice.
Then the limits
\begin{align*}
\alpha_{a}&= \lim_{N\to \infty}
\frac{1}{2N}\abs{a_{[-N,N)}}_{1},&
\alpha_{b}&= \lim_{N\to \infty}
\frac{1}{2N}\abs{b_{[-N,N)}}_{1},&
\alpha_{c}&= \lim_{N\to \infty}
\frac{1}{2N}\abs{c_{[-N,N)}}_{1}
\end{align*}
exist, and all three converge to the same value $\alpha = \alpha_{a} = \alpha_{b} = \alpha_{c}$.
\end{prop}

\begin{proof}
The existence of each limit follows from Proposition~\ref{prop:2-bal}, and the mutual balancedness (Proposition~\ref{prop:mutually-bal}) forces the limits to coincide.
\end{proof}

We call the common limit $\alpha\in [0, 1]$ defined here the \emph{slope} of a symbolic Sturmian lattice $(a_{i}, b_{j}, c_{k})$.
We denote by $\SL(\alpha)$ the class of Sturmian lattices with slope $\alpha$.
A $2$-color 
Sturmian lattice $(a(i), b(j), c(k))$ has \emph{slope} $\alpha$ if $\pi(a(i), b(j), c(k))$ has slope $\alpha$.
We denote by $\SL(\kappa, \alpha)$ the class of Sturmian lattices with passage $\kappa$ and slope $\alpha$.
By definition, $\alpha$ coincides with the natural density
\[
\lim_{N\to \infty}\frac{1}{2N}
\Card\{i\in [-N, N)\mid a(i+1) - a(i) = \kappa + 1\}
\]
of the wider corridors in each direction, and $1 - \alpha$ is of the narrower corridors.
Returning to the geometric Sturmian lattice, that is, to the case $\kappa\ge 1$, we can then compute the average of all corridor widths as
\[
(1 - \alpha)\kappa + \alpha(\kappa + 1) = \kappa + \alpha.
\]
We call the inverse $q:= 1 / (\kappa + \alpha)\in (0, 1]$ the \emph{frequency} of $\SL(\kappa, \alpha)$, which is the average number of straight lines that appear per length one.

\begin{cor}
\label{cor:Turtle-slope}
For any tiling by Smith Turtle~\cite{SMKGS:23_1}, its Ammann bars form a geometric Sturmian lattice of slope $(-1 + \sqrt{5}) / 2$ or $(3 - \sqrt{5}) / 2$.
\end{cor}

\begin{proof}
In \cite{Akiyama-Araki:23}, it is shown that its Ammann bars in Figure~\ref{fig:intro} form a Sturmian lattice in our definition.
Further it is also proved that the possible pair $(q, \kappa)$ of frequency and passage is one of $\bigl((5 - \sqrt{5}) / 10, 3\bigr)$ and $\bigl((5 + \sqrt{5}) / 10, 1\bigr)$.
By calculating $\alpha = 1/q - \kappa$ in each case, we obtain the desired values.
\end{proof}

The following proposition is a rephrase of Lemma~\ref{lem:3-col-case}.

\begin{prop}
\label{prop:2-bal-case}
Let $(a_{i}, b_{j}, c_{k})$ be a symbolic Sturmian lattice.
If $(a_{i})$ is strictly $2$-balanced, then $(b_{j})$ and $(c_{k})$ in the other two directions must be $1$-balanced periodic words.
\end{prop}

\begin{proof}
The discussion is essentially the same as Lemma~\ref{lem:3-col-case}.
Suppose that
\[
\abs{a_{[i,i+n)}}_{1} - \abs{a_{[i',i'+n)}}_{1} = 2
\]
holds for some $n\in \N$ and $i, i'\in \Z$ (here we take the minimum $n$).
Then we have for every $j\in \Z$
\[
\abs{b_{[j,j+n)}}_{1} = \abs{a_{[i,i+n)}}_{1} - 1 = \abs{a_{[i',i'+n)}}_{1} + 1 = \text{constant}
\]
by Proposition~\ref{prop:mutually-bal}.
Thus, $(b_{j})$ must be periodic.

To see the $1$-balancedness, suppose that $(b_{j})$ is strictly $2$-balanced.
Take the minimum $m\in \N$ that satisfies
\[
\abs{b_{[j,j+m)}}_{1} - \abs{b_{[j',j'+m)}}_{1} = 2
\]
for some $j, j'\in \Z$.
Note that $m < n$ holds since $(b_{j})$ has period $n$.
Applying the discussion above, one can see that $(a_{i})$ has period $m$.
However, it deduces that
\[
\abs{a_{[i,i+n-m)}}_{1} - \abs{a_{[i',i'+n-m)}} = 2,
\]
which violates the minimality of $n$.
\end{proof}

Finally, we restate the remarks on isometries and shifts in the symbolic setting.

\begin{rem}
\label{rem:rotate-symb}
The first remark is a symbolic version of Remark~\ref{rem:rotate-geom}.
For a $2$-color symbolic Sturmian lattice $(a_{i}, b_{j}, c_{k})$,
\[
\sigma(a_{i}, b_{j}, c_{k}) =
(b_{-i-1}, c_{-j-1}, a_{-k-1}) =
(\tilde{b}_{i}, \tilde{c}_{j}, \tilde{a}_{k})
\]
given by mirror words, is the image of $(a_{i}, b_{j}, c_{k})$ rotated $\pi/3$.
\[
\tau(a_{i}, b_{j}, c_{k}) =
(a_{-i-1}, b_{-j-1}, c_{-k-1}) =
(\tilde{a}_{i}, \tilde{b}_{j}, \tilde{c}_{k})
\]
describes the mirror image with respect to direction $a$.

The second remark is about shift-equivalence, relating to Remark~\ref{rem:shift}.
We identify two Sturmian lattices
\[
(a_{i}, b_{j}, c_{k})\sfteq (a_{i_{0}+i}, b_{j_{0}+j}, c_{k_{0}+k})
\]
if $i_{0} + j_{0} + k_{0} = 0$.
They are also said to be \emph{shift-equivalent}.
\end{rem}

\section{Classification of Sturmian lattices}
\label{sec:main}

\begin{figure}[htb]\centering
\includegraphics[width = \linewidth]
{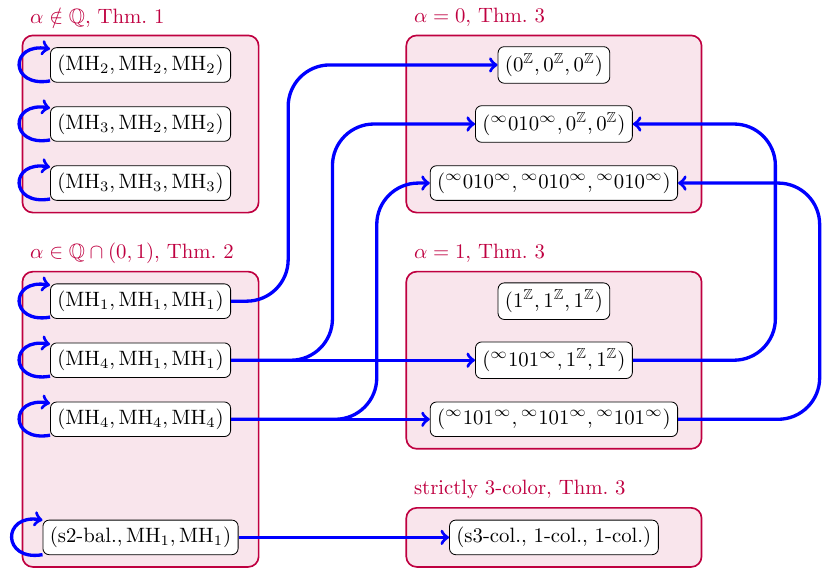}
\caption{Transition graph for the SL substitution $\Psi$}
\label{fig:transition}
\end{figure}

We classify all Sturmian lattices.
Figure~\ref{fig:transition} briefly summarizes the result.
For instance, $(\MH{3}, \MH{2}, \MH{2})$ is the class that one of the three sequences (say, $(a_{i})\in \A^{\Z}$) is of $(\MH{3})$ in the sense of Markoff~\cite{Markoff:79}, and that the other two sequences (say, $(b_{j})$ and $(c_{k})$) are of $(\MH{2})$.
Also, we say $(a_{i})$ is of (s2-bal.) if $(a_{i})$ is strictly $2$-balanced, and we say $(a_{i})$ is of (s3-col.) if $(a(i))$ is strictly $3$-color.
The blue arrows indicate the transitions by the SL substitution $\Psi$, which will be introduced in \S\ref{sec:negative}.

Our discussion can be divided into three parts:
\begin{itemize}
\item
In \S\ref{sec:main:irrational}, we deal with irrational slopes.
\item
In \S\ref{sec:main:rational}, we treat rational slopes but for $0$ and $1$.
\item
In \S\ref{sec:main:others}, we discuss the remaining cases: $3$-color, and $2$-color with slope $\alpha\in \{0, 1\}$.
\end{itemize}

In Theorem~\ref{thm:rational}, we describe the statement in terms of symbolic Sturmian lattices.
Through the following proposition, this simultaneously completes the classification of geometric Sturmian lattices as well.

\begin{prop}
\label{prop:proj}
For any symbolic Sturmian lattice $(a_{i}, b_{j}, c_{k})$, there exists a $2$-color Sturmian lattice $(a(i), b(j), c(k))$ such that
\[
\pi(a(i), b(j), c(k)) =
(a_{i}, b_{j}, c_{k}).
\]
Moreover, if $(a_{i}, b_{j}, c_{k})$ is strictly $2$-color, then such a $(a(i), b(j), c(k))$ is unique up to combinatorial equivalence.
\end{prop}

\begin{proof}
If $(a_{i}, b_{j}, c_{k}) = (0^{\Z}, 0^{\Z}, 0^{\Z})$, then only the two lattices
\begin{align*}
a(i)&= i\kappa\pm \frac{1}{2},&
b(j)&= j\kappa,&
c(k)&= k\kappa
\end{align*}
correspond to it. From now on, we assume otherwise.

Applying an isometry and a shift, we may assume without loss of generality that $b_{0} = 1$ and $c_{-1} = 0$. If a Sturmian lattice with $\pi(a(i), b(j), c(k)) = (a_{i}, b_{j}, c_{k})$ exists, then it must take the following form: for an arbitrary $\kappa\in \R$,
\begin{align*}
a(i)&= i\kappa + \abs{a_{[0,i)}}_{1} + a(0),&
b(j)&= j\kappa + \abs{b_{[0,j)}}_{1},&
c(k)&= k\kappa + \abs{c_{[0,k)}}_{1}
\end{align*}
with $a(0) = \pm1/2$.
Therefore, since
\begin{align*}
a(0) + b(0) + c(0)&= a(0),\\
a(0) + b(1) + c(-1)&=
a(0) + \abs{b_{0}}_{1} - \abs{c_{-1}}_{1} =
a(0) + 1,
\end{align*}
the axiom of a Sturmian lattice forces $a(0) = -1/2$.

It remains to verify that the $(a(i), b(j), c(k))$ determined above is indeed a Sturmian lattice.
In fact, suppose $i + j + k = 0$ and consider the single quantity
\begin{equation}
\label{eq:prf:two-windows}
\abs{a_{[0,i)}}_{1} + \abs{b_{[0,j)}}_{1} + \abs{c_{[0,k)}}_{1} =
\abs{a_{[0,i)}}_{1} + \abs{b_{[1,j)}}_{1} + \abs{c_{[-1,k)}}_{1} + 1,
\end{equation}
Here we used the equalities $b_{0} = 1$ and $c_{-1} = 0$.
The axiom \eqref{eq:prop:axiom} of a symbolic Sturmian lattice, applied to the starting index triples $(0, 0, 0)$ and $(0, 1, -1)$ (each ending at $(i, j, k)$), shows that the left-hand expression of \eqref{eq:prf:two-windows} lies in $\{-1, 0, 1\}$ while the right-hand expression lies in $\{0, 1, 2\}$.
Since the two are equal, their common value lies in the intersection $\{0, 1\}$.
Therefore
\[
a(i) + b(j) + c(k) =
\abs{a_{[0,i)}}_{1} + \abs{b_{[0,j)}}_{1} + \abs{c_{[0,k)}}_{1} - \frac{1}{2}
\in \left\{-\frac{1}{2}, +\frac{1}{2}\right\},\]
which yields the claim.
\end{proof}

The next result states a construction and rigidity of Sturmian lattices.

\begin{prop}
\label{prop:construction}
Let $(b_{j}), (c_{k})\in \A^{\Z}\setminus \{0^{\Z}\}$ be mutually balanced.
Then there exists $(a_{i})\in \A^{\Z}$ such that $(a_{i}, b_{j}, c_{k})$ is a symbolic Sturmian lattice.
Moreover, there are two choices
\begin{align*}                                                     
a&= \cdots a_{-2}0.1a_{1}\cdots,&
a'&= \cdots a_{-2}1.0a_{1}\cdots
\end{align*}
if and only if $b$ is the mirror word of $c$, i.e.,
\[
b = \cdots b_{-2}b_{-1}.b_{0}b_{1}\cdots =
\cdots c_{1}c_{0}.c_{-1}c_{-2}\cdots = \tilde{c}.
\]
\end{prop}

\begin{proof}
Suppose that $b$ and $c$ are mutually balanced.
We again read these as
\begin{align*}
b(j)&= j\kappa + \abs{b_{[0,j)}}_{1},&
c(k)&= k\kappa + \abs{c_{[0,k)}}_{1},
\end{align*}
and look for a suitable $a(0)\in \{\pm1/2\}$.

First, when $b = \tilde{c}$, that is, $b_{n-1} = c_{-n}$ holds for all $n\in \Z$, then
\[
b(n) + c(-n) =
\abs{b_{[0,n)}}_{1} - \abs{c_{[-n,0)}}_{1} =
0,
\]
so either choice $a(0) = \pm1/2$ is admissible.

If there exists $n\in \Z$ with $b_{n-1}\ne c_{-n}$, then, by a shift that leaves $(a_{i})$ unchanged, we may assume $b_{0}\ne c_{-1}$; and, by a rotation through $\pi$, we may further assume $b_{0} = 1\ne 0 = c_{-1}$.
In this case, the same argument as in Proposition~\ref{prop:proj} forces $a(0) = -1/2$.
With this choice, the same two-window estimate as in \eqref{eq:prf:two-windows}---here applying mutual balancedness to $b$ and $c$ instead of the axiom \eqref{eq:prop:axiom}---gives
\[
\abs{b_{[0,n)}}_{1} + \abs{c_{[0,-n)}}_{1} =
\abs{b_{[1,n)}}_{1} + \abs{c_{[-1,-n)}}_{1} + 1\in \{0, 1\},
\]
so that
\[
a(0) + b(n) + c(-n) =
\abs{b_{[0,n)}}_{1} + \abs{c_{[0,-n)}}_{1} - \frac{1}{2}\in
\left\{-\frac{1}{2}, +\frac{1}{2}\right\}.
\]
The other $a(i)$ are determined in the same way, and this $(a(i))$ is $2$-color by Lemma~\ref{lem:3-col-case}.
\end{proof}

\subsection{Irrational case}
\label{sec:main:irrational}
We begin with irrational slopes.
We fix $\alpha\in (0, 1)\setminus \Q$ as the slope.
Theorem~\ref{thm:irrational} states that any Sturmian lattice of an irrational slope can be described with Sturmian words.
For unified treatment of lower and upper words, we prepare three copies of $\R$:
\[
\R^{-}\sqcup \R\sqcup \R^{+}:=
\R\times \{{-}, \emptyset, {+}\} =
\{x^{-}\mid x\in \R\}\sqcup
\R\sqcup
\{x^{+}\mid x\in \R\}
\]
equipped with the lexicographical order, i.e.,
\[
(x - \varepsilon)^{+} <
x^{-} < x < x^{+} <
(x + \varepsilon)^{-}
\]
for $x\in \R$ and $\varepsilon > 0$.
We also use the additive action $\R\curvearrowright \R^{\pm}$:
\[
c + x^{\pm}:= (c+x)^{\pm}
\]
for $c\in \R$, and the multiplicative action $\R\setminus \{0\}\curvearrowright \R^{\pm}$:
\begin{align*}
c\cdot x^{\pm}&:= (cx)^{\pm},&
-c\cdot x^{\pm}&:= (-cx)^{\mp}
\end{align*}
for $c > 0$.
We denote by
\[
\st(x^{-}) = \st(x) = \st(x^{+}):= x\in \R
\]
the ``standard part function'' $\st\colon \R^{-}\sqcup \R\sqcup \R^{+}\to \R$.
We introduce the rounding function $\round{\bullet}\colon \R^{-}\sqcup (\R\setminus \Z)\sqcup \R^{+}\to \Z + 1/2$ as follows:
\begin{align*}
\round{x^{+}}&:= \floor{x} + \frac{1}{2},&
\round{x^{-}}&:= \ceil{x} - \frac{1}{2}
\end{align*}
for $x\in \R$, and
\[
\round{x}:= \floor{x} + \frac{1}{2} =
\ceil{x} - \frac{1}{2}
\]
for $x\in \R\setminus \Z$.
This function is odd, i.e., $\round{-x} = -\round{x}$ for $x\notin \Z$.
We define
\begin{align*}
\R_{\alpha}&=
(\Z + \Z\alpha)^{-}\sqcup
(\Z + \Z\alpha)^{+}\sqcup
(\Z + \Z\alpha)^{c}\\&:=
\{\rho^{-}\in \R^{-}\mid \rho\in \Z + \Z\alpha\}\sqcup
\{\rho^{+}\in \R^{+}\mid \rho\in \Z + \Z\alpha\}\sqcup
\{\rho\in \R\mid \rho\notin \Z + \Z\alpha\}.
\end{align*}
There is a canonical map from $\R_{\alpha}$ to the set of Sturmian words of slope $\alpha$: for $\rho\in \R_{\alpha}$, the word $(w_{n})\in \A^{\Z}$ defined by
\[
w_{n} = \round{(n+1)\alpha + \rho} - \round{n\alpha + \rho}
\]
is the Sturmian word of slope $\alpha$ and intercept $\rho$.
If $\rho\notin (\Z + \Z\alpha)^{-}$, then $(w_{n})$ is lower mechanical.
If $\rho\notin (\Z + \Z\alpha)^{+}$, then $(w_{n})$ is upper mechanical.
In particular, $\rho\notin \Z + \Z\alpha$ if and only if $(w_{n})$ is type of $(\MH2)$.

For $\bm{\eta} = (\eta_{0}, \eta_{1}, \eta_{2})\in \R^{3}$ and $\bm{\rho} = (\rho_{0}, \rho_{1}, \rho_{2})\in \R_{\alpha}^{3}$, we denote by
\[
\SL(\kappa, \alpha\mid \bm{\eta}, \bm{\rho}) =
\SL\left(
\begin{array}{@{\,}c|ccc@{\,}}
\kappa & \eta_{0} & \eta_{1} & \eta_{2}\\
\alpha & \rho_{0} & \rho_{1} & \rho_{2}
\end{array}\right):=
(a(i), b(j), c(k))
\]
the triple of the following three bi-infinite $\R$-valued sequences
\begin{equation}
\label{eq:mechanical}
\begin{aligned}
a(i)&:=
i\kappa + \eta_{0} + \round{i\alpha + \rho_{0}},\\
b(j)&:=
j\kappa + \eta_{1} + \round{j\alpha + \rho_{1}},\\
c(k)&:=
k\kappa + \eta_{2} + \round{k\alpha + \rho_{2}}.
\end{aligned}
\end{equation}
We call $\bm{\eta}$ and $\bm{\rho}$ the \emph{basepoint} and \emph{intercept} of $\SL(\kappa, \alpha\mid \bm{\eta}, \bm{\rho})$, respectively.
Note that the associated symbolic Sturmian lattice $(a_{i}, b_{j}, c_{k}) = \pi(a(i), b(j), c(k))$ is given by
\begin{equation}
\label{eq:prf:round-symb}
\begin{aligned}
a_{i}&= \round{(i+1)\alpha + \rho_{0}} - \round{i\alpha + \rho_{0}},\\
b_{j}&= \round{(j+1)\alpha + \rho_{1}} - \round{j\alpha + \rho_{1}},\\
c_{k}&= \round{(k+1)\alpha + \rho_{2}} - \round{k\alpha + \rho_{2}},
\end{aligned}
\end{equation}
which are bi-infinite Sturmian words of slope $\alpha$.
We denote by
\[
\SL(\alpha\mid \bm{\rho}) =
\SL(\alpha\mid \rho_{0}, \rho_{1}, \rho_{2}):=
(a_{i}, b_{j}, c_{k}).
\]

\begin{thm}
\label{thm:irrational}
Let $(a(i), b(j), c(k))\in \SL(\kappa, \alpha)$ be a Sturmian lattice of irrational slope $\alpha\in (0, 1)$.
Then all three words $(a_{i}), (b_{j}), (c_{k})\in \A^{\Z}$ are Sturmian words of slope $\alpha$.
Moreover, there exist $\bm{\eta} = (\eta_{0}, \eta_{1}, \eta_{2})\in \R^{3}$, $\bm{\rho}:= (\rho_{0}, \rho_{1}, \rho_{2})\in \R_{\alpha}^{3}$ with
\begin{align*}
\eta_{0} + \eta_{1} + \eta_{2}&= 0,&
\st(\rho_{0}) + \st(\rho_{1}) + \st(\rho_{2})&= 0
\end{align*}
such that $(a(i), b(j), c(k)) = \SL(\kappa, \alpha\mid \bm{\eta}, \bm{\rho})$.
If $\rho_{0}, \rho_{1}, \rho_{2}\in (\Z + \Z\alpha)^{\pm}$, then there are six possibilities
\begin{align*}
&(\rho_{0}^{+}, \rho_{1}^{-}, \rho_{2}^{-}),&
&(\rho_{0}^{-}, \rho_{1}^{+}, \rho_{2}^{-}),&
&(\rho_{0}^{-}, \rho_{1}^{-}, \rho_{2}^{+}),\\
&(\rho_{0}^{-}, \rho_{1}^{+}, \rho_{2}^{+}),&
&(\rho_{0}^{+}, \rho_{1}^{-}, \rho_{2}^{+}),&
&(\rho_{0}^{+}, \rho_{1}^{+}, \rho_{2}^{-})
\end{align*}
to choose the lower or upper rounding scheme.
\end{thm}

The proof of Theorem~\ref{thm:irrational} will be completed at the end of this subsection.
If $\alpha$ is irrational, then all three words $(a_i)$, $(b_j)$, and $(c_k)$ are Sturmian.
Indeed, Proposition~\ref{prop:2-bal-case} shows that they are non-periodic and $1$-balanced.
Hence they are bi-infinite Sturmian words by the Morse--Hedlund theorem; see \cite{Lothaire:02, Fogg}.

Without loss of generality, we may take $\bm{\eta} = (\eta_{0},\eta_{1},\eta_{2})\in\R^{3}$ with $\eta_{0}+\eta_{1}+\eta_{2}=0$.
Indeed,
\[
a(0) + b(0) + c(0) =
(\eta_{0} + \eta_{1} + \eta_{2}) + \round{\rho_{0}} + \round{\rho_{1}} + \round{\rho_{2}}\in \{\pm1/2\}
\]
and $\round{\rho_{0}} + \round{\rho_{1}} + \round{\rho_{2}}\in \Z + 1/2$ imply that $\eta_{0} + \eta_{1} + \eta_{2}\in \Z$.
Moreover, since
\begin{equation}
\label{eq:identification}
\SL\left(
\begin{array}{@{\,}c|ccc@{\,}}
\kappa & \eta_{0} & \eta_{1} & \eta_{2}\\
\alpha & \rho_{0} + 1 & \rho_{1} & \rho_{2}
\end{array}\right) =
\SL\left(
\begin{array}{@{\,}c|ccc@{\,}}
\kappa & \eta_{0} + 1 & \eta_{1} & \eta_{2}\\
\alpha & \rho_{0} & \rho_{1} & \rho_{2}
\end{array}\right)
\end{equation}
holds, we may take $\bm{\eta}\in \R^{3}$ with $\eta_{0} + \eta_{1} + \eta_{2} = 0$.
For simplicity, recalling Remark~\ref{rem:rotate-geom}, we may set $\bm{\eta}$ to be $(\eta_{0}, \eta_{1}, \eta_{2}) = (0, 0, 0)$ by acting translation.

We claim that $\bm{\rho} = (\rho_{0}, \rho_{1}, \rho_{2})\in \R_{\alpha}^{3}$ yields $\st(\rho_{0}) + \st(\rho_{1}) + \st(\rho_{2}) = 0$.

\begin{lem}
\label{lem:rho}
Let $\rho_{0}, \rho_{1}, \rho_{2}\in \R$.
Then
\[
\round{-x - y + \rho_{0}} + \round{x + \rho_{1}} + \round{y + \rho_{2}}\in \left\{\frac{1}{2}, -\frac{1}{2}\right\}
\]
holds for almost all $(x, y)\in \R^{2}$ if and only if $\rho_{0} + \rho_{1} + \rho_{2} = 0$.
\end{lem}

The exceptional null set in Lemma~\ref{lem:rho} is $N = N_{0}\cup N_{1}\cup N_{2}$ where
\begin{align*}
N_{0}&= \{(x, y)\in \R^{2}\mid x + y - \rho_{0}\in \Z\},\\
N_{1}&= \{(x, y)\in \R^{2}\mid x + \rho_{1}\in \Z\},\\
N_{2}&= \{(x, y)\in \R^{2}\mid y + \rho_{2}\in \Z\}.
\end{align*}

\begin{figure}[htb]\centering
\subfigure[$\round{-x - y + \rho_{0}}$]{%
\includegraphics[height=5cm]{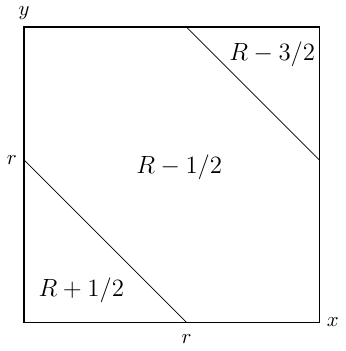}
}\qquad%
\subfigure[$\round{x + \rho_{1}} + \round{y + \rho_{2}}$]{%
\includegraphics[height=5cm]{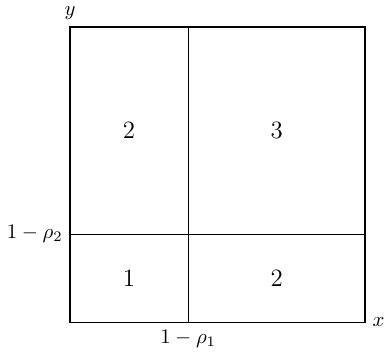}
}
\caption{Values of the rounding functions on $[0, 1)^{2}$}
\label{fig:round}
\end{figure}

\begin{proof}
We define
\[
f(x, y) = \round{-x - y + \rho_{0}} + \round{x + \rho_{1}} + \round{y + \rho_{2}}.
\]
Since $f(x, y)$ is periodic with respect to $x$ and $y$ of period one, we may assume that $x, y\in [0, 1)$.
Moreover, the set of values of $f$ is invariant under the transformation $(\rho_{0}, \rho_{1}, \rho_{2})\mapsto (\rho_{0}', \rho_{1}', \rho_{2}')$ with $\rho_{0} + \rho_{1} + \rho_{2} = \rho_{0}' + \rho_{1}' + \rho_{2}'$, and it suffices to consider $\rho_{1}, \rho_{2}\in [0, 1)$.
We have
\[
\round{-x - y + \rho_{0}} = \begin{cases}
R + 1/2 & \text{if }x + y < r\\
R - 1/2 & \text{if }r < x + y < 1 + r\\
R - 3/2 & \text{if }x + y > 1 + r
\end{cases}
\]
for $R\in \Z$ and $r\in [0, 1)$ with $R + r = \rho_{0}$, and
\[
\round{x + \rho_{1}} + \round{y + \rho_{2}} = \begin{cases}
1 & \text{if }x < 1 - \rho_{1}\text{ and }y < 1 - \rho_{2}\\
3 & \text{if }x > 1 - \rho_{1}\text{ and }y > 1 - \rho_{2}\\
2 & \text{otherwise}
\end{cases}.
\]
See Figure~\ref{fig:round}.
Since 
$\round{x+\rho_1} + \round{y+\rho_2}$ takes three values $\{1,2,3\}$ 
in any ball centered at $(1-\rho_1, 1-\rho_2)$, 
if the line $x+y=r$ does not pass through the point
$(1-\rho_1, 1-\rho_2)$, then $f(x,y)$ must take at least three values. 
On the other hand, if it passes through $(1-\rho_1,1-\rho_2)$, then the values of $f$ are in $\{R + 3/2, R + 5/2\}$.
Thus $\abs{f(x, y)} = 1/2$ holds for almost all $(x, y)\in [0, 1)^{2}$ if and only if both
\begin{align*}
(1 - \rho_{1}) + (1 - \rho_{2})&= r,&
R&= -2
\end{align*}
hold, and we have $\rho_{0} + \rho_{1} + \rho_{2} = 0$ as required.
\end{proof}

The following proposition gives a specific description of Sturmian lattices generated by three Sturmian words of type $(\MH2)$.

\begin{prop}
\label{prop:rho}
Let $\alpha$ be irrational and $\rho_{0}, \rho_{1}, \rho_{2}\notin \Z + \Z\alpha$.
Then $\SL(\kappa, \alpha\mid \bm{0}, \bm{\rho})$ forms a Sturmian lattice if and only if $\rho_{0} + \rho_{1} + \rho_{2} = 0$ holds.
\end{prop}

\begin{proof}
Suppose that $\rho_{0} + \rho_{1} + \rho_{2} = 0$ and let
\begin{align*}
a(i)&= i\kappa + \round{i\alpha + \rho_{0}},&
b(j)&= j\kappa + \round{j\alpha + \rho_{1}},&
c(k)&= k\kappa + \round{k\alpha + \rho_{2}}.
\end{align*}
It suffices to check \eqref{eq:def:axiom}, namely,
\[
\bigl|a(i) + b(j) + c(k)\bigr| = \frac{1}{2}
\]
for $i + j + k = 0$.
If $i + j + k = 0$, then Lemma~\ref{lem:rho} implies \eqref{eq:def:axiom} everywhere since
\[
\{(j\alpha, k\alpha)\in \R^{2}\mid j, k\in \Z\}\cap N = \emptyset.
\]
The converse is proved in the same way.
\end{proof}

\begin{figure}[htb]\centering
\subfigure[$-$, $+$, $+$]{%
\includegraphics[width = .2\linewidth]{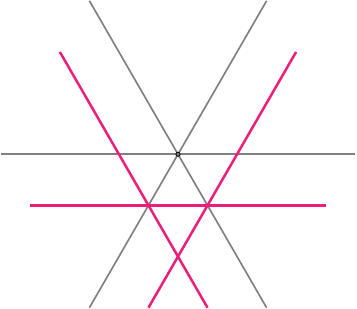}
}\quad%
\subfigure[$-$, $-$, $+$]{%
\includegraphics[width = .2\linewidth]{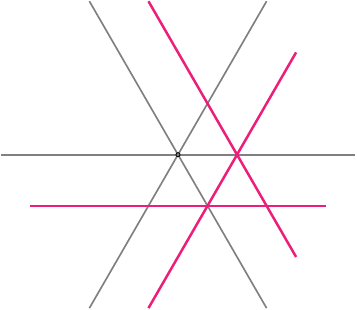}
}\quad%
\subfigure[$+$, $-$, $+$]{%
\includegraphics[width = .2\linewidth]{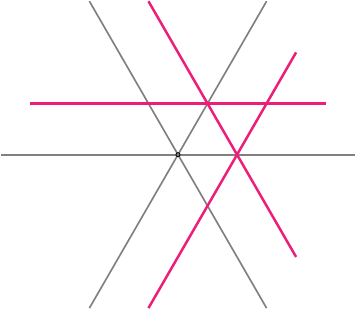}
}\quad%
\subfigure[$+$, $-$, $-$]{%
\includegraphics[width = .2\linewidth]{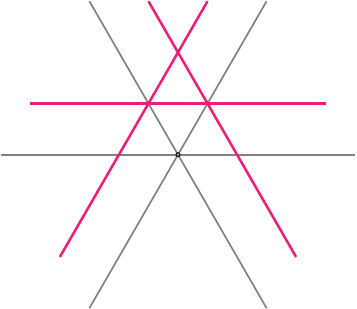}
}\\%
\subfigure[$+$, $+$, $-$]{%
\includegraphics[width = .2\linewidth]{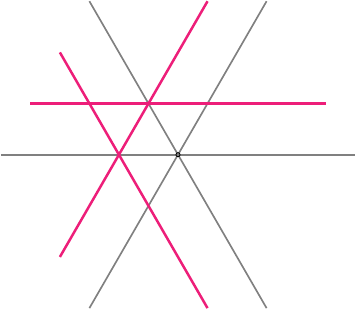}
}\quad%
\subfigure[$-$, $+$, $-$]{%
\includegraphics[width = .2\linewidth]{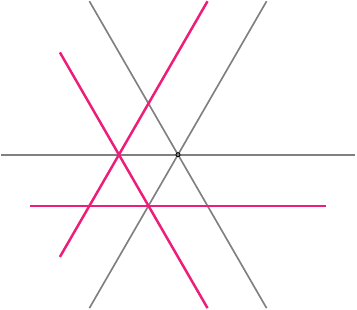}
}\quad%
\subfigure[$-$, $-$, $-$*]{%
\includegraphics[width = .2\linewidth]{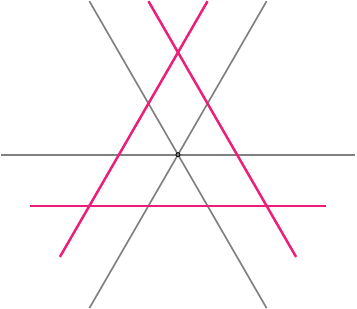}
}\quad%
\subfigure[$+$, $+$, $+$*]{%
\includegraphics[width = .2\linewidth]{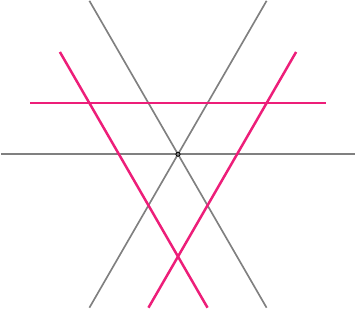}
}
\caption{The positions of lines $a = a(0)$, $b = b(0)$, and $c = c(0)$ (the symbols * indicate forbidden patterns)}
\label{fig:cf-choice}
\end{figure}

We next consider the case that at least one of the three Sturmian words is of type $(\MH3)$.
The situation is almost the same as before: $\st(\rho_{0}) + \st(\rho_{1}) + \st(\rho_{2}) = 0$ must be satisfied.
The remaining issue is the choice between the lower and upper forms, and this is resolved by Proposition~\ref{prop:construction}.
\begin{itemize}
\item
Let $\rho_{0}\in (\Z + \Z\alpha)^{\pm}$ but $\rho_{1}, \rho_{2}\notin \Z + \Z\alpha$.
Recalling \eqref{eq:identification} and
\[
\SL\left(
\begin{array}{@{\,}c|ccc@{\,}}
\kappa & \eta_{0} & \eta_{1} & \eta_{2}\\
\alpha & \rho_{0} + \alpha & \rho_{1} & \rho_{2}
\end{array}\right)\sfteq
\SL\left(
\begin{array}{@{\,}c|ccc@{\,}}
\kappa & \eta_{0} - \kappa & \eta_{1} & \eta_{2}\\
\alpha & \rho_{0} & \rho_{1} & \rho_{2}
\end{array}\right),
\]
we may assume that $\st(\rho_{0}) = 0$.
Then $\rho:= \rho_{1} = -\rho_{2}$ so that $b = \tilde{c}$.
By Proposition~\ref{prop:construction}, the choice at $a(0)$ occurs and both options are possible; we obtain two Sturmian lattices with respect to intercepts $(0^{+}, \rho, -\rho)$ and $(0^{-}, \rho, -\rho)$, which are the images of each other under reflection (see Remark~\ref{rem:rotate-symb}).
\item
Let $\rho_{0}, \rho_{1}\in (\Z + \Z\alpha)^{\pm}$ (and thus $\rho_{2}\in (\Z + \Z\alpha)^{\pm}$).
We may assume again that $\st(\rho_{0}) = \st(\rho_{1}) = \st(\rho_{2}) = 0$.
In this case, the equations
\begin{align*}
a(0^{+}) = b(0^{+}) = c(0^{+})&= +\frac{1}{2},&
a(0^{-}) = b(0^{-}) = c(0^{-})&= -\frac{1}{2}
\end{align*}
hold, so there are six admissible possibilities, shown as (a)--(f) in Figure~\ref{fig:cf-choice}, according to the choice of $0^{+}$ or $0^{-}$ in each direction.
For example, ``$-$, $+$, $+$'' in (a) means that we set $\bm{\rho}=(0^{-},0^{+},0^{+})$.
Two further cases, (g) and (h), are exceptional:
\[
\abs{a(0^{+}) + b(0^{+}) + c(0^{+})} =
\abs{a(0^{-}) + b(0^{-}) + c(0^{-})} = \frac{3}{2}.
\]
The six admissible Sturmian lattices form a single orbit under the rotation group $C_{6}$ of order six.
Here we take (a) for the initial state so that
\begin{align*}
a&= \tilde{u}0.1u,&
b&= \tilde{u}1.0u,&
c&= \tilde{u}1.0u,
\end{align*}
where $u = u_{1}u_{2}\cdots\in \A^{\N}$ is the \emph{left special sequence} of slope $\alpha$, that is, the right-infinite word $u_{n} = \round{(n+1)\alpha} - \round{n\alpha}$ ($n\ge 1$), characterized by the property that both $0u$ and $1u$ are Sturmian words of slope $\alpha$; see \cite[Chapter~6]{Fogg}.
Since $\tilde{a} = b = c$, there are two choices in the sense of Proposition~\ref{prop:construction}: switching at $b(0)$ or at $c(0)$.
If we switch $b$: $\tilde{u}1.0u\mapsto \tilde{u}0.1u$, then the Sturmian lattice changes to (b).
Similarly, one passes cyclically through the six cases
\[
\text{(a)}\leftrightarrow \text{(b)}\leftrightarrow \cdots \leftrightarrow \text{(f)}\leftrightarrow \text{(a)}.
\]
\end{itemize}

%
%

\subsection{Rational case}
\label{sec:main:rational}
We next discuss the classification of (symbolic) Sturmian lattices with rational slopes (except for $0$ and $1$).

\begin{thm}
\label{thm:rational}
Let $(a_{i}, b_{j}, c_{k})\in \SL(\alpha)$ be a symbolic Sturmian lattice of rational slope $\alpha\in (0, 1)$.
Then all three words $(a_{i}), (b_{j}), (c_{k})\in \A^{\Z}$ have the form
\begin{equation}
\label{eq:thm:rational}
(f_{n}) = \cdots cw_{-1}cw_{0}cw_{1}c\cdots
\end{equation}
for some $w_{s}\in \{01, 10\}$ and $w_{0} = f_{-1}.f_{0}$.
Here $c$ is the central word of slope $\alpha$.
More precisely, there are two cases (up to isometry and shift-equivalence):
\begin{itemize}
\item
$(b_{j})$ and $(c_{k})$ are periodic $1$-balanced words, and $(a_{i})$ may be $2$-balanced:
\begin{align*}
(a_{i})&= \cdots cw_{-1}cw_{0}cw_{1}c\cdots\\
(b_{j})&= {}^{\infty}(1c0).(1c0)^{\infty}\\
(c_{k})&= {}^{\infty}(0c1).(0c1)^{\infty}
\end{align*}
with any choice of $w_{s}\in \{01, 10\}$.
\item
$(a_i)$, $(b_j)$, and $(c_k)$ are the same skew $1$-balanced word:
\begin{align*}
&\left\{\begin{aligned}
(a_{i})&= {}^{\infty}(0c1).0c0(1c0)^{\infty}\\
(b_{j})&= {}^{\infty}(0c1)0c0.(1c0)^{\infty}\\
(c_{k})&= {}^{\infty}(0c1)0c0.(1c0)^{\infty}
\end{aligned}\right.,&
&\left\{\begin{aligned}
(a_{i})&= {}^{\infty}(1c0).1c1(0c1)^{\infty}\\
(b_{j})&= {}^{\infty}(1c0)1c1.(0c1)^{\infty}\\
(c_{k})&= {}^{\infty}(1c0)1c1.(0c1)^{\infty}
\end{aligned}\right..
\end{align*}
\end{itemize}
\end{thm}

\begin{proof}
First we verify that all three words must have the form \eqref{eq:thm:rational}.
It suffices to show the claim for a strictly $2$-balanced word $(a_{i})$.
By Proposition~\ref{prop:2-bal-case}, the other words $(b_{j})$ and $(c_{k})$ are both $1$-balanced periodic, i.e., of type $(\MH1)$ by Lemma~\ref{lem:Markoff}.
Thus these two words have period $|c| + 2$ and $(b_{j})$ is the mirror word of $(c_{k})$.
Applying Proposition~\ref{prop:construction}, we find periodic positions at which we may choose $a_{i-1}a_i=01$ or $10$, and one sees that $(a_i)$ has the form
\[
\cdots d_{-1}w_{-1}d_0w_0d_1w_1\cdots
\]
with arbitrary choices of $w_i\in\{01,10\}$ and some words $d_i$ of length $|c|$.
Recalling that $(a_{i})$ and $(b_{j}) = (0c1)^{\Z}$ are mutually balanced, and $(0c1)^{\Z}$ and $(1c0)^{\Z}$ are shift equivalent, four pairs of subwords
\begin{align*}
&(0d_{i}0, 0c1),&
&(0d_{i}0, 1c0),&
&(1d_{i}1, 0c1),&
&(1d_{i}1, 1c0)
\end{align*}
of $(a_{i})$ and $(b_{j})$ are also mutually balanced.
We claim that $d_{i} = c$.
Indeed, if $d_i\ne c$, then there exist words $w$, $u$, $v$ such that $d_i=w0u$ and $c=w1v$, or $d_i=w1u$ and $c=w0v$.
However, $1w1$ and $0w0$ are not mutually balanced, contradicting the mutual balancedness of $(a_i)$ and $(b_j)$.

It remains to classify the case where all three words are $1$-balanced.
By Lemma~\ref{lem:Markoff}, a periodic $1$-balanced word is of type $(\MH1)$ or $(\MH4)$.
The case where both $(b_j)$ and $(c_k)$ are of type $(\MH1)$ has already been treated; in this case there are infinitely many choices for $(a_i)$.
If $(b_j)$ is of type $(\MH1)$ but $(a_i)$ is not, then $(c_k)$ is uniquely determined by Proposition~\ref{prop:construction}, and it must also be of type $(\MH1)$.
The remaining case is that $(a_i)$, $(b_j)$, and $(c_k)$ are all of type $(\MH4)$.
Since the two skew $1$-balanced words
\begin{align*}
&{}^{\infty}(0c1)0c0(1c0)^{\infty},&
&{}^{\infty}(1c0)1c1(0c1)^{\infty}
\end{align*}
are not mutually balanced, the three words must be of the same type.
Recalling the discussion in the irrational case with $\rho_0=\rho_1=\rho_2$, there are six possibilities, forming a single orbit under rotation.
\end{proof}

\subsection{Other cases}
\label{sec:main:others}

We conclude this section by classifying $3$-color Sturmian lattices and $2$-color Sturmian lattices with slopes $0$ or $1$.
The result is analogous to Theorem~\ref{thm:rational}.
For a $3$-color Sturmian lattice as well, we assign the letters $\{\bar{1}, 0, 1\}$ to the corridor widths:
\[
a_{i}:= \begin{cases}
\bar{1} & \text{if }a(i+1) - a(i) = \kappa - 1,\\
0 & \text{if }a(i+1) - a(i) = \kappa,\\
1 & \text{if }a(i+1) - a(i) = \kappa + 1.
\end{cases}
\]
\color{black}

\begin{thm}
\label{thm:others}
Let $(a(i), b(j), c(k))$ be a Sturmian lattice whose slope is $0$, $1$, or undefined.
Then the coordinate sequences $(a(i))$, $(b(j))$, and $(c(k))$ have the form
\[
f(n) = n\kappa + \varepsilon(n)
\]
for some $\varepsilon\colon \Z\to \{-1/2, 1/2\}$.
More precisely, there are two cases (up to isometry and shift-equivalence):
\begin{itemize}
\item
$(b(j))$ and $(c(k))$ are equidistant, while $(a(i))$ may be $3$-color,
i.e., $(a_{i}, b_{j}, c_{k}) = \pi(a(i), b(j), c(k))$ is of
\begin{align*}
(a_{i})&= \cdots \bar{1}0^{u(-1)}10^{u(0)}\bar{1}0^{u(1)}1\cdots,&
(b_{j}) = (c_{k})&= 0^{\Z}
\end{align*}
with any choice of $u(n)\in \Z_{\ge0}$.
\item
$(a(i))$, $(b(j))$, and $(c(k))$ are equidistant except for one corridor in each direction,
i.e., $(a_{i}, b_{j}, c_{k}) = \pi(a(i), b(j), c(k))$ is one of the following:
\begin{align*}
&\left\{\begin{aligned}
(a_{i})&= {}^{\infty}1.01^{\infty}\\
(b_{j})&= {}^{\infty}10.1^{\infty}\\
(c_{k})&= {}^{\infty}10.1^{\infty}
\end{aligned}\right.,&
&\left\{\begin{aligned}
(a_{i})&= {}^{\infty}0.10^{\infty}\\
(b_{j})&= {}^{\infty}01.0^{\infty}\\
(c_{k})&= {}^{\infty}01.0^{\infty}
\end{aligned}\right..
\end{align*}
\end{itemize}
\end{thm}

\begin{proof}
Recalling Lemma~\ref{lem:3-col-case}, we already see that if $(a(i))$ is $3$-color, then $(b(j))$ and $(c(k))$ must be equidistant.
It remains to treat the $2$-color case with slopes $0$ and $1$.
As mentioned in the proof of Lemma~\ref{lem:Markoff}, all the possible $1$-balanced words with slope $0$ (resp.\ $1$) are $0^{\Z}$ and ${}^{\infty}010^{\infty}$ (resp.\ $1^{\Z}$ and ${}^{\infty}101^{\infty}$).
If the periodic word $0^{\Z}$ or $1^{\Z}$ appears in at least one direction, then the case falls into the first case in Theorem~\ref{thm:others}.
If the word in each direction is the ``skew'' word ${}^{\infty}010^{\infty}$ or ${}^{\infty}101^{\infty}$, then the situation is reduced to the second case.
\end{proof}

\section{Super Sturmian lattices}
\label{sec:SSL}

We introduce ``super Sturmian lattices'', which describe a hierarchical structure of Sturmian lattices.
There are two types: ``negative'' and ``positive''.
We will use the former to construct an example of aperiodic tile sets, see \S\ref{sec:sqrt2-2}.

\subsection{Negative super Sturmian lattices}
\label{sec:negative}

\begin{figure}[htb]\centering
\subfigure[the original Sturmian lattice]{%
\includegraphics[width = .48\linewidth]{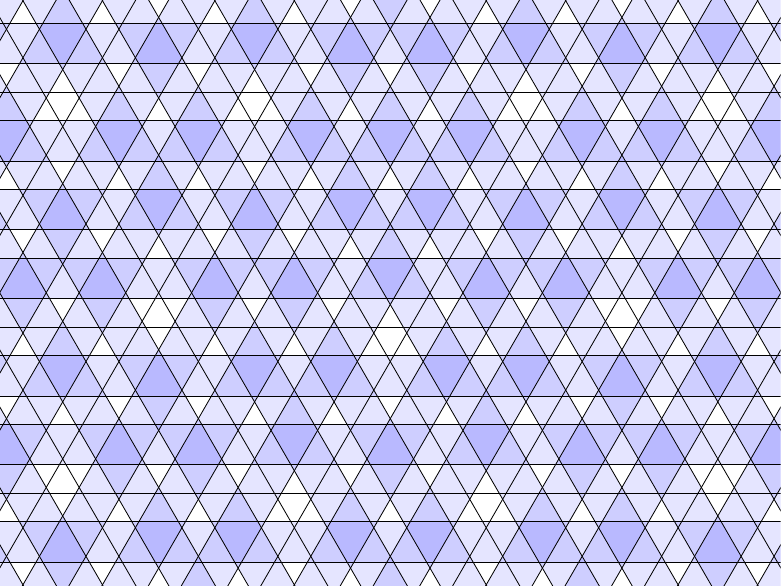}
}\hfill%
\subfigure[the super Sturmian lattice]{%
\includegraphics[width = .48\linewidth]{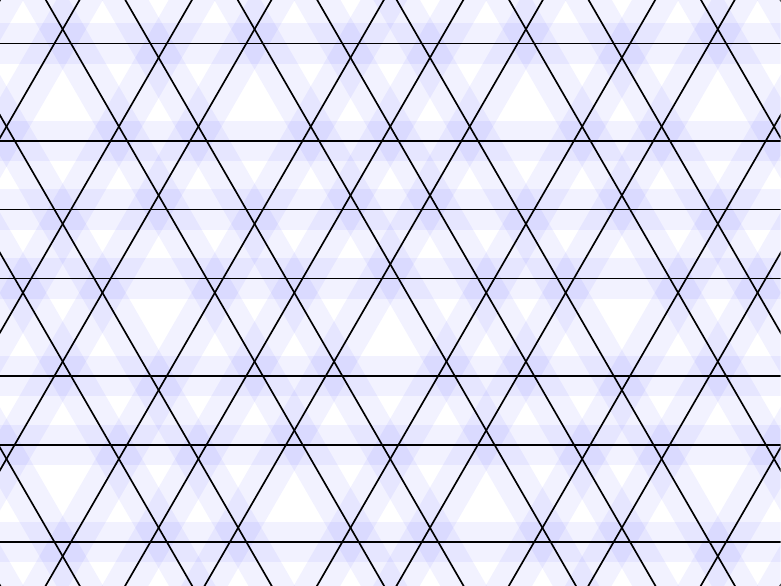}
}
\caption{Example of a super Sturmian lattice}
\label{fig:SSL}
\end{figure}

We consider $2$-color Sturmian lattices with passage $\kappa\ge 1$ and slope $\alpha\in (0, 1]$; it has infinitely many wider corridors.
Geometrically speaking, the ``super Sturmian lattice'' is the collection of middle lines of wider corridors in a given Sturmian lattice.
Figure~\ref{fig:SSL}(b) depicts an example of a super Sturmian lattice.
For $r\in \R$, we denote by
\[
r\SL(\kappa, \alpha):=
\{(r\cdot a(i), r\cdot b(j), r\cdot c(k))\in (\R^{\Z})^{3}\mid
(a(i), b(j), c(k))\in \SL(\kappa, \alpha)\}
\]
the set of $r$-multiple of Sturmian lattices.

We state a sufficient condition for a triple to be a Sturmian lattice.

\begin{lem}
\label{lem:SL}
For a constant $\eps > 0$ and a triple $(a(i), b(j), c(k))$ of real sequences with
\[
a(i+1) - a(i),~
b(j+1) - b(j),~
c(k+1) - c(k) > 2\eps,
\]
we define
\[
V:= \bigl\{
(i, j, k)\in \Z^{3}
\bigm|
\abs{a(i) + b(j) + c(k)} = \eps
\bigr\}.
\]
If $(a(i), b(j), c(k))$ satisfies all of the following conditions:
\begin{itemize}
\item
$\forall j, k\in \Z$, $\exists i = i(j,k)\in \Z$ such that $(i, j, k)\in V$,
\item
$\forall k, i\in \Z$, $\exists j = j(k,i)\in \Z$ such that $(i, j, k)\in V$,
\item
$\forall i, j\in \Z$, $\exists k = k(i, j)\in \Z$ such that $(i, j, k)\in V$,
\end{itemize}
then there exists $s\in \Z$ such that
\[
V = \{(i, j, k)\in \Z^{3}\mid i + j + k = s\}.
\]
In particular, $(a(i - s), b(j), c(k))$ is a Sturmian lattice.
\end{lem}

This lemma produces a Sturmian lattice even without knowing in advance that the tiny triangles are arranged in a lattice pattern.

\begin{proof}
We begin with a key property of the set $V$: for $(i', j', k), (i'', j'', k)\in V$,
\begin{equation}
\label{eq:prf:V}
i' < i''\implies j' > j''.
\end{equation}
Indeed, when $i' + 1\le i''$, the increasing monotonicity of $(a(i))$ gives $a(i'+1)\le a(i'')$, so that
\[
a(i'') - a(i')\ge
\bigl(a(i'') - a(i'+1)\bigr) +
\bigl(a(i'+1) - a(i')\bigr) > 2\eps.
\]
Hence
\begin{align*}
b(j') - b(j'')&\ge
\bigl(a(i' ) + b(j' ) + c(k)\bigr) -
\bigl(a(i'') + b(j'') + c(k)\bigr) +
\bigl(a(i'') - a(i' )\bigr)\\&>
-\eps - \eps + 2\eps = 0.
\end{align*}
Using the monotonicity of $(b(j))$ in turn, we obtain $j' > j''$.

Now we turn to the proof of the lemma. First, we show that the map
\[
V\ni (i, j, k)\mapsto i + j + k\in \Z
\]
is constant.
By assumption, there exists $s\in \Z$ with $(s, 0, 0)\in V$.
Take an arbitrary $(i, j, k)\in V$.
When $j\ge 0$ and $k\ge 0$, the assumption again allows us to choose a chain of points:
\begin{equation}
\label{eq:prf:V2}
\begin{array}{cccc}
&&&(i, j, k),\\
&&&\vdots\\
&&&(*, j, 1),\\
(s, 0, 0), & (*, 1, 0), & \dots, & (*, j, 0)\phantom{,}
\end{array}\in V.
\end{equation}
In general, if $(i+t, j, k), (i, j+1, k)\in V$ then $t = 1$.
Indeed, that $t\ge 1$ follows by repeating the first argument cyclically, so suppose $t\ge 2$ and derive a contradiction.
Taking an integer $i'\in (i, i+t)$, the assumption provides $j'\in \Z$ with
\[
(i+t, j, k), (i', j', k), (i, j+1, k)\in V.
\]
But then $i+t > i' > i$ forces $j < j' < j+1$ by \eqref{eq:prf:V}, which is absurd.
Similarly, if $(i, j+t, k), (i+1, j, k)\in V$ then $t = 1$.
Consequently, the component sum is constant along the chain \eqref{eq:prf:V2}, and in particular $i + j + k = s$.
Hence
\[
V\subset \{(i, j, k)\in \Z^{3}\mid i + j + k = s\}.
\]
For the reverse inclusion, take any element of the right-hand side, written $(s-j-k, j, k)$.
By assumption there exists $i\in \Z$ with $(i, j, k)\in V$, and the inclusion just proved gives $i + j + k = s$, that is, $i = s - j - k$.
Thus $(s-j-k, j, k)\in V$.
\end{proof}

\begin{prop}
\label{prop:SSL}
Let $(a(i), b(j), c(k))\in \SL(\kappa, \alpha)$ be a Sturmian lattice of passage $\kappa\ge 1$ and of slope $\alpha\ne 0$.
If three bi-infinite $\R$-valued sequences $(A(I))_{I\in \Z}$, $(B(J))_{J\in \Z}$, and $(C(K))_{K\in \Z}$ are strictly increasing, and satisfy that
\begin{equation}
\label{eq:SSL}
\begin{aligned}
\{A(I)\in \R\mid I\in \Z\}&=
\left\{
\frac{a(i+1) + a(i)}{2}\in \R
~\middle|~
a_{i} = 1
\right\},\\
\{B(J)\in \R\mid J\in \Z\}&=
\left\{
\frac{b(j+1) + b(j)}{2}\in \R
~\middle|~
b_{j} = 1
\right\},\\
\{C(K)\in \R\mid K\in \Z\}&=
\left\{
\frac{c(k+1) + c(k)}{2}\in \R
~\middle|~
c_{k} = 1
\right\},
\end{aligned}
\end{equation}
then there exists $S\in \Z$ such that the triple
\[
\bigl(
A(I - S),
B(J),
C(K)
\bigr)_{I,J,K}
\]
forms a (possibly $3$-color) Sturmian lattice with $\eps = \kappa/2$.
\end{prop}

\begin{proof}
First, for $I\in \Z$, take $i, i'\in \Z$ with
\begin{align*}
A(I+1)&= \frac{a(i'+1) + a(i')}{2},&
A(I)&= \frac{a(i+1) - a(i)}{2}.
\end{align*}
By the construction of the sequence $(A(I))$, clearly $i < i'$.
Then
\begin{align*}
A(I+1) - A(I)&=
\left(
a(i') + \frac{\kappa + 1}{2}
\right) - \left(
a(i+1) - \frac{\kappa + 1}{2}
\right)\\&=
a(i') - a(i+1) + (\kappa + 1) >
\kappa = 2\eps,
\end{align*}
and similarly $B(J+1) - B(J) > 2\eps$ and $C(K+1) - C(K) > 2\eps$.

For integers $J, K\in \Z$, we wish to find $I\in \Z$ satisfying
\[
\abs{A(I) + B(J) + C(K)} =
\frac{\kappa}{2} = \eps.
\]
By a shift and a translation, we may assume without loss of generality that
\begin{align*}
B(J)&= \frac{b(1) + b(0)}{2} = 0,&
C(K)&= \frac{c(0) + c(-1)}{2} = 0.&
&\left(\begin{aligned}
b(0) = c(-1)&= -\tfrac{\kappa + 1}{2}\\
b(1) = c( 0)&= +\tfrac{\kappa + 1}{2}
\end{aligned}\right)
\end{align*}
Noting that $b(0) + c(0) = 0$, the axiom of a geometric Sturmian lattice gives $a(0) = \pm1/2$.
By a rotation through $\pi$, we may assume without loss of generality that $a(0) = -1/2$.
Then
\[
-\frac{1}{2}\le
a(1) + b(0) + c(-1) =
a(1) - (\kappa + 1)
\]
and, since $(a(i), b(j), c(k))$ is $2$-color,
\[
a(1) - a(0) =
a(1) + \frac{1}{2}\in
\{\kappa, \kappa + 1\}.
\]
Therefore $a(1) = \kappa + 1/2$ is forced, and $\{a(0)\le a\le a(1)\}$ is a wide corridor.
Setting
\[
A(I) = \frac{a(1) + a(0)}{2} = \frac{\kappa}{2},
\]
we obtain $A(I) + B(J) + C(K) = \eps$.
The other directions are analogous.
\end{proof}

\begin{definition}
\label{def:SSL}
For a Sturmian lattice $(a(i), b(j), c(k))\in \SL(\kappa, \alpha)$ with $\kappa\ge 1$ and $\alpha\ne 0$, we define the \emph{(negative) super Sturmian lattice}
\[
\Psi(a(i), b(j), c(k)):= (A(I), B(J), C(K)),
\]
by \eqref{eq:SSL}, up to shift-equivalence.
We call $\Psi$ the \emph{(negative) coding map}.
\end{definition}

\color{red!50!black}
For a Sturmian lattice $(a(i), b(j), c(k))\in \SL(\kappa, \alpha)$ with an irrational slope $\alpha\notin \Q$, let us find the parameters $r_{1}$, $\kappa_{1}$, and $\alpha_{1}$ of the super Sturmian lattice
\[
(A(I), B(J), C(K)) =
\Psi(a(i), b(j), c(k))\in
r_{1}\SL(\kappa_{1}, \alpha_{1}).
\]
By Theorem~\ref{thm:irrational}, the associated symbolic Sturmian lattice $(a_{i}, b_{j}, c_{k})$ satisfies that, if $10^{d}1\prec (a_{i})$, then $d\in \{d_{1}-1, d_{1}\}$ where $d_{1} = \floor{1/\alpha}$.
Thus there are two types of corridors in the super Sturmian lattice; the narrower is of width
\[
1\cdot (\kappa + 1) +
(d_{1}-1)\cdot \kappa =
\kappa\left(
\frac{1}{\kappa} + d_{1}
\right) =: \kappa\kappa_{1},
\]
and the wider is of width $\kappa(\kappa_{1} + 1)$.
By calculating the natural densities of subwords $10^{d_{1}-1}$ and $10^{d_{1}}$ in the Sturmian word $(a_{i})$ of slope $\alpha$, we have
\[
\alpha =
\frac{1}{(1 - \alpha_{1})d_{1} + \alpha_{1}(d_{1} + 1)} =
\frac{1}{d_{1} + \alpha_{1}},
\]
which implies that $\alpha_{1} = 1/\alpha - d_{1}$.
\color{black}

This observation suggests that super Sturmian lattices have a strong relation with the natural extension
\[
\Psi(\kappa, \alpha):= \left(
\frac{1}{\kappa} + \floor{\frac{1}{\alpha}},
\frac{1}{\alpha} - \floor{\frac{1}{\alpha}}
\right)
\]
of the regular continued fraction map, defined on
\[
X = \bigl([1, \infty)\setminus \Q\bigr)
\times
\bigl((0, 1]\setminus \Q\bigr),
\]
see \cite{Dajani-Kraaikamp:02_0, NIT:77, Nakada:81}.
We identify $X$ with $\N^{\Z}$ through $(\kappa, \alpha)\leftrightarrow (d_{n})_{n\in \Z}$, where
\begin{align*}
\kappa&= d_{0} +
\frac{1}{d_{-1} +
\cfrac{1}{d_{-2} +
\cfrac{1}{\ddots}}},&
\alpha&=
\frac{1}{d_{1} +
\cfrac{1}{d_{2} +
\cfrac{1}{\ddots}}}
\end{align*}
are the continued fractions of $\kappa = [d_{0}; d_{-1}, d_{-2}, \dots]$ and $\alpha = [d_{1}, d_{2}, \dots]$.
Then $\Psi$ is conjugate to the left shift on $\N^{\Z}$.
For $n\in \Z$, we write $(\kappa_{n}, \alpha_{n}):= \Psi^{n}(\kappa, \alpha)$, i.e.,
\begin{align*}
\kappa_{n}&:= [d_{n}; d_{n-1}, d_{n-2}, \dots],&
\alpha_{n}&:= [d_{n+1}, d_{n+2}, \dots].
\end{align*}
\color{red!50!black}
By applying $\Psi$ repeatedly to a Sturmian lattice $(a(i), b(j), c(k))\in \SL(\kappa, \alpha)$ with $\alpha\notin \Q$, we have
\[
\Psi^{n}(a(i), b(j), c(k))\in
\kappa_{0}\cdots \kappa_{n-1}\SL(\kappa_{n}, \alpha_{n})
\]
for positive integers $n\ge 1$.
If the passage $\kappa$ is irrational, then the inverse map $\Psi^{-1}$ can be defined and
\[
\Psi^{-n}(a(i), b(j), c(k)) =
\frac{1}{\kappa_{-1}\cdots \kappa_{-n}}
\SL(\kappa_{-n}, \alpha_{-n})
\]
holds for negative integers $-n\le -1$.

Furthermore, we can explicitly specify super Sturmian lattices, including information on the intercepts.
\color{black}

\begin{prop}
\label{prop:negative}
The super Sturmian lattice of $\SL(\kappa, \alpha\mid \bm{\eta}, \bm{\rho})$ is
\[
\Psi\SL(\kappa, \alpha\mid \bm{\eta}, \bm{\rho}) =
\kappa \SL\left(
\kappa_{1}, \alpha_{1}
~\middle|~
\frac{1}{\kappa}\bm{\eta},
-\frac{1}{\alpha}\bm{\rho}
\right).
\]
Equivalently, $\Psi\SL(\kappa, \alpha\mid \bm{\eta}, \bm{\rho}) = (A(I), B(J), C(K))$ where
\begin{align*}
\frac{1}{\kappa}A(I)&=
\frac{I + \eta_{0}}{\kappa} +
\round{\frac{I - \rho_{0}}{\alpha}},\\
\frac{1}{\kappa}B(J)&=
\frac{J + \eta_{1}}{\kappa} +
\round{\frac{J - \rho_{1}}{\alpha}},\\
\frac{1}{\kappa}C(K)&=
\frac{K + \eta_{2}}{\kappa} +
\round{\frac{K - \rho_{2}}{\alpha}}.
\end{align*}
\end{prop}

{\color{HMDcomment}As this formula is not used later, we omit the proof.}

\subsection{Positive super Sturmian lattices}
\label{sec:positive}

We also give a brief explanation of the positive case.
For simplicity, we only consider the case where the slope $\alpha$ is irrational.

\begin{definition}
\label{def:SSL-pos}
For a Sturmian lattice $(a(i), b(j), c(k))\in \SL(\kappa, \alpha)$ with $\kappa\ge 1$ and $\alpha\ne 1$, there exists a unique Sturmian lattice
\[
\Psi^{*}(a(i), b(j),c(k)):=
(A^{*}(I), B^{*}(J), C^{*}(K)),
\]
up to shift-equivalence, such that
\begin{align*}
\{A^{*}(I)\in \R\mid I\in \Z\}&=
\left\{
\frac{a(i+1) + a(i)}{2}\in \R
~\middle|~
a_{i} = 0
\right\},\\
\{B^{*}(J)\in \R\mid J\in \Z\}&=
\left\{
\frac{b(j+1) + b(j)}{2}\in \R
~\middle|~
b_{j} = 0
\right\},\\
\{C^{*}(K)\in \R\mid K\in \Z\}&=
\left\{
\frac{c(k+1) + c(k)}{2}\in \R
~\middle|~
c_{k} = 0
\right\}.
\end{align*}
We call it the \emph{positive super Sturmian lattice} of $(a(i), b(j), c(k))$, and $\Psi^{*}$ is the \emph{positive coding map}.
\end{definition}

\begin{figure}[htb]\centering
\subfigure[the original Sturmian lattice]{%
\includegraphics[width = .48\linewidth]{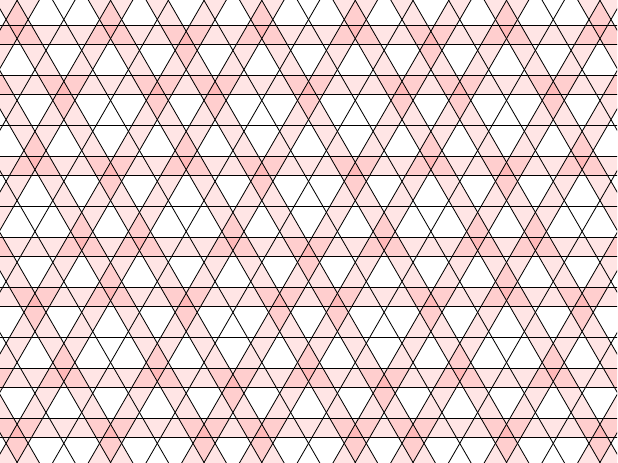}
}\hfill%
\subfigure[the super Sturmian lattice]{%
\includegraphics[width = .48\linewidth]{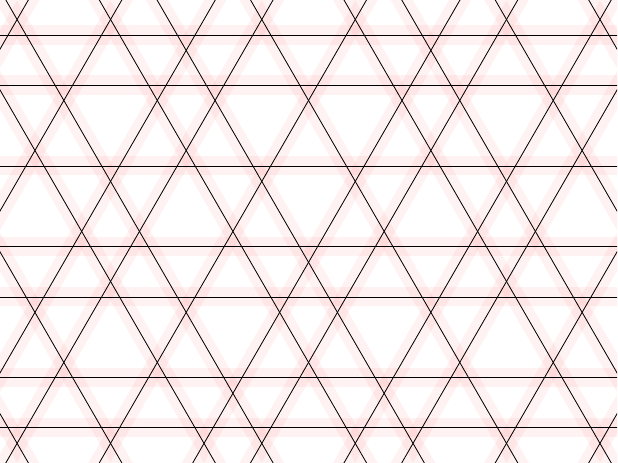}
}
\caption{Example of a positive super Sturmian lattice}
\label{fig:posSSL}
\end{figure}

See Figure~\ref{fig:posSSL} for an example.

In contrast to the case for negative super Sturmian lattices and the associated substitution $\Psi$, the wider and the narrower corridors have switched their roles.
We prepare the dual parameter systems
\begin{align*}
\kappa^{*}&:= \kappa + 1,&
\alpha^{*}&:= 1 - \alpha,&
\bm{\eta}^{*}&:=  \bm{\eta},&
\bm{\rho}^{*}&:= -\bm{\rho},
\end{align*}
and rewrite $\SL(\kappa, \alpha\mid \bm{\eta}, \bm{\rho}) = (a(i), b(j), c(k))$ as
\begin{align*}
a(i)&= i\kappa^{*} + \eta^{*}_{0} - \round{i\alpha^{*} + \rho^{*}_{0}},\\
b(j)&= j\kappa^{*} + \eta^{*}_{1} - \round{j\alpha^{*} + \rho^{*}_{1}},\\
c(k)&= k\kappa^{*} + \eta^{*}_{2} - \round{k\alpha^{*} + \rho^{*}_{2}}.
\end{align*}
We denote by $\SL^{*}(\kappa^{*}, \alpha^{*}\mid \bm{\eta}^{*}, \bm{\rho}^{*}):= (a(i), b(j), c(k))$.
These parameter systems are ``dual'' in the following sense:
\begin{itemize}
\item
Passage $\kappa$ is the narrower width, whereas $\kappa^{*}$ is the wider.
\item
Slope $\alpha$ is the natural density of the wider corridors, whereas $\alpha^{*}$ is that of the narrower.
\end{itemize}

As in Proposition~\ref{prop:negative}, we have

\begin{prop}
\label{prop:positive}
The positive super Sturmian lattice of $\SL^{*}(\kappa^{*}, \alpha^{*}\mid \bm{\eta}^{*}, \bm{\rho}^{*})$ is
\[
\Psi^{*}\SL^{*}(\kappa^{*}, \alpha^{*}\mid \bm{\eta}^{*}, \bm{\rho}^{*}) =
\kappa^{*}\SL^{*}\left(
\kappa^{*}_{1}, \alpha^{*}_{1}
~\middle|~
\frac{1}{\kappa^{*}}\bm{\eta}^{*},
\frac{1}{\alpha^{*}}\bm{\rho}^{*}
\right)
\]
where
\begin{align*}
\kappa^{*}_{1}&:=
\ceil{\frac{1}{\alpha^{*}}} - \frac{1}{\kappa^{*}},&
\alpha^{*}_{1}&:=
\ceil{\frac{1}{\alpha^{*}}} - \frac{1}{\alpha^{*}}.
\end{align*}
Equivalently, $\Psi^{*}\SL^{*}(\kappa^{*}, \alpha^{*}\mid \bm{\eta}^{*}, \bm{\rho}^{*}) = (A^{*}(I), B^{*}(J), C^{*}(K))$ where
\begin{align*}
\frac{1}{\kappa^{*}}A^{*}(I)&=
\frac{-I + \eta^{*}_{0}}{\kappa^{*}} - \round{\frac{-I + \rho^{*}_{0}}{\alpha^{*}}},\\
\frac{1}{\kappa^{*}}B^{*}(J)&=
\frac{-J + \eta^{*}_{1}}{\kappa^{*}} - \round{\frac{-J + \rho^{*}_{1}}{\alpha^{*}}},\\
\frac{1}{\kappa^{*}}C^{*}(K)&=
\frac{-K + \eta^{*}_{2}}{\kappa^{*}} - \round{\frac{-K + \rho^{*}_{2}}{\alpha^{*}}}.
\end{align*}
\end{prop}

We obtain the following transition:
\begin{equation}
\label{eq:negativeNE}
\Psi^{*}(\kappa^{*}, \alpha^{*}) =
(\kappa^{*}_{1}, \alpha^{*}_{1}) =
\left(
\ceil{\frac{1}{\alpha^{*}}} -
\frac{1}{\kappa^{*}},
\ceil{\frac{1}{\alpha^{*}}} -
\frac{1}{\alpha^{*}}
\right).
\end{equation}
The map $\Psi^{*}$ in \eqref{eq:negativeNE} is defined almost everywhere on $X = (1, \infty)\times (0, 1)$, and this is conjugate to a certain sub-shift on $\N^{\Z} = \{(d^{*}_{n})_{n\in \Z}\mid d^{*}_{n}\in \N\}$ where
\begin{align*}
\kappa^{*}&= d^{*}_{0} -
\cfrac{1}{d^{*}_{-1} -
\cfrac{1}{d^{*}_{-2} -
\cfrac{1}{\ddots}}},&
\alpha^{*}&= 
\cfrac{1}{d^{*}_{1} -
\cfrac{1}{d^{*}_{2} -
\cfrac{1}{\ddots}}}.
\end{align*}
Hereafter we simply write $\kappa^{*} = [d^{*}_{0}; d^{*}_{-1}, d^{*}_{-2}, \dots]^{*}$ and $\alpha^{*} = [d^{*}_{1}, d^{*}_{2}, \dots]^{*}$.

\begin{lem}
\label{lem:CFE-relation}
Let $\alpha = [d_{1}, d_{2}, \dots]\in (0, 1)$ be irrational.
Then the negative continued fraction expansion of $\alpha^{*} = 1 - \alpha$ is
\[
\alpha^{*} = [~
\underbrace{2, \dots, 2, d_{2}+2}_
{d_{1}},~
\underbrace{2, \dots, 2, d_{4}+2}_
{d_{3}},~
\underbrace{2, \dots, 2, d_{6}+2}_
{d_{5}},~\dots
~]^{*}.
\]
In particular, $\alpha^{*} = [d^{*}_{1}, d^{*}_{2}, \dots]^{*}$ is purely (resp.\ eventually) periodic if and only if $\alpha = [d_{1}, d_{2}, \dots]$ is purely (resp.\ eventually) periodic.
\end{lem}

\begin{proof}
See \cite[Proposition~2]{Dajani-Kraaikamp:00} and \cite[Aufgaben~3, p.~131]{Zagier}.
\end{proof}

\begin{rem}
\label{rem:sgn}
Recalling Proposition~\ref{prop:negative}, the negative coding map $\Psi$ acts linearly on the $4$-dimensional parameter space
\[
\bigl\{
(\bm{\eta}, \bm{\rho})\in \R^{6}
\bigm|
\eta_{0} + \eta_{1} + \eta_{2} = 0,~
\rho_{0} + \rho_{1} + \rho_{2} = 0
\bigr\}
\]
as
\[
\Psi(\bm{\eta}, \bm{\rho}):= \left(
\frac{1}{\kappa}\bm{\eta},
-\frac{1}{\alpha}\bm{\rho}
\right),
\]
which has a negative determinant.
On the other hand, the positive coding map $\Psi^{*}$ acts as
\[
\Psi^{*}(\bm{\eta}^{*}, \bm{\rho}^{*}):=
\left(
\frac{1}{\kappa}\bm{\eta},
\frac{1}{\alpha}\bm{\rho}
\right),
\]
which has a positive determinant.
That is the reason why we name them ``negative'' or ``positive''.
\end{rem}

%

\subsection{Fundamental Sturmian lattices}
\label{sec:exp-const}

If a Sturmian lattice has a quadratic slope $\alpha$, then it is asymptotically self-similar.
This property is inherited from the corresponding self-similarity of Sturmian words.
We introduce ``fundamental Sturmian lattices'', which describe the self-similarity of Sturmian lattices most concisely.

A Sturmian word whose slope $\alpha = [\overline{d_{1}, d_{2}, \dots, d_{\ell}}]$ is purely periodic, that is, a reduced quadratic irrational in the sense of Gauss, is of special interest.
It is a fixed point of a certain
substitution (see \cite[Theorem~2.3.25]{Fogg})
whose expansion constant $\lambda = \lambda(\alpha)$ is the dominant eigenvalue $>1$ of
\[
M(\alpha) = \begin{pmatrix}
p_{\ell-1} & p_{\ell}\\
q_{\ell-1} & q_{\ell}
\end{pmatrix} =
\begin{pmatrix}
0& 1\\
1& d_{1}
\end{pmatrix}
\begin{pmatrix}
0& 1\\
1& d_{2}
\end{pmatrix}
\cdots
\begin{pmatrix}
0& 1\\
1& d_{\ell}
\end{pmatrix}.
\]
This expansion constant is a fundamental unit of the integer ring of the quadratic field $\Q(\alpha)$.
The self-similar structure of Sturmian words gives rise to the self-similar structure of Sturmian lattices and the expansion constant coincides with
\[
\lambda = \prod_{j=1}^{\ell} [d_{j};\overline{d_{\ell+j-1},d_{\ell+j-2},\dots, d_{j}}]
\]
under the 
choice 
$\kappa=[d_{\ell};\overline{d_{\ell-1},d_{\ell-2},\dots, d_{1},d_{\ell}}]$.
The expansion constant of a Sturmian lattice 
of a general quadratic slope 
$\alpha=[d_{1},d_{2},\dots, d_{n},\overline{d_{n+1},d_{n+2},\dots, d_{n+\ell}}]$
is determined by the purely periodic part 
$[\overline{d_{n+1},d_{n+2},\dots, d_{n+\ell}}]$.

For given algebraic integer $\lambda\ge 1$, we wish to find the ``simplest'' Sturmian lattice whose expansion constant is $\lambda$.
Defining the norm $N(\lambda)$ of $\lambda$ by the product of $\lambda$ and its Galois conjugate, $N(\lambda) = \pm 1$ holds.
\begin{itemize}
\item
If $N(\lambda) = -1$, then $\lambda$ has minimal $\Z$-polynomial
\[
\lambda^{2} - h\lambda - 1 = 0
\qquad(h\ge 1).
\]
We consider the set $\SL(\kappa, \alpha)$ where
\begin{align*}
\kappa&:= \lambda = [h; \bar{h}],&
\alpha&:= \frac{1}{\lambda} =
\lambda - h = [\bar{h}].
\end{align*}
\item
If $N(\lambda) = +1$, then $\lambda$ has minimal $\Z$-polynomial
\[
\lambda^{2} - h\lambda + 1 = 0
\qquad(h\ge 3).
\]
We consider the set $\SL^{*}(\kappa^{*}, \alpha^{*})$ where
\begin{align*}
\kappa^{*}&:=
\lambda = [h; \bar{h}]^{*},&
\alpha^{*}&:=
\frac{1}{\lambda} =
h - \lambda = [\bar{h}]^{*}.
\end{align*}
\end{itemize}

We call an element of $\SL(\lambda, 1/\lambda)$ (if $N(\lambda) = -1$) or $\SL^{*}(\lambda, 1/\lambda)$ (if $N(\lambda) = +1$) a \emph{fundamental Sturmian lattice} of expansion constant $\lambda$.
By discussion on the previous subsections \S\ref{sec:negative} and \S\ref{sec:positive}, the coding map $\Psi$ or $\Psi^{*}$ describes the self-similarity of fundamental Sturmian lattices.

We give a geometric characterization of fundamental Sturmian lattices.
For a Sturmian lattice $(a(i), b(j), c(k))$, a \emph{Sturmian sub-lattice} of $(a(i), b(j), c(k))$ is a Sturmian lattice $(a(i_{I}), b(j_{J}), c(k_{K}))$ for some sub-sequences $(i_{I}), (j_{J}), (k_{K})\in \Z^{\Z}$.

\begin{prop}
\label{prop:fundamental}
Let $\kappa = [d_{\ell}; \overline{d_{\ell-1}, \dots, d_{1}, d_{\ell}}]$ and $\alpha = [\overline{d_{1}, d_{2}, \dots, d_{\ell}}]$, and let $\lambda = \lambda(\alpha) > 1$ the associated expansion constant.
Then, any Sturmian lattice $(a(i), b(j), c(k))\in \SL(\kappa, \alpha)$ has a fundamental Sturmian sub-lattice of slope $1/\lambda$.
\end{prop}

We use a lemma for sub-lattices of Sturmian lattices.

\begin{lem}
\label{lem:sub-lattice}
Let $n\in \N$ and let $(a(i), b(j), c(k))\in \SL(\kappa, \alpha)$ be a Sturmian lattice with irrational $\kappa$ and $\alpha$.
Then $(a(ni), b(nj), c(nk))_{i,j,k}$ is a Sturmian sub-lattice of passage $n\kappa + \floor{n\alpha}$ and slope $n\alpha - \floor{n\alpha}$.
\end{lem}

\begin{proof}
The statement follows by an easy calculation:
\begin{align*}
a(ni)&=
ni\kappa + \round{ni\alpha + \rho_{0}}\\&=
i(n\kappa + \floor{n\alpha}) +
\round{i\{n\alpha\} + \rho_{0}},
\end{align*}
where $\{n\alpha\} = n\alpha - \floor{n\alpha}$ is the fractional part of $n\alpha$.
\end{proof}

\begin{proof}[Proof of Proposition~\ref{prop:fundamental}]
For $n\in \Z/\ell\Z$, we define
\begin{align*}
\kappa_{n}&= [d_{n}; \overline{d_{n-1}, d_{n-2}, \dots, d_{n-\ell}}],&
\alpha_{n}&= [\overline{d_{n+1}, d_{n+2}, \dots, d_{n+\ell}}].
\end{align*}
In particular, $\kappa_{0} = \kappa$ and $\alpha_{0} = \alpha$.
Since
\begin{align*}
\frac{-1}{\kappa_{n-1}}
\begin{pmatrix}
\kappa_{n-1}\\ -1
\end{pmatrix}&= \begin{pmatrix}
0 & 1\\ 1 & d_{n}
\end{pmatrix}\begin{pmatrix}
\kappa_{n}\\ -1
\end{pmatrix},&
\frac{1}{\alpha_{n-1}}
\begin{pmatrix}
\alpha_{n-1}\\ 1
\end{pmatrix}&= \begin{pmatrix}
0 & 1\\ 1 & d_{n}
\end{pmatrix}\begin{pmatrix}
\alpha_{n}\\ 1
\end{pmatrix},
\end{align*}
we have
\begin{align*}
\frac{(-1)^{\ell}}{\lambda}
\begin{pmatrix}
\kappa\\ -1
\end{pmatrix}&= \begin{pmatrix}
p_{\ell-1} & p_{\ell}\\
q_{\ell-1} & q_{\ell}
\end{pmatrix}\begin{pmatrix}
\kappa\\ -1
\end{pmatrix},&
\lambda
\begin{pmatrix}
\alpha\\ 1
\end{pmatrix}&= \begin{pmatrix}
p_{\ell-1} & p_{\ell}\\
q_{\ell-1} & q_{\ell}
\end{pmatrix}\begin{pmatrix}
\alpha\\ 1
\end{pmatrix}.
\end{align*}
Simply writing $N(\lambda) = (-1)^{\ell} = \mp1$, it follows that
\begin{align*}
\pm\frac{1}{\lambda}&=
q_{\ell-1}\kappa - q_{\ell},&
\lambda&=
q_{\ell-1}\alpha + q_{\ell},
\end{align*}
and $\lambda = h\pm 1/\lambda$ implies that
\begin{align*}
\lambda&=
q_{\ell-1}\kappa  + (h - q_{\ell}),&
\pm\frac{1}{\lambda}&=
q_{\ell-1}\alpha - (h - q_{\ell}).
\end{align*}

If $N(\lambda) = -1$, then $1/\lambda\in (0, 1)$ and $h - q_{\ell}\in \Z$ follow that $h - q_{\ell}  = \floor{q_{\ell-1}\alpha}$, and hence
\begin{align*}
\lambda&= q_{\ell-1}\kappa + \floor{q_{\ell-1}\alpha},&
\frac{1}{\lambda}&= q_{\ell-1}\alpha - \floor{q_{\ell-1}\alpha}
\end{align*}
as claimed.
If $N(\lambda) = +1$, then $h - q_{\ell} - 1 = \floor{q_{\ell-1}\alpha}$ implies that
\begin{align*}
\lambda - 1&= q_{\ell-1}\kappa + \floor{q_{\ell-1}\alpha},&
1 - \frac{1}{\lambda}&= q_{\ell-1}\alpha - \floor{q_{\ell-1}\alpha}.
\end{align*}
Thus $\SL(\lambda - 1, 1 - 1/\lambda) = \SL^{*}(\lambda, 1/\lambda)$ is a sub-lattice of $\SL(\kappa, \alpha)$.
\end{proof}

\section{Aperiodic tile sets}
\label{sec:aperiodic}

We consider the Euclidean plane $\R^{2}$. Let $F$ be a finite set. 
A \emph{tile} is a pair $\T = (T, a)$ where $T\subset\R^2$ is a non-empty compact set in $\R^2$ which coincides with the closure of its interior and $a\in F$. 
Note that it is important not to assume any further topological property of $T$ in this paper.
We say that $T$ is the \emph{support} of $\T$ and $a$ is the \emph{color} of $\T$.
A \emph{prototile set} $\AA = \{\T_{1}, \dots \T_{m}\}$ is a finite set of tiles with $\T_{i} = (T_{i}, a_{i})$.
An isometry $g\in \Isom(\R^{2})$ acts on $\T = (T, a)$ by $g((T,a)) = (g(T), a)$.
A collection
\[
\P = \{g_{i,j}(\T_i)\mid g_{i,j}\in J_{i},\ i\in \{1,2,\dots,m\}\}
\]
with $J_{i}\subset \Isom(\R^2)$
is called a \emph{patch} if $g_{i,j}(T_i)$'s have disjoint interiors.
The \emph{support} of a patch $\P$ is defined by
\[
\supp(\P)=\bigcup_{i=1}^m \bigcup_{g_{i,j}\in J_i} g_{i,j}(T_i).
\]
The color is often used to distinguish tiles having congruent supports.
A \emph{translation} of a patch $\P$ by a vector $u \in \R^2$ is defined as
\[
\P+u=\{ (T+u,a)\mid (T,a)\in \P \}
\]
where $T + u = \{x + u\mid x\in T\}$.
A \emph{tiling} by $\AA$ is a patch $\P$ that $\supp(\P)=\R^2$ 
and we use symbol $\TT$
for a tiling. 
This definition of tiles and 
tilings is standard and found in \cite[Definition~5.2]{BG}, 
which generalize the ones in \cite{GS}.
We say that $\AA$ \emph{admits a tiling} if there exists a tiling $\TT$ by $\AA$.
An element of
\[
\Per(\TT) = \{u\in \R^2\mid \TT + u = \TT\}
\]
is called a \emph{period}.
A tiling is \emph{non-periodic} if $\Per(\TT) = \{0\}$ holds.
A prototile set $\AA$ is called \emph{aperiodic} if $\AA$ admits a tiling but all tilings by $\AA$ are non-periodic.
A tiling has \emph{finite local complexity} (FLC for short) 
for any $R>0$, there exists $n\in \N$ that 
the number of patches whose support is contained in the ball $B(x,R)$ centered at any $x\in \R^2$ 
does not exceeds $n$ up to isometry\footnote{Note that FLC is often defined with respect to translations, a subgroup of $\mathrm{Isom}(\R^2)$.}.
We also use a notion of \emph{patch-tile}. 
Given a tiling, if we can partition it into a finite set of patches up to rigid motion that each patch is understood as a new tile, then we call it a 
patch-tile, see \cite{Akiyama-Araki:23}. 

The prototiles often have markers or decorations, 
which restrict the adjacency of prototiles. 
In this paper, our prototiles are polygons or a union of polygons. 
The marker is 
given by some segments written on the ``surface'' of prototiles. 
We assume a matching rule that 
the segment continues across the edges of prototiles
and form a line. This matching rule is introduced by R.~Ammann 
and called \emph{Ammann bars} (see \cite[Chapter~10.4]{GS}, \cite{AGS:92}), 
as we saw in Figures~\ref{fig:Turtle}~and~\ref{fig:intro}.
Such markers or decorations are understood as 
colors associated with tiles in the finite set $F$ together with
corresponding 
adjacency conditions.
Note that Smith Turtle is particularly special: 
the Ammann bars are dispensable.
We can show that there is no tiling which does not obey this rule, see \cite
[Lemma~2]{Akiyama-Araki:23}.

\section{Bounded Displacement correspondences}
\label{sec:BD}

We introduce a ``bounded displacement correspondence'' between $2$-dimensional lattices.
It is a many-to-many correspondence pairing uniformly bounded lattice points.
We will use this idea to construct aperiodic tile sets on Sturmian lattices in \S\ref{sec:Voronoi}.
In \S\ref{sec:Laczkovich}, we derive some results from study of Laczkovich, which theoretically guarantee the existence of such tile sets.
In \S\ref{sec:cross}, we give necessary BD correspondences to construct the desired tile sets explicitly.

\subsection{Bounded Displacement equivalence}
\label{sec:BDeq}

A set $X \subset \R^d$ is said to be \emph{relatively dense} if there exists $R>0$ such that for any $x\in \R^d$, $\overline{B(x,R)}\cap X\neq \emptyset$.
A set $X \subset \R^d$ is said to be \emph{uniformly discrete} if there exists $r>0$ such that for any $x\in \R^d$, $\Card(B(x,r)\cap X)\le 1$. A set $X\subset \R^d$ is a \emph{Delone set} if it is relatively dense and uniformly discrete.
Let $X$ be a Delone set.
It follows directly from the definition that for any $x\in \R^d$ 
\[
0 < \liminf_{R\to \infty}
\frac{\Card(B(x,R)\cap X)}{\mu(B(x,R))} 
\le \limsup_{R\to \infty}
\frac{\Card(B(x,R)\cap X)}{\mu(B(x,R))} < \infty
\]
where $\mu$ is the $d$-dimensional Lebesgue measure.
If there exists $\delta=\delta(X)$ with
\[
\delta = \lim_{R\to \infty}
\frac{\Card(B(x,R)\cap X)}{\mu(B(x,R))},
\]
then it is called the \emph{natural density} of $X$.
Two Delone sets $X$ and $Y$ are said to be \emph{bounded displacement equivalent} (abbreviated as BD equivalent) if there exists $r>0$ and a bijection $\phi\colon X\to Y$ such that
\begin{equation}
\label{eq:BD}
\bigl|\phi(x)-x\bigr| \le r.
\end{equation}
This defines an equivalence relation on the set of Delone sets and we denote by 
$X \BDsim Y$, and we call the associated map $\phi$ a \emph{BD bijection}.
Note that
$X \BDsim Y$ implies $\delta(X)=\delta(Y)$.
A Delone set $X$ is {\it uniformly spread}
if $X\BDsim \delta^{-1/d} \Z^d$. 
The characterization of uniformly spread sets is an important problem.
The case $d=1$ is easy.

\begin{prop}
\label{Onedim}
For $d=1$, a Delone set $X\subset \R$ is uniformly spread if and only if there exists $\delta>0$ and a positive constant $C$ such that for any $m\in \Z$, we have
\begin{equation}
\label{eq:BR}
-C < \Card\bigl(X\cap [0,m]\bigr) - \delta |m| < C.
\end{equation}
In this case, $X\BDsim \delta^{-1} \Z$ holds.
\end{prop}

Here $[0, m]$ is the closed 
interval between $0$ and $m$ regardless of the sign of $m\in \Z$.

\begin{proof}
If $X\BDsim \delta^{-1} \Z$, 
then $X = \{x_j\in \R\mid j\in \Z\}$ and
\begin{equation}
\label{eq:onedim-proof}
-r\le x_j-\delta^{-1}j \le r.
\end{equation}
For integer $m > 0$, we have
\begin{align*}
\Card\bigl(X\cap [0,m]\bigr)&\le
\Card\bigl(\delta^{-1}\Z\cap [-r, m+r]\bigr)\\&\le
\ceil{(m + 2r)\delta}\\&\le
\ceil{m\delta} + \ceil{2r\delta}\\&\le
m\delta + \ceil{2r\delta} + 1
\end{align*}
and similarly
\[
\Card\bigl(X\cap [0,m]\bigr)\ge
m\delta - \ceil{2r\delta} - 1.
\]
Thus we obtain \eqref{eq:BR} with $C =  \ceil{2r\delta} + 1$. 
The case $m<0$ is similar.

Conversely we assume \eqref{eq:BR}. We wish to list the elements of $X$ in the increasing order $X=\{\cdots < x_{-1}< x_0<x_1 < \cdots \}$ so that \eqref{eq:onedim-proof}.
By \eqref{eq:BR}, if 
$m\ge C/\delta$ then $X \cap [0,m]\neq \emptyset$. Thus we can take $x_0$
with $|x_0|\le C/\delta + 1$.
Assume that $x_k$ is defined for $0\le k< k_0$ and satisfies
$\abs{x_k - \delta^{-1}k}\le C/\delta + 1$. 
If $m\ge (k_0 + C)/\delta$ then $\Card\bigl(X \cap [0,m]\bigr)\ge k_0 + 1$
and
\[
X' = \bigl(X\cap [0, m]\bigr)\setminus \{x_0,\dots, x_{k_0-1}\}\neq \emptyset.
\]
Thus we define $x_{k_0}$ by the
minimum element of $X'$. Since if we take 
$m < (k_0 - C)/\delta$ then $\Card(X \cap [0,m]) < k_0$, 
we see
\[
\floor{\frac{k_0-C}{\delta}} < x_{k_0} \le 
\ceil{\frac{k_0+C}{\delta}}.
\]
We continue in this way to define
the index of the elements of $X$ in both positive and negative directions.
It follows from the construction that
\[
\abs{x_j - \delta^{-1}j}\le \frac{C}{\delta} + 1
\]
for any $j\in \Z$, and hence $X\BDsim \delta^{-1} \Z$.
\end{proof}

\subsection{Laczkovich's discussion}
\label{sec:Laczkovich}

The bounded displacement equivalence for $d=2$ is pretty
non-trivial
and studied by Laczkovich~\cite{Laczkovich:90} 
in connection to the classical problem of ``squaring a circle''.
Let $e_1=(1,0)$ and $e_2=(0,1)$ and $P=[0,1)\times [0,1)$.
Denote by $\H$ the collection of any
finite union of squares of the form
$P+ i e_1+ j e_2$ ($i, j\in \Z$).
Laczkovich~\cite{Laczkovich:90} gave a criterion that a Delone 
set is uniformly spread in $\R^2$.

\begin{prop}[{\cite[Lemma~3.6]{Laczkovich:90}, \cite[Lemma~2.3]{Laczkovich:92}}]
\label{prop:Lac}
For a Delone set $X\subset \R^{2}$, if 
there exists $C>0$ that for any $H\in \H$ we have
\begin{equation}
\label{eq:LC}
\bigl|\Card(X \cap H) - \delta(X) \mu (H)\bigr| < C p(H)
\end{equation}
holds then $X \BDsim \delta^{-1/2} \Z^2$.
\end{prop}

Here $\mu$ is the 2-dimensional Lebesgue measure (i.e., volume) and $p$
is the 1-dimensional Lebesgue measure of the boundary (i.e., the perimeter).
Note that one can change the lattice $\Z^2$ to any unimodular lattice,
since the condition \eqref{eq:LC} is invariant 
under the action of unimodular matrix, up to the choice of the constant $C$.
The proof relies on the M.~Hall's infinite version of marriage theorem.
See \cite{Laczkovich:92} for the corresponding result for $d\ge 3$.

\begin{cor}
If two Delone sets $X$ and $Y$ satisfy the condition \eqref{eq:LC}, then 
$\delta(X)=\delta(Y)$ implies $X \BDsim Y$.
\end{cor}

The condition \eqref{eq:LC} is not additive.
Even when $H_1\cap H_2=\emptyset$, two inequalities
\[
\bigl|\Card(X \cap H_1)- \delta(X) \mu (H_1)\bigr| < C p(H_1)
\]
and 
\[
\bigl|\Card(X \cap H_2) - \delta(X) \mu (H_2)\bigr| < C p(H_2)
\]
do not immediately imply
\[
\bigl|\Card(X \cap (H_1\cup H_2)) - \delta(X) \mu (H_1 \cup H_2)\bigr| < C 
p(H_1\cup H_2),
\]
since $p(H_1\cup H_2)$ is less than $p(H_1)+p(H_2)$
when $H_1$ and $H_2$ share their boundaries. 
This is a 
difficulty to apply this criterion. Hereafter we show that this condition is naturally satisfied in our setting related to Sturmian lattices.

\begin{figure}[htb]\centering
\begin{overpic}[
width = .9\linewidth, page = 1]
{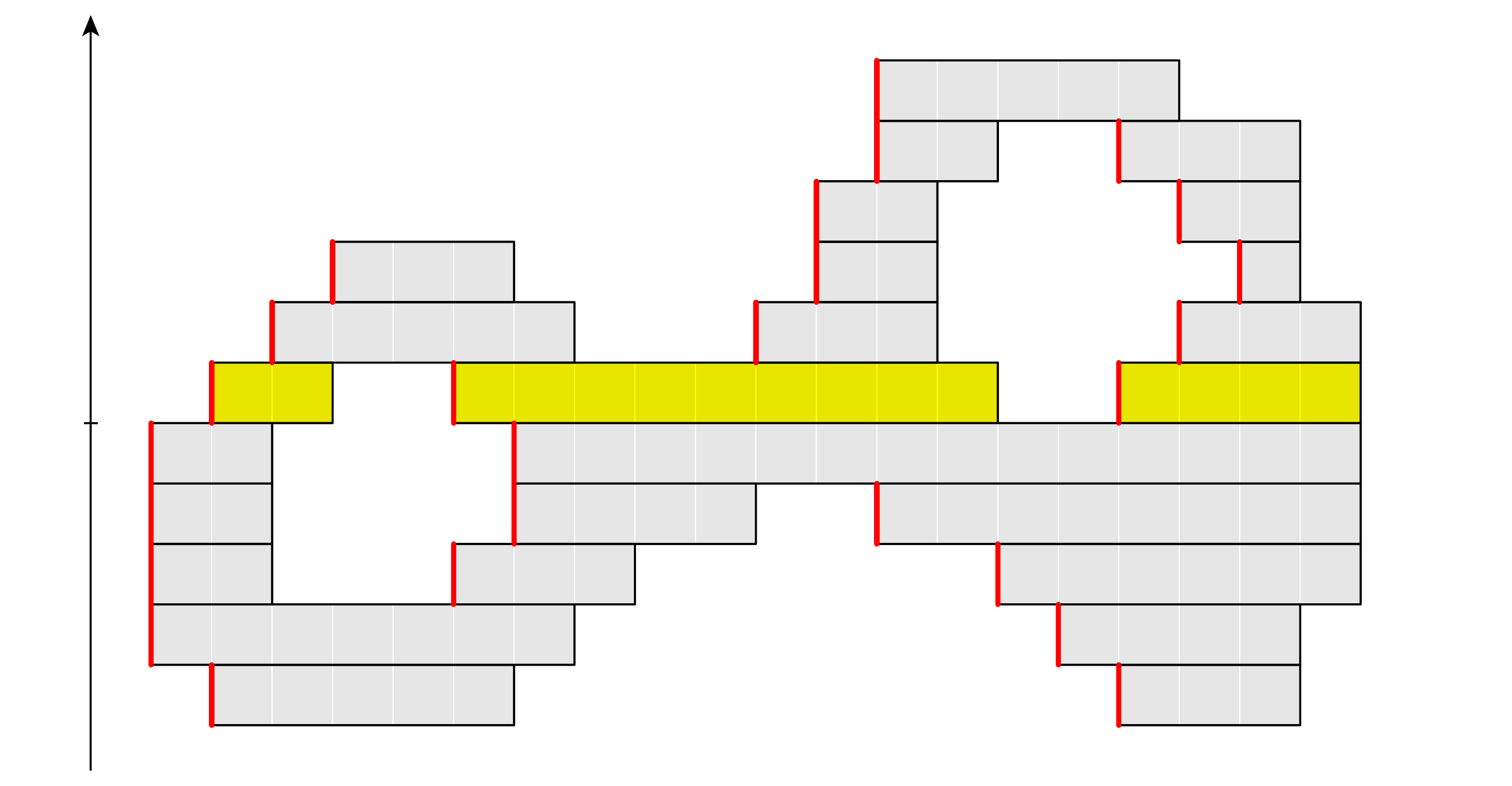}
\put( 5,24){\makebox(0,0)[r]{$y$}}
\put(18,25.8){\makebox(0,0){$H_{y,1}$}}
\put(48,25.8){\makebox(0,0){$H_{y,2}$}}
\put(82,25.8){\makebox(0,0){$H_{y,3}$}}
\put(95,25.8){\makebox(0,0){$H_{y}$}}
\end{overpic}
\caption{Example of $H\in \H$}
\label{fig:Laczkovich}
\end{figure}

We start with an easy but important case.
If $X, Y\subset \R$ satisfy \eqref{eq:BR} with 
natural densities $\delta_1$ and $\delta_2$, then we have
\begin{align*}
X&\BDsim \delta_1^{-1}\Z,&
Y&\BDsim \delta_2^{-1}\Z
\end{align*}
with the common constant $C>0$. 
Then we 
obtain
\begin{equation}
\label{ProductForm}
X \times Y \BDsim \delta_1^{-1}\Z \times \delta_2^{-1}\Z
\end{equation}
with the constant $2C$. 
To see the power of Proposition~\ref{prop:Lac}, 
we prove that 
\[
\delta_1^{-1}\Z \times \delta_2^{-1}\Z \BDsim
(\delta_1 \delta_2)^{-1/2} \Z^2
\]
as an example. Without loss of generality, it suffices to show that
\begin{equation}
\label{Target}
\delta_1^{-1}\Z \times \Z \BDsim
(\delta_1)^{-1/2} \Z^2
\end{equation}
by applying similitude. 
There exists $C_1$ such that 
\[
\Bigl|\Card\left((\delta_1^{-1}\Z \times \{y\})\cap ([m,n]\times \{y\})\right)-\delta_1 |n-m|\Bigr| < C_1
\]
holds for all $y\in \Z$ and all $m,n\in \Z$.
For any $H\in \H$, the strip $H_y=H\cap (\R\times [y,y+1))$ is decomposed into disjoint connected components
\[
H_y = H_{y,1} \sqcup H_{y,2} \sqcup \dots \sqcup H_{y,s(y)}.
\]
Figure~\ref{fig:Laczkovich} helps us understand this setting.
Then
\[
\bigl|\Card\bigl(H_{y,i} \cap (\delta_1^{-1}\Z \times \Z)\bigr)- \delta_1 \mu(H_{y,i})\bigr| < C_1.
\]
Therefore
\[
\bigl|\Card\bigl(H_y \cap (\delta_1^{-1}\Z \times \Z)\bigr)- \delta_1 \mu(H_y)\bigr| < C_1 s(y)
\]
and we obtain
\[
\bigl|\Card\bigl(H \cap (\delta_1^{-1}\Z \times \Z)\bigr)- \delta_1 \mu(H)\bigr| < C_1 \sum_{y}s(y).
\]
Since $\sum_{y}s(y)$ gives the total length of the left boundaries of $H$, we see $2\sum_{y}s(y) < p(H)$. 
Thus \eqref{eq:LC} holds and we obtain \eqref{Target}. The proof works 
in the same way for any lattice in $\R^2$ of density $\delta$, and it is
bounded displacement equivalent to $\delta^{-1/2}\Z^2$.

We use the multi-set notation in \cite[\S2]{Lagarias-Wang:03}. A multi-set is defined as a set of points in which each point is assigned a non-negative integer multiplicity. For two multi-sets $X$, $Y$, the multi-set union
$X \vee Y$ gives the multi-set defined by the sum of the corresponding multiplicities. The usual set is a multi-set that all of its points are 
of multiplicity one. We use this identification for brevity.

This notation is not essential in this paper, 
but it makes the 
description easier. For example, the union of two Delone sets may not be uniformly discrete and it may have coincidence of points, but 
bounded displacement equivalence is extended in a natural way to 
multi-sets; by treating a point of multiplicity $m$ as distinct $m$ points. Since the proof of Proposition~\ref{prop:Lac} depends on the counting argument of 
cardinality of appropriately defined neighbors to apply marriage theorem, 
it immediately extends to multi-set settings.

\begin{prop}
\label{prop:add-sub}
Let $X, Y\subset \R^{2}$ be uniformly spread Delone sets with natural densities 
$\delta_{1} = \delta(X)$ and $\delta_{2} = \delta(Y)$ respectively.
Then we have
\[
X\vee Y\BDsim \delta^{-1/2}\Z^{2} \text{ with }  \delta = \delta_{1} + \delta_{2}.
\]
If we assume further that $Y\subset X$, then 
\[
X\setminus Y\BDsim \delta^{-1/2}\Z^{2} \text{ with } 
\delta = \delta_{1} - \delta_{2}.
\]
\end{prop}

\begin{proof}
The first statement directly
follows from Proposition~
\ref{prop:Lac} since in the multi-set notation, the function 
$p(H)$ is the sum of the contributions from
$X$ and $Y$ independently. Thus, if
$X$ and $Y$ satisfy \eqref{eq:LC}, then so does $X\vee Y$.
The second statement follows in the same line. There exist $C_1,C_2$ such that
for any
$H\in \H$, 
we have
\[
\bigl|\Card(X \cap H)- \delta_1 \mu(H)\bigr| < C_1 p(H)
\]
and
\[
\bigl|\Card(Y \cap H) - \delta_2 \mu(H)\bigr| < C_2 p(H).
\]
Therefore
\[
\bigl|\Card\bigl((X\setminus Y) \cap H\bigr) - (\delta_1-\delta_2) \mu(H)\bigr| < (C_1+C_2) p(H)
\]
and \eqref{eq:LC} is satisfied.
\end{proof}

\begin{prop}
\label{prop:existence-sub}
Let $X\subset \R^{2}$ be uniformly spread Delone set and let $\delta\in (0, 1)$.
Then there exists a uniformly spread Delone subset $Y\subset X$ with natural density $\delta(Y) = \delta$.
\end{prop}

\begin{proof}
We may assume that $X = \Z^{2}$ without loss of generality, and let
\[
Y = \bigl\{\phi(i,j) = \bigl(\floor{\tfrac{1}{\delta}i}, j\bigr)\in \Z^{2}\bigm| i, j\in \Z\bigr\}.
\]
Since $\bigl|\floor{x} - x\bigr|\le 1$ for $x\in \R$, we have
\[
\bigl|\phi(i,j) - (\tfrac{1}{\delta}i, j)\bigr|\le \sqrt{2},
\]
which shows $Y\BDsim \frac{1}{\delta}\Z\times \Z\BDsim \delta^{-1/2}\Z^{2}$.
\end{proof}

\begin{cor}
\label{cor:BD-correspondence}
Let $X, Y\subset \R^{2}$ be uniformly spread Delone sets with natural densities $\delta_{1} = \delta(X)$ and $\delta_{2} = \delta(Y)$.
Then $\delta_{1}/\delta_{2} = p/q\in \mathbb{Q}$ implies the existence of countable partitions
\begin{align*}
X&= \bigsqcup_{k\ge 1}X_{k},&
Y&= \bigsqcup_{k\ge 1}Y_{k}
\end{align*}
of $X$ and $Y$, which satisfy that $\Card(X_{k}) = p$, $\Card(Y_{k}) = q$ for each $k$, and that $X_{k}\vee Y_{k}$ is uniformly bounded.
\end{cor}

\begin{proof}
Since there is a bijection $\phi\colon X\to \delta_1^{-1/2}\Z^2$ that induces $X\BDsim \delta_1^{-1/2}\Z^2$, 
a disjoint union
\[
\Z^2=\bigsqcup_{i,j\in \Z} \bigl(\{p i,p i+1,\dots, p i+p-1\}\times \{j\} \bigr)
\]
gives rises to 
\[
X_{i,j}=\phi^{-1} \left(\delta_1^{-1/2}\bigl(
\{p i,p i+1,\dots, p i+p-1\}\times \{j\} \bigr)\right) 
\]
and the subset
\[
X(p)=\phi^{-1}\left(\delta_1^{-1/2} (
p \Z \times \Z) \right)
\]
with the natural density 
$\delta(X(p))=\delta_1/p$. Similarly we define
$Y_{i,j}$ and $Y(q)$. Since $X(p)$, $Y(q)$
are uniformly spread and $\delta(X(p))=\delta(Y(q))$ we have
\[
X(p)\BDsim Y(q).
\]
By construction, each $X_{i,j}\vee Y_{i,j}$ is uniformly bounded.
\end{proof}

Note that Corollary~\ref{cor:BD-correspondence} is equivalent to the statement: there exist a uniformly spread Delone set $\Lambda$, a constant $r > 0$, and two surjections $\phi\colon X\xrightarrow{p:1} \Lambda$, $\psi\colon Y\xrightarrow{q:1} \Lambda$ such that
\begin{align*}
\abs{x - \phi(x)}&\le r,&
\abs{y - \psi(y)}&\le r
\end{align*}
for any $x\in X$ and $y\in Y$.
These maps $\phi$ and $\psi$ are said to be \emph{bounded displacement surjective} (abbreviated as BD surjective).
We often write $X\stackrel{p}{\to} \Lambda \stackrel{q}{\gets} Y$, where $\Lambda$ frequently serves as a dummy and may be left unspecified.
We call such a pair of maps a \emph{bounded displacement pair} (abbreviated as a BD pair).

\subsection{Cross-BD pairs}
\label{sec:cross}

\begin{figure}[htb]\centering
\includegraphics[width = \linewidth]
{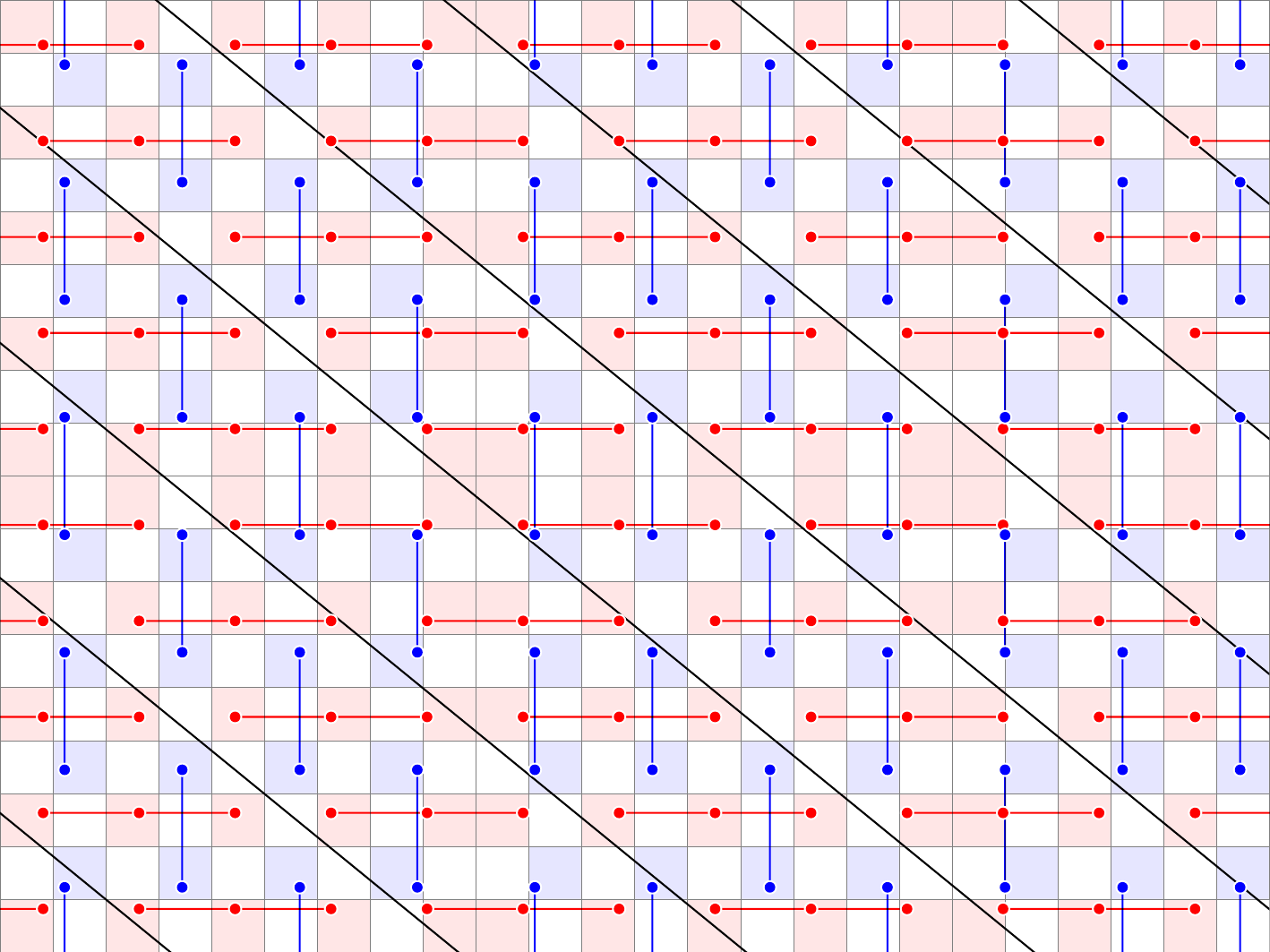}
\caption{Configuration of a cross-BD pair $X\stackrel{3}{\to} \Lambda\stackrel{2}{\gets} Y$}
\label{fig:sqrt6-BD}
\end{figure}

We construct a BD pair $X\stackrel{p}{\to} \Lambda\stackrel{q}{\gets} Y$ between two lattices
\begin{align*}
X&= \lambda\Z\times \frac{1}{\lambda p}\Z,&
Y&= \frac{1}{\mu q}\Z\times \mu\Z
\end{align*}
with $\lambda\mu\le 1$ and $p, q\in \N$, called a \emph{cross-BD pair}.
Figure~\ref{fig:sqrt6-BD} shows an example for
\begin{align*}
X&= \frac{1}{\sqrt{2}}\Z\times \frac{1}{\sqrt{2}}\Z,&
Y&= \frac{1}{\sqrt{3}}\Z\times \frac{1}{\sqrt{3}}\Z;
\end{align*}
one can find many ``crosses'' composed of five lattice points!
The horizontally aligned red dots belong to one lattice $X$, and the vertically aligned blue dots belong to the other lattice $Y$.
We will use this example in the configuration of an aperiodic tile set, see \S\ref{sec:sqrt6}.

Duneau and Oguey~\cite{DO:91} studied explicit constructions of BD bijections between arbitrary two $d$-dimensional unimodular lattices $L_{1}$ and $L_{2} = AL_{1}$ ($A\in \operatorname{SL}_{d}(\R)$).
Their key idea is quite simple; applying Gauss's diagonalization procedure, one can decompose the matrix $A = M_{n}M_{n-1}\cdots M_{1}$ into a product of shear matrices $M_{i}$.
Since an explicit BD bijection between $L_{1}$ and $M_{1}L_{1}$ can be given easily, they finally get a desired map $f\colon L_{1}\to L_{2}$ and it is piecewise affine.
We are inspired by their idea and construct cross-BD pairs.


We define one-dimensional lattices
\begin{align*}
X_{i}&= \{\lambda i\}\times \frac{1}{\lambda p}\Z,&
Y_{j}&= \frac{1}{\mu q}\Z\times \{\mu j\}
\end{align*}
for $i, j\in \Z$, and define a family of unbounded strips
\[
B_{k} =
\bigl\{(x, y)\in \R^{2}\bigm| \mu x + \lambda y\in [k, k+1)\bigr\}
\]
for $k\in \Z$.
We have partitions
\begin{align*}
X&= \bigsqcup_{i\in \Z}X_{i},&
Y&= \bigsqcup_{j\in \Z}Y_{j},&
\R^{2}&= \bigsqcup_{k\in \Z}B_{k}
\end{align*}
of the given lattices and the whole plane.
For fixed $k\in \Z$, it is easy to check that
\begin{align*}
\Card(X_{i}\cap B_{k})&= p,&
\Card(Y_{j}\cap B_{k})&= q
\end{align*}
hold for any $i, j\in \Z$.
By our assumption $\lambda \mu\le 1$, there exists $c_{k}\in \Z$ such that $c_{k}\lambda \mu\in [k, k+1)$.
In particular, if $i + j = c_{k}$ then
\[
(\lambda i, \mu j)\in
\biggl\{(x, y)\in \R^{2}\biggm|
\frac{x}{\lambda} + \frac{y}{\mu} = c_{k}\biggr\}\subset B_{k};
\]
thus $(X_{i}\vee Y_{c_{k} - i})\cap B_{k}$ is uniformly bounded with respect to $i$ and $k$.
We finally get suitable partitions
\begin{align*}
X&= \bigsqcup_{k, i\in \Z}
(X_{i}\cap B_{k}),&
Y&= \bigsqcup_{k, i\in \Z}
(Y_{c_{k} - i}\cap B_{k}).
\end{align*}
In other words, we constructed two BD surjections $X\stackrel{p}{\to} \Lambda$ and  $Y\stackrel{q}{\to} \Lambda$ where
\[
\Lambda = \bigl\{
(\lambda i, \mu j)\in \R^{2}
\bigm|
i + j = c_{k},\ i, k\in \Z\bigr\}
\]
is a (uniformly spread) Delone set with natural density one.
We call the point $(\lambda i, \mu j)\in \Lambda$ the \emph{nexus} of component $X_{i}\vee Y_{j}$.

\section{Nuts and Bolts}
\label{sec:Voronoi}

\begin{figure}[htb]\centering\hfill
\subfigure[$\Nut(S) = (0, 0, +, +)$]{%
\includegraphics[page = 2,
width = .30\linewidth]
{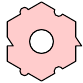}
}\hfill%
\subfigure[$\Nut(M) = (0, 1, -, +)$]{%
\includegraphics[page = 3,
width = .30\linewidth]
{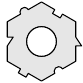}
}\hfill%
\subfigure[$\Nut(L) = (1, 1, +, +)$]{%
\includegraphics[page = 4,
width = .30\linewidth]
{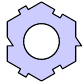}
}\hfill\mbox{}\\ \hfill
\subfigure[Bolt $S$]{%
\includegraphics[page = 5,
width = .30\linewidth]
{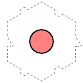}
}\hfill%
\subfigure[Bolt $M$]{%
\includegraphics[page = 6,
width = .30\linewidth]
{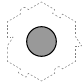}
}\hfill%
\subfigure[Bolt $L$]{%
\includegraphics[page = 7,
width = .30\linewidth]
{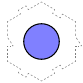}
}\hfill\mbox{}
\caption{Nuts and Bolts}
\label{fig:Nuts-and-Bolts}
\end{figure}

\begin{figure}[htb]\centering
\includegraphics[page = 1,
width = \linewidth, pagebox = artbox]
{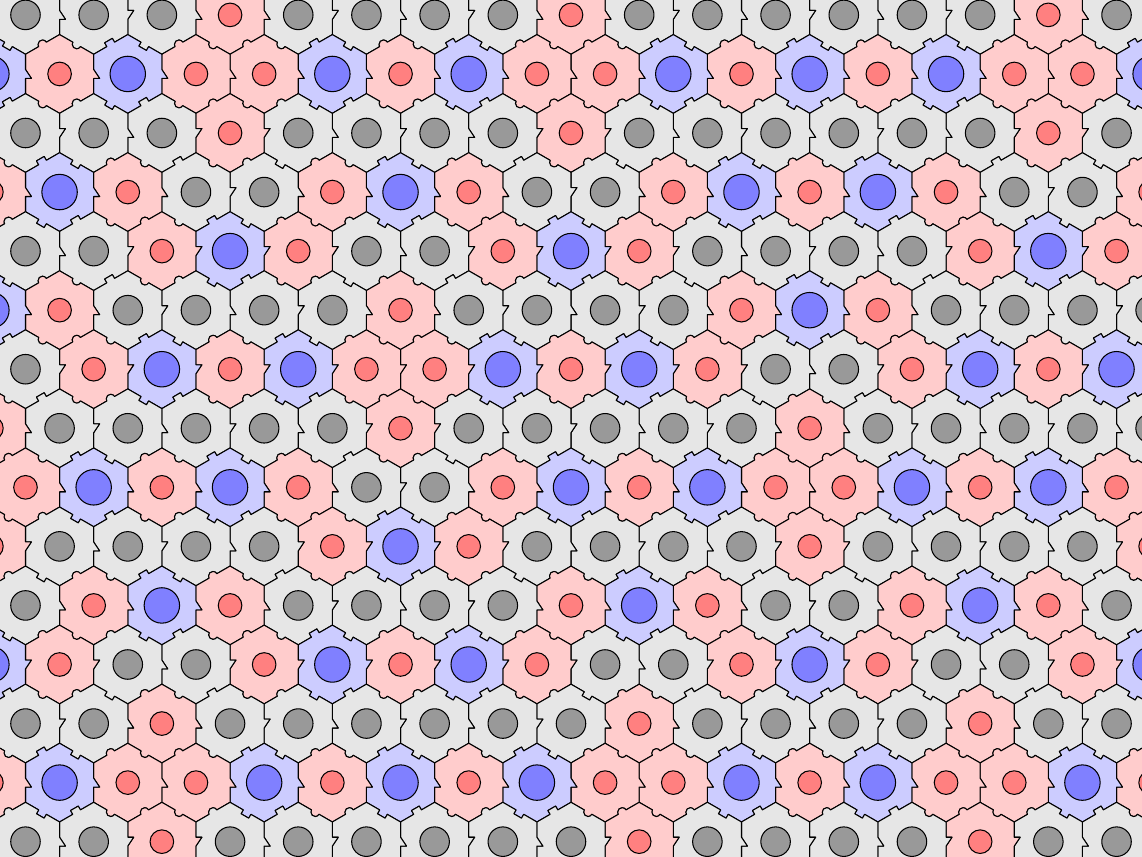}
\caption{Nut-and-Bolt tiling}
\label{fig:Nut-tiling}
\end{figure}

Our purpose is to construct aperiodic tile sets, called ``$\SL(\alpha)$-tile sets'', that realize a Sturmian-lattice structure with a fixed irrational slope $\alpha$.
The $\SL(\alpha)$-tiles are built from three kinds of ``nuts'' and three kinds of ``bolts'', shown in Figure~\ref{fig:Nuts-and-Bolts}: from a deformed regular hexagon we cut out a disk of sufficiently small radius, splitting it into two parts.
Of these, the one homeomorphic to an annulus is called a \emph{Nut}, and the one homeomorphic to a disk a \emph{Bolt}.
They play the following roles, respectively:
\begin{description}
\item[Nuts]
The bumps and dents on the edges of the hexagon enforce the Sturmian-lattice structure.
\item[Bolts]
Applying a BD correspondence enforces a specific slope of the Sturmian lattice.
\end{description}
Figure~\ref{fig:Nut-tiling} shows an example of a Nut-and-Bolt tiling.
The existence of $\SL(\alpha)$-tile sets is guaranteed for every quadratic irrational $\alpha\in (0, 1)$ by Theorem~\ref{thm:Voronoi-aperiodic}(2).

\subsection{Nuts}
\label{sec:Nuts}

\begin{figure}[htb]\centering
\includegraphics[page = 1,
width = .24\linewidth]
{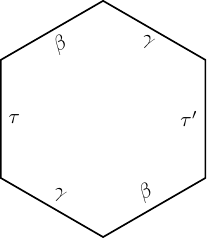}
\caption{The labels of Nut $(\beta, \gamma, \tau, \tau')$}
\label{fig:Nut-label}
\end{figure}

We prepare an uncountably infinite family of tiles, continuously parametrized by a real number $\kappa\in \R$:
\begin{equation}
\label{eq:Nuts}
\begin{aligned}
\tilde{\mathcal{N}}&:=
\bigl\{(\beta, \gamma, \tau, \tau')\in \R^{2}\times \{\pm1/2\}^{2}\bigm|
\beta + \tau' = \gamma + \tau\bigr\}\\&=
\left\{\begin{aligned}
&(\kappa, \kappa, +\tfrac{1}{2}, +\tfrac{1}{2}),~
(\kappa, \kappa + 1, -\tfrac{1}{2}, +\tfrac{1}{2}),\\
&(\kappa, \kappa, -\tfrac{1}{2}, -\tfrac{1}{2}),~
(\kappa + 1, \kappa, +\tfrac{1}{2}, -\tfrac{1}{2})
\end{aligned}~\middle|~
\kappa\in \R
\right\}.
\end{aligned}
\end{equation}
These represent regular hexagons whose edges are deformed as in Figure~\ref{fig:Nuts-and-Bolts}(a)--(c), and the four parameters describe how each edge is deformed.
The correspondence between the edges and the parameters is shown in advance in Figure~\ref{fig:Nut-label}.

There are two kinds of edge deformation, ``symmetric'' and ``asymmetric''.
An edge carrying a bump is called a \emph{bump edge}, and an edge carrying a dent a \emph{dent edge}.
Around each hexagon, bump edges and dent edges alternate;
here the $\tau$-edge is a dent edge and the $\tau'$-edge is a bump edge.
In particular, a dent edge and the bump edge facing it form a matching pair that meshes with each other.
\begin{description}
\item[Symmetric edges]
They form two of the three pairs of opposite edges.
{\color{HMDcomment}For each $\kappa\in \R$, the shape of a symmetric edge is uniquely determined.}
\item[Asymmetric edges]
They form the remaining pair of opposite edges.
In Figure~\ref{fig:Nuts-and-Bolts}(a)--(c), the orientation $\Delta$ corresponds to $+1/2$ and the orientation $\nabla$ to $-1/2$.
\end{description}

In a hexagonal tiling by these Nuts, the asymmetric edges are clearly parallel.
From now on we fix them in the vertical direction.
Moreover, since the bump edges and dent edges of each Nut alternate, if the right (left) vertical edge of one Nut is a bump edge, then so is the right (left) vertical edge of every other Nut.
For this reason, in Nut tilings we use only translations and reflection across a horizontal axis.
Note that the reflection of $(\beta, \gamma, \tau, \tau')\in \tilde{\mathcal{N}}$ across a horizontal axis is $(\gamma, \beta, -\tau, -\tau')\in \tilde{\mathcal{N}}$, so a Nut is characterized by $\{\beta, \gamma\}$ up to isometry; writing this Nut as $N_{\beta, \gamma}$, we have
\[
\tilde{\mathcal{N}}/\Isom(\R^{2}) =
\{N_{\kappa, \kappa}, N_{\kappa, \kappa + 1}\mid \kappa\in \R\}.
\]

\begin{prop}
\label{prop:Nut-tiling}
There is a bijection between axiomatic Sturmian lattices and Nut tilings.
\end{prop}

\begin{proof}
For a given Sturmian lattice $(a(i), b(j), c(k))$, let us make a Nut tiling
\[
\mathscr{N}(a(i), b(j), c(k)) =
\{N(j,k) + \Triplet{*}{j}{k}\mid (j, k)\in \Z^{2}\}.
\]
Here $\Triplet{*}{j}{k} = \Triplet{-j-k}{j}{k}\in \R^{2}$ denotes the triangular coordinates of Figure~\ref{fig:triplet}, which means that the doubly indexed Nuts are arranged as in Figure~\ref{fig:Nut-adjacency}.
For each $(j, k)\in \Z^{2}$, we define the Nut $N(j, k)\in \tilde{\mathcal{N}}$ by
\[
N(j,k) = (\beta_{j}, \gamma_{k}, \tau_{j,k}, \tau_{j-1,k+1})
\]
where
\begin{align*}
\beta_{j}&= b(j+1) - b(j)\in \R,\\
\gamma_{k}&= c(k+1) - c(k)\in \R,\\
\tau_{j,k}&= a(i-1) + b(j+1) + c(k)\in \left\{+\frac{1}{2}, -\frac{1}{2}\right\}
\end{align*}
with $i + j + k = 0$.
Indeed, this satisfies
\[
\beta_{j} + \tau_{j-1,k+1} =
a(i-1) + b(j+1) + c(k+1) =
\gamma_{k} + \tau_{j,k}
\]
and hence defines a Nut by \eqref{eq:Nuts}.
Moreover, the tiling arranged as in Figure~\ref{fig:Nut-adjacency} satisfies the matching rule.

\begin{figure}[htb]\centering
\includegraphics[page = 1,
width = .8\linewidth, pagebox = artbox]
{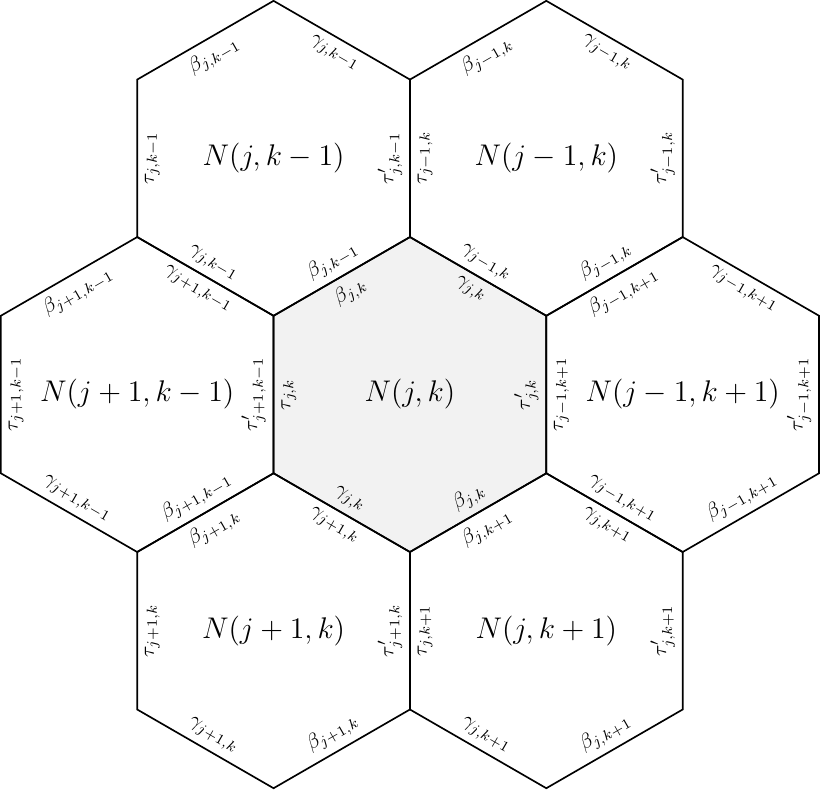}
\caption{Nut adjacency}
\label{fig:Nut-adjacency}
\end{figure}

Conversely, starting from a Nut tiling, we reconstruct a Sturmian lattice.
We wish to find a Sturmian lattice consistent with the markers of the Nuts.
Again indexing the Nut tiling doubly as in Figure~\ref{fig:Nut-adjacency} and writing each Nut as
\[
N(j,k) = (\beta_{j,k}, \gamma_{j,k}, \tau_{j,k}, \tau'_{j,k})
\qquad(j, k\in \Z)
\]
we obtain a system of equations for the triple of sequences $(a(i), b(j), c(k))_{i,j,k}$:
\begin{align}
b(j+1) - b(j)&= \beta_{j,k}\in \R,
\label{eq:Nut-to-SL1} \\
c(k+1) - c(k)&= \gamma_{j,k}\in \R,
\label{eq:Nut-to-SL2} \\
a(i-1) + b(j+1) + c(k)&= \tau_{j,k}\in \left\{+\frac{1}{2}, -\frac{1}{2}\right\}.
\label{eq:Nut-to-SL3}
\end{align}
From the matching rule of the Nut tiling, we have
\begin{align*}
\beta_{j,k}&= \beta_{j,k+1},&
\gamma_{j,k}&= \gamma_{j+1,k},&
\tau'_{j,k}&= \tau_{j-1,k+1}.
\end{align*}
In particular, $\beta_{j}:= \beta_{j,k}$ does not depend on $k$, and $\gamma_{k}:= \gamma_{j,k}$ does not depend on $j$.
We want to find one solution $(a(i), b(j), c(k))$ of this system (it is clear that such a solution gives a Sturmian lattice).

First, we define $b(j)$ and $c(k)$ recursively as follows.
\begin{align*}
b(0)&= 0,&
b(j+1)&:= b(j) + \beta_{j},\\
c(0)&= 0,&
c(k+1)&:= c(k) + \gamma_{k}.
\end{align*}
This clearly satisfies \eqref{eq:Nut-to-SL1} and \eqref{eq:Nut-to-SL2} as well.

It remains to reconstruct $a(i)$.
For each $j, k\in \Z$ with $i + j + k = 0$, we define
\[
a_{j,k}(i-1):= \tau_{j,k} - b(j+1) - c(k)
\]
and we show that this does not depend on $j$ and $k$.
Indeed, applying \eqref{eq:Nuts} to $N(j,k)$ gives
\begin{align*}
a_{j-1,k+1}&=
\tau'_{j,k} - b(j) - c(k+1)\\&=
\tau_{j,k} - b(j+1) - c(k) = a_{j,k}.
\end{align*}
Therefore $a(i-1):= a_{j,k}(i-1)$ satisfies \eqref{eq:Nut-to-SL3}.
\end{proof}

From the way the bijection $\mathscr{N}$ of Proposition~\ref{prop:Nut-tiling} is constructed, the following is clear.

\begin{cor}
\label{cor:Nut-non-periodic}
A Sturmian lattice $(a(i), b(j), c(k))$ is non-periodic if and only if the corresponding Nut tiling $\mathscr{N}(a(i), b(j), c(k))$ is non-periodic.
In particular, the Nut tiling corresponding to a Sturmian lattice of irrational slope is non-periodic.
\end{cor}

From the classification theorems of Sturmian lattices (Theorems~\ref{thm:irrational}--\ref{thm:others}), it follows that a single Nut tiling uses only finitely many kinds of Nut.
From now on, we consider the Nut tilings by
\[
\mathcal{N}:=
\{N_{\kappa, \kappa},
  N_{\kappa, \kappa + 1},
  N_{\kappa + 1, \kappa + 1}\}.
\]
This corresponds to the Sturmian lattices $(a(i), b(j), c(k))$ whose $(b(j), c(k))$ is two-colored.
Rewriting the symmetric-edge labels $\kappa$ and $\kappa + 1$ as $0$ and $1$ respectively, we write
\begin{align*}
\Nut(S)&:= N_{0,0},&
\Nut(M)&:= N_{0,1},&
\Nut(L)&:= N_{1,1}.
\end{align*}
This rewriting operation corresponds to the combinatorial-equivalence transformation of Proposition~\ref{prop:comb-equiv}.

\subsection{Bolts}
\label{sec:Bolts}
Corresponding to the three kinds of Nut, we prepare Bolts $S$, $M$, and $L$.
A Bolt $B\in \{S, M, L\}$ must be combined with $\Nut(B)$.
A \emph{Bolt tiling} is a set of Bolts in the plane
\[
\mathscr{B} = \{B(j,k)\mid (j, k)\in \Z^{2}\}
\]
such that there exists a Nut tiling
\[
\mathscr{N} = \{N(j,k)\mid (j, k)\in \Z^{2}\}
\]
in which each Bolt $B(j,k)$ is paired with the Nut $N(j,k)$.
From now on, by taking its centroid we regard a Bolt as a colored point in the plane $\R^{2}$.
A Bolt tiling is a unimodular triangular lattice, which we write as follows.
\[
\mathscr{B}:=
\left\{
B(j,k) = \Bigl(
\Triplet{*}
{~j + \frac{1}{2}}
{~k + \frac{1}{2}},~\ell(j,k)
\Bigr)
~\middle|~
\begin{aligned}
(j, k)&\in \Z^{2},\\
\ell(j,k)&\in \{S, M, L\}
\end{aligned}
\right\}.
\]
We denote by $\mathcal{S}$, $\mathcal{M}$, $\mathcal{L}$ the sets of all Bolts of type $S$, $M$, $L$ in a Bolt tiling, respectively.
Setting the natural density of the Bolt tiling to one, we write the natural density of each cell as $\delta(\mathcal{S})$, $\delta(\mathcal{M})$, $\delta(\mathcal{L})$.

\begin{prop}
\label{prop:Bolt-nd}
Let $(a(i), b(j), c(k))$ be a Sturmian lattice whose $(b(j), c(k))$ is two-colored.
Then, in the corresponding Nut tiling $\mathscr{N}(a(i), b(j), c(k))$,
\begin{itemize}
\item
if $(a(i))$ is exactly three-colored, then $\delta(\mathcal{S}) = 1$ or $\delta(\mathcal{L}) = 1$;
\item
if $(a(i), b(j), c(k))$ is two-colored, then, writing the slope as $\alpha\in [0, 1]$, we have
\begin{equation}
\label{eq:Bolts-nd}
\delta(\mathcal{S}) :
\delta(\mathcal{M}) :
\delta(\mathcal{L}) =
(1 - \alpha)^{2} :
2\alpha(1 - \alpha) :
\alpha^{2}.
\end{equation}
\end{itemize}
\end{prop}

\begin{proof}
If $(a(i))$ is three-colored, then by Theorem~\ref{thm:others} the other two directions are one-colored.
Otherwise, the claim is clear from the definition of the slope.
\end{proof}

We define the \emph{slope} of a Bolt tiling as the slope of the corresponding Sturmian lattice.
Then the natural density of each cell turns out to be a complete invariant for the slope.

Moreover, $\mathcal{S}$, $\mathcal{M}$, $\mathcal{L}$ are uniformly spread sets, as outlined below:
indeed, in a Bolt tiling $\mathscr{B} = \Z\times \Z$ of slope $\alpha$, noting for instance that
\[
\mathcal{S} = \{B(j, k)\in \mathscr{B}\mid b_{j} = c_{k} = 0\}
\]
we find in fact that
\[
\mathcal{S}\BDsim
\frac{1}{1 - \alpha}\Z\times
\frac{1}{1 - \alpha}\Z.
\]
One can check that the other sets are uniformly spread as well.
We derive a more precise statement from Rayleigh's theorem (Lemma~\ref{lem:Beatty}); this is used in constructing explicit examples of tile sets, see \S\ref{sec:ex1}.

\begin{quote}\bfseries
Our goal is to construct a finite set of patch-tiles composed of Bolts.
\end{quote}

By combining the three kinds of Bolt in suitable proportions, we make patch-tiles; this forces the slope of the Sturmian lattice.
Here the BD correspondence (see \S\ref{sec:BD}) plays a central role.
We write the BD-equivalence class of a patch-tile consisting of $x$ copies of Bolt $S$, $y$ copies of $M$, and $z$ copies of $L$ as $xS + yM + zL$.
The $\SL(\alpha)$-tile set $\AA = \AA(\alpha)$ that we construct has the form
\begin{equation}
\label{eq:A-form}
\AA/{\BDsim} = \left\{
\begin{aligned}
T' &= x'S + y'M + z'L,\\
T''&= x''S + y''M + z''L,
\end{aligned}~
\Nut(S),~\Nut(M),~\Nut(L)
\right\}
\end{equation}
Let us look at a worked example.

\begin{ex}
\label{ex:sqrt6}
If a tile set $\AA$ that admits a tiling has the form
\[
\AA/{\BDsim} = \{3S + 2L,~M,~\Nut(S),~\Nut(M),~\Nut(L)\}
\]
then it is an $\SL(\sqrt{6} - 2)$-tile set.
Indeed, in a tiling $\TT$ by $\AA$ the natural densities of the Bolts are computed as
\[
\delta(\mathcal{S}) :
\delta(\mathcal{L}) = 3 : 2.
\]
By the argument in \S\ref{sec:Nuts}, $\TT$ corresponds to some Sturmian lattice.
Writing its slope as $\alpha\in [0, 1]$, by \eqref{eq:Bolts-nd} it satisfies
\[
(1 - \alpha)^{2} : \alpha^{2} = 3 : 2
\]
which forces $\alpha = \sqrt{6} - 2$.
For a concrete example of a tile set, see \S\ref{sec:sqrt6}.
\end{ex}

More generally, for a given irrational slope $\alpha$, if we find $x', y', z', x'', y'', z''\in \Z_{\ge 0}$ such that
\[
\begin{pmatrix}
(1 - \alpha)^{2}\\
2\alpha(1 - \alpha)\\
\alpha^{2}
\end{pmatrix} = \delta'\begin{pmatrix}
x'\\ y'\\ z'
\end{pmatrix} + \delta''\begin{pmatrix}
x''\\ y''\\ z''
\end{pmatrix}
\]
holds for some $\delta', \delta'' > 0$, then one can construct an $\SL(\alpha)$-tile set $\AA$ with \eqref{eq:A-form}.
By Propositions~\ref{prop:existence-sub}~and~\ref{prop:add-sub}, one can respectively divide lattices $\mathcal{S}$, $\mathcal{M}$, and $\mathcal{L}$ into two uniformly spread sets
\begin{align*}
\mathcal{S}&\BDsim \mathcal{S}'\vee \mathcal{S}'',&
\mathcal{M}&\BDsim \mathcal{M}'\vee \mathcal{M}'',&
\mathcal{L}&\BDsim \mathcal{L}'\vee \mathcal{L}''
\end{align*}
with natural densities
\begin{align*}
\delta(\mathcal{S}') : \delta(\mathcal{M}') : \delta(\mathcal{L}')&=
x' : y' : z',&
\delta(\mathcal{S}'') : \delta(\mathcal{M}'') : \delta(\mathcal{L}'')&=
x'' : y'' : z''.
\end{align*}
By applying Corollary~\ref{cor:BD-correspondence}
to these lattices finite times, we finally get a set of patch-tiles $\AA$.
However, we do not yet know
\begin{itemize}
\item
whether there exist such suitable integers $x', y', \dots, z''$, or
\item
whether $\AA$ enforces slope $\alpha\notin \Q$ (if we allow $\AA$ to provide another rational slope, then we cannot conclude that $\AA$ is aperiodic).
\end{itemize}
To resolve these remaining problems, we perform the necessary calculations.
\label{rough:1}

We denote by $\R P^{2}$ the real projective plane and let \[
\Delta:= \bigl\{[x:y:z]\in \R P^{2}\bigm| x, y, z\ge 0\bigr\}.
\]We interpret each element $[x:y:z]\in \Delta$ as the equivalence class of patches consisting of Bolts $S$, $M$, and $L$ in ratio $x : y : z$, where the ``ratio'' refers to
\begin{itemize}
\item
the cardinalities if the patch is bounded (e.g., a patch-tile);
\item
the natural densities if unbounded (e.g., a tiling).
\end{itemize}
By abuse of notation, we often write these equivalence classes as $T$ or $\TT$.

For example, an (equivalence class of) patch-tile $T$ has the form
\[
T = [x : y : z]\in \Delta
\]
with $x, y, z\in \Z_{\ge 0}$ (or $x, y, z\in \Q_{\ge 0}$), which is a rational point in $\Delta$.
For a finite set
\[
\AA = \{T_{i} = [x_{i} : y_{i} : z_{i}]\in \Delta\mid i = 1, \dots, n\}
\]
of (equivalence classes of) tiles with normalization $x_{i} + y_{i} + z_{i} = 1$, the element
\[
\TA(\delta_{1}, \dots, \delta_{n}):=
\left[
\sum_{i=1}^{n}\delta_{i}x_{i} :
\sum_{i=1}^{n}\delta_{i}y_{i} :
\sum_{i=1}^{n}\delta_{i}z_{i}
\right]\in \Delta
\]
corresponds to a tiling in which $T_{i}$ has natural density $\delta_{i}\ge 0$ with $\delta_{1} + \dots + \delta_{n} = 1$.
In particular, the set of tilings realized by $\AA$ corresponds to a convex subset
\[
\TA:= \bigl\{\TA(\delta_{1}, \dots, \delta_{n})\in \Delta\bigm|
\delta_{i}\ge 0\bigr\}
\]
with rational vertices $T_{1}, \dots, T_{n}$.
The element
\[
\TTH(\alpha):= \bigl[(1 - \alpha)^{2} : 2\alpha(1 - \alpha) : \alpha^{2}\bigr]\in \Delta
\]
corresponds to Bolt tilings with respect to slope $\alpha$.

\begin{figure}[htb]\centering
\includegraphics[width = .5\linewidth]{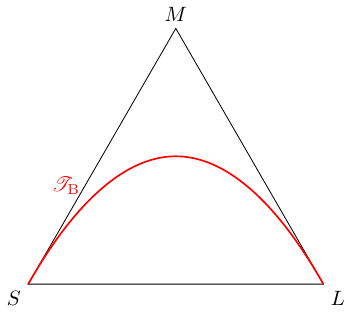}
\caption{Triangular region $\Delta$ and $\TTH$}
\label{fig:delta}
\end{figure}

Due to our matching rule, a Bolt tiling must form a Sturmian lattice with some slope $\alpha$, so a tile set $\AA$ admits a tiling only in $\TA\cap \TTH$ where
\[
\TTH:= \bigl\{\TTH(\alpha)\in \Delta\bigm| \alpha\in [0, 1]\bigr\}
\]
is the parabolic curve with end points $S = [1:0:0]$ and $L = [0:0:1]$ shown in Figure~\ref{fig:delta}.
This means that $\AA$ does not admit a tiling if $\TA\cap \TTH = \emptyset$.
Conversely, if $\TA\cap \TTH\ne \emptyset$, then for any $\TT\in \TA\cap \TTH$ one can find a finite tile set $\AT$ that admits $\TT$.
Precisely speaking, for each $i\in \{1, \dots, n\}$, we can take finitely many patch-tiles $T_{i}^{1}, \dots, T_{i}^{J(i)}$ that contain Bolts in the given ratio 
$x_{i} : y_{i} : z_{i}$.

To show this, let
\[
\left[
\sum_{i=1}^{n}\delta_{i}x_{i} :
\sum_{i=1}^{n}\delta_{i}y_{i} :
\sum_{i=1}^{n}\delta_{i}z_{i}
\right] =
\bigl[(1 - \alpha)^{2} : 2\alpha(1 - \alpha) : \alpha^{2}\bigr]\in \TA\cap \TTH.
\]
One can then take a Delone subset $\mathcal{S}_{1}\subset \mathcal{S}$ with natural density $\delta(\mathcal{S}_{1}) = \delta_{1}x_{1}$ such that each $\mathcal{S}_{1}$ and $\mathcal{S}\setminus \mathcal{S}_{1}$ is uniformly spread by Propositions~\ref{prop:existence-sub}~and~\ref{prop:add-sub}.
A similar discussion allows us to take $\mathcal{S}_{2}\subset \mathcal{S}\setminus \mathcal{S}_{1}$ with $\delta(\mathcal{S}_{2}) = \delta_{2}x_{2}$ such that each $\mathcal{S}_{2}$ and $\mathcal{S}\setminus (\mathcal{S}_{1}\vee \mathcal{S}_{2})$ is uniformly spread.
By iterating this operation, we finally get a partition $\mathcal{S} = \bigsqcup_{i}\mathcal{S}_{i}$ with $\delta(\mathcal{S}_{i}) = \delta_{i}x_{i}$ and each $\mathcal{S}_{i}$ is uniformly spread.
The same result holds for the other Bolts, i.e., there are partitions $\mathcal{M} = \bigsqcup_{i}\mathcal{M}_{i}$ and $\mathcal{L} = \bigsqcup_{i}\mathcal{L}_{i}$ by uniformly spread Delone sets with natural densities
\begin{align*}
\delta(\mathcal{S}_{i})&= \delta_{i}x_{i},&
\delta(\mathcal{M}_{i})&= \delta_{i}y_{i},&
\delta(\mathcal{L}_{i})&= \delta_{i}z_{i}.
\end{align*}
By applying Corollary~\ref{cor:BD-correspondence} for Bolt sets $\mathcal{S}_{i}$, $\mathcal{M}_{i}$, and $\mathcal{L}_{i}$, one can take partitions
\begin{align*}
\mathcal{S}_{i}&= \bigsqcup_{j\ge 1}\mathcal{S}_{i}^{j},&
\mathcal{M}_{i}&= \bigsqcup_{j\ge 1}\mathcal{M}_{i}^{j},&
\mathcal{L}_{i}&= \bigsqcup_{j\ge 1}\mathcal{L}_{i}^{j}
\end{align*}
by finite sets which satisfy that \[
\Card(\mathcal{S}_{i}^{j}) :
\Card(\mathcal{M}_{i}^{j}) :
\Card(\mathcal{L}_{i}^{j}) =
x_{i} : y_{i} : z_{i}
\]for each $j$, and that patch-tile $T_{i}^{j} = \bigcup(\mathcal{S}_{i}^{j}\cup \mathcal{M}_{i}^{j}\cup \mathcal{L}_{i}^{j})$, the union of finitely many Bolts, is uniformly bounded with respect to $j$.
Since we have FLC for the resulting tiling by our matching rule, $\{T_{i}^{j}\}_{j\ge 1}$ is finite up to isometry and we may denote it by $\{T_{i}^{j}\}_{j = 1, \dots, J(i)}$.
Finally we get a finite tile set\[
\AT = \bigcup_{i=1}^{n}
\{T_{i}^{j}\mid j = 1, \dots, J(i)\}
\]that admits a tiling $\TT$.
If $\TA\cap \TTH = \{\TTH(\alpha)\}$ and $\alpha$ is irrational, then $\AA_{\TTH(\alpha)}$ must be aperiodic since $\SL(\alpha)$ has no periods.

In fact, we give aperiodic tile sets that form a Sturmian lattice with 
arbitrary quadratic irrational slope $\alpha$:

\begin{thm}
\label{thm:Voronoi-aperiodic}
Let $\alpha\in [0, 1]$.
\begin{enumerate}
\item
If $\alpha$ is rational, then there exists unique $T\in \Delta$ such that
\[
\AA/{\BDsim} = \{T, \Nut(S), \Nut(M), \Nut(L)\}
\]
enforces Sturmian lattices with only rational slope $\alpha$.
\item
If $\alpha$ is quadratic irrational, then there exists $T', T''\in \Delta$ such that
\[
\AA/{\BDsim} = \{T', T'', \Nut(S), \Nut(M), \Nut(L)\}
\]
enforces Sturmian lattices with only irrational slope $\alpha$ (or its Galois conjugate).
In particular, there exists an aperiodic tile 
set $\AA(\alpha):= \AA_{\TTH(\alpha)}$ whose underlying expansion constant is a quadratic unit
of any real quadratic field $\Q(\alpha)$.
\item
Otherwise, there do not exist tile sets $\AA$ that enforces a Sturmian lattice of slope $\alpha$.
\end{enumerate}
\end{thm}

We call such a tile set $\AA(\alpha)$ in Theorem~\ref{thm:Voronoi-aperiodic}(2) an \emph{$\SL(\alpha)$-tile set}.
We give a few remarks.

\begin{rem}
\label{rem:upper-bound}
The case (1) is obvious since there is a periodic lattice with rational slope $\alpha$ by Theorem~\ref{thm:rational}.
We may take a fundamental domain of the Sturmian lattice and this gives a periodic tiling.

In the case (2), note that this result does not give an upper bound on the number of patch-tiles since $T'$ and $T''$ are equivalent classes of patch-tiles.
The actual cardinality of the patch-tile set 
depends on the BD correspondence used 
for construction.
This also does not exclude the possibility of aperiodic monotile, i.e., the number of patch-tiles is one.

The case (3) says that if there exists $\AA$ that admits 
a Sturmian lattice with such a slope $\alpha$, 
then it also admits another Sturmian lattice with rational slope.
For example, 
$\AA = \{S, M, L, \Nut(S), \Nut(M), \Nut(L)\}$ trivially forms Sturmian lattices with any slope $\alpha\in [0, 1]$.
The problem to construct aperiodic tile sets with a 
cubic or higher dimensional irrational expansion factor is challenging but 
it is beyond the scope of this method, although Sturmian lattices exist
for any irrational slope.
\end{rem}

\begin{figure}[htb]\centering
\subfigure[$\AA = \{T\}$]{%
\includegraphics[width = .3\linewidth]{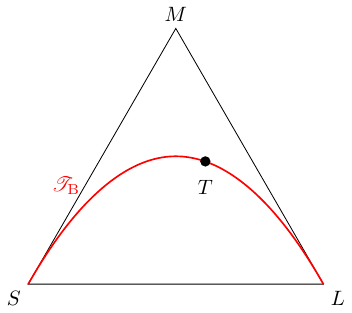}
}\quad%
\subfigure[$\AA = \{T_{1}, T_{2}\}$]{%
\includegraphics[width = .3\linewidth]{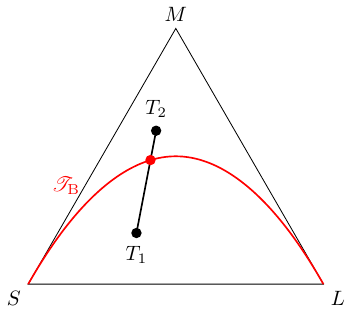}
}\quad%
\subfigure[$\AA = \{T_{1}, T_{2}, T_{3}\}$]{%
\includegraphics[width = .3\linewidth]{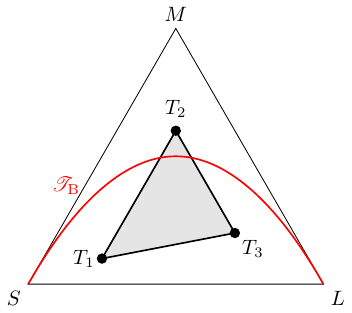}
}
\caption{Examples for $\TA$ with cardinality $1$, $2$, and $3$ of $\AA$}
\label{fig:Voronoi-TA}
\end{figure}

\begin{proof}[Proof of  Theorem~\ref{thm:Voronoi-aperiodic}(3).]
We observe a patch-tile set $\AA$ that satisfies $\TA\cap \TTH\ne \emptyset$.
Recall that $\TTH$ has rational end points $S$ and $L$.
There are three cases (see Figure~\ref{fig:Voronoi-TA}).
\begin{enumerate}
\renewcommand{\labelenumi}{(\alph{enumi})}
\item
If a vertex $T$ of $\TA$ is on $\TTH$, then $T = \TTH(\alpha)$ is rational and so is $\alpha$.
\item
If an edge of $\TA$ intersects with $\TTH$, then the intersection $T = \TTH(\alpha)$ is rational or quadratic and so is $\alpha$.
\item
If $\TA$ is a convex polygon and $\TTH$ passes through its interior, then $\TA$ cuts out an arc of $\TTH$ with positive length.
In particular, there are infinitely many rational points $\TTH(q)$ with $q\in \Q$ in the arc.
\end{enumerate}
These observations yield the case (3).
\end{proof}

For (2), the situation is reduced to the case~(b).
We provide a specific algorithm for determining coefficients $x', y', z', x'', y'', z''\in \Z_{\ge0}$ that
\[
\AA(\alpha) = \bigl\{T' = [x' : y' : z'],~T'' = [x'' : y'' : z'']\bigr\}
\]
is an $\SL(\alpha)$-tile set for given quadratic $\alpha\in (0, 1)$.

For $p_{i} = [x_{i} : y_{i} : z_{i}]\in \R P^{2}$ ($i = 1, 2$), $p_{1}$ is \emph{orthogonal} to $p_{2}$ if
\[
x_{1}x_{2} + y_{1}y_{2} + z_{1}z_{2} = 0
\]
holds, and we denote by $p_{1}\perp p_{2}$.
%
We give a criterion for \emph{$\SL(\alpha)$-tiles}, the elements of $\AA(\alpha)$.

\begin{prop}
\label{prop:Voronoi-formula}
Let $\alpha\in (0, 1)$ be a quadratic irrational number with the minimal monic polynomial
\[
f(\alpha):= \alpha^{2} + u\alpha + v = 0
\]
over $\Q$, and let
\[
\bm{n}(\alpha):= \left[
v : \left(\frac{u}{2} + v\right) : (u + v + 1)
\right]\in \Q P^{2}.
\]
Then for
\[
\AA = \bigl\{T' = [x' : y' : z'],~T'' = [x'' : y'' : z'']\bigr\}
\]
with $x_{i}, y_{i}, z_{i}\in \Q_{\ge 0}$, the following are equivalent as long as $\TA\cap \TTH\ne \emptyset$.
\begin{enumerate}
\item
$\TTH(\alpha)\in \TA\cap \TTH$.
\item
Both $T'$ and $T''$ are orthogonal to $\bm{n}(\alpha)$.
\end{enumerate}
\end{prop}

\begin{proof}
First notice that one can characterize a vector $\bm{n}(\alpha)$ as follows:
\begin{equation}
\label{eq:n-characterization}
\bm{n}(\alpha)\in \Q P^{2}\text{ orthogonal to }\TTH(\alpha) = \bigl[
(1 - \alpha)^{2} :
2\alpha(1 - \alpha) :
\alpha^{2}
\bigr].
\end{equation}
Indeed, supposing that $\bm{n} = (n_{1}, n_{2}, n_{3})\in \Q P^{2}$ is orthogonal to $\TTH(\alpha)$:
\[
n_{1}(1 - \alpha)^{2} +
n_{2}\cdot 2\alpha(1 - \alpha) +
n_{3}\alpha^{2} = 0,
\]
we have
\[
(n_{1} - 2n_{2} + n_{3})\alpha^{2} +
2(n_{2} - n_{1})\alpha + n_{1} = 0.
\]
Our assumption $\alpha\notin \Q$ implies that $n_{1} - 2n_{2} + n_{3}\ne 0$, so we may assume that $n_{1} - 2n_{2} + n_{3} = 1$ without loss of generality.
By the uniqueness of minimal polynomial, we obtain that
\begin{align*}
u&= 2(n_{2} - n_{1}),&
v&= n_{1},
\end{align*}
and hence $n_{1} = v$, $n_{2} = u/2 + v$ and $n_{3} = 2n_{2} - n_{1} + 1 = u + v + 1$ as required.

For $\AA = \{T', T''\}$, we define the subspace $\Pi(\AA) = \operatorname{span}\{T_{1}, T_{2}\}$ of $\R P^{2}$.
Then one can see that (1) and (2) are respectively equivalent to
\begin{enumerate}
\renewcommand{\labelenumi}{(\arabic{enumi})$'$}
\item
$\TTH(\alpha)\in \Pi(\AA)$, and
\item
$\bm{n}(\alpha)\in \Pi(\AA)^{\perp}$.
\end{enumerate}
$(2)'\Rightarrow (1)'$ is clear since $\TTH(\alpha)$ is orthogonal to $\bm{n}(\alpha)$.
Conversely, suppose $(1)'$.
Then both $\bm{n}(\alpha)$ and $\bm{n}:= T_{1}\times T_{2}$ in $\Q P^{2}$ are orthogonal to $\TTH(\alpha)$.
By characterization \eqref{eq:n-characterization} of $\bm{n}(\alpha)$, they must coincide up to scalar multiplication.
It follows that $\bm{n}(\alpha) = \bm{n}\in \Pi(\AA)^{\perp}$ as required.
\end{proof}

We naturally identify hexagonal cells $S$, $M$, and $L$ with elements
\begin{align*}
S&= [1:0:0],&
M&= [0:1:0],&
L&= [0:0:1]
\end{align*}
in $\Delta$, and write
\[
xS + yM + zL:= [x : y : z]\in \Delta.
\]
We define map $\varphi\colon \Delta\to \R$ by
\begin{align*}
\varphi(xS + yM + zL)&:=
xv + y(u/2 + v) + z(u + v + 1)\\&=
[x : y : z]\cdot \bm{n}(\alpha).
\end{align*}
Even though $\varphi$ is not well-defined, it does not matter since we only focus on whether the value of $\varphi$ is zero or not.
For convenience, we also write
\[
\varphi(S) : \varphi(M) : \varphi(L) =
v : (u/2 + v) : (u + v + 1).
\]

Define
\[
\AA(\alpha):= \{xS + yM + zL\in \Delta\mid
x, y, z\in \Z_{\ge 0}\text{ and }\varphi(xS + yM + zL) = 0\}.
\]
We see that $\AA(\alpha)$ coincides with the set of $\SL(\alpha)$-tiles by $\bm{n}(\alpha)\perp \AA(\alpha)$.

We finish the discussion for Theorem~\ref{thm:Voronoi-aperiodic}(2), we will see the strength of Proposition~\ref{prop:Voronoi-formula} to obtain $\AA$.

\begin{proof}[Proof of  Theorem~\ref{thm:Voronoi-aperiodic}(2).]
Let $\alpha\in (0, 1)$ be quadratic irrational.
We claim that at least two of $\varphi(S)$, $\varphi(M)$, $\varphi(L)$ have opposite signs.
Assume that all $\varphi(S)$, $\varphi(M)$, $\varphi(L)$ are have the same signs.
By Proposition~\ref{prop:Voronoi-formula}, vector $\bm{n}(\alpha)$ is orthogonal to $\TTH(\alpha) = [\delta_{S} : \delta_{M} : \delta_{L}]$, i.e.,
\[
\varphi(S)\delta_{S} + \varphi(M)\delta_{M} + \varphi(L)\delta_{L} = 0.
\]
Since $\alpha\in (0, 1)$, we have $\delta_{S}, \delta_{M}, \delta_{L} > 0$, and hence $\varphi(S) = \varphi(M) = \varphi(L) = 0$ is the only possibility.
However, since $f$ does not have rational roots, we obtain $\varphi(S) = v = f(0)\ne 0$ and $\varphi(L) = u + v + 1 = f(1)\ne 0$, which is a contradiction.

Thus there exists $A, B\in \{S, M, L\}$ with $a = \varphi(A) > 0 > \varphi(B) = -b$.
Take a patch-tile $T' = bA + aB$, which satisfies that 
\[
\varphi(T') = b\varphi(A) + a\varphi(B) = 0.
\]
Let $C\in \{S, M, L\}\setminus \{A, B\}$ be the other cell.
We define patch-tile $T''$ by
\[
T'' = \begin{cases}
C & \text{if }\varphi(C) = 0\\
cB + bC & \text{if }\varphi(C) = c > 0\\
cA + aC & \text{if }\varphi(C) =-c < 0
\end{cases}.
\]
By construction, both $T'$ and $T''$ are $\SL(\alpha)$-tiles.
We claim that $\AA = \{T', T''\}$ admits a tiling of quadratic slope $\alpha$.
By changing the labels, we may assume that
\begin{equation}
\label{eq:thm:formula}
\begin{aligned}
T'&= x'A + y'B,&
T''&= x''A + z''C
\end{aligned}
\end{equation}
with $x', y', z''\in \Q_{> 0}$, $x''\in \Q_{\ge 0}$, and $x' + y' = x'' + z'' = 1$.
Then
\[
\TTH(\alpha) = [\delta_{A} : \delta_{B} : \delta_{C}] =
\left[
* : \frac{\delta_{B}}{y'}y' : \frac{\delta_{C}}{z''}z''
\right] =
\TA\left(\frac{\delta_{B}}{y'}, \frac{\delta_{C}}{z''}\right)\in \TA
\]
as claimed.
\end{proof}

\section{Examples I}
\label{sec:ex1}

We apply our method in the previous section to several specific slopes.
Table~\ref{tab:number} summarizes the number of patch-tiles.
Note that there is redundancy in our tile enumeration; we should be able to eliminate more unnecessary tiles.

More examples are given in \S\ref{sec:ex2}.
It is worth comparing the two examples presented in \S\ref{sec:sqrt2} and \S\ref{sec:sqrt2-2}.
Despite using the same slope, the latter example succeeds in reducing the cardinality to three.
Similar optimization may be possible with other quadratic slopes.
Moreover, \S\ref{sec:exTurtle} provides an interpretation of Smith Turtle \cite{SMKGS:23_1} in terms of our theory.
This raises the possibility that a new aperiodic monotile could be found in this way.

\begin{table}[htb]\centering
\caption{Number of tiles in \S\ref{sec:ex1} and \S\ref{sec:ex2}, up to isometry
}
\label{tab:number}
\begin{tabular}{c|c|crcrr|r}
See & Slope $\alpha$ &
\multicolumn{2}{c}{Type I} &
\multicolumn{2}{c}{Type II} & 
Nuts & Total \\ \hline
\S\ref{sec:sqrt2} & $\sqrt{2} - 1$ &
$2S + \phantom{1}L$ & $6$ &
$M$ & $1$ & $3$ & $10$ \\
\S\ref{sec:sqrt6} & $\sqrt{6} - 2$ &
$3S + 2L$ & $25$ &
$M$ & $1$ & $3$ & $29$ \\
\S\ref{sec:sqrt7} & $\sqrt{7} - 2$ &
$2S + 3L$ & $20$ &
$2M + L$ & $1$ & $3$ & $24$ \\ \hline
\S\ref{sec:sqrt2-2} & $\sqrt{2} - 1$ &
$12S + 6M + 6L$ & $2$ &
$M$ & $1$ & $-$ & $3$ \\
\S\ref{sec:exTurtle}, \cite{SMKGS:23_1} &
$\frac{5 - \sqrt{5}}{10}$ &
$6S + 2M$ & $1$ &
$2M + 6L$ & $1$ & $-$ & $2$ \\ \hline
\end{tabular}
\end{table}

In particular, we have

\begin{thm}
\label{thm:lambda}
Let $\lambda > 1$ be a quadratic integer with height $h\ge 3$.
Then there exists an absolute constant $C > 0$ and an aperiodic tile set $\AA$ for tilings with expansion constant $\lambda$ with
\[
\Card(\AA)\le 2h + C.
\]
\end{thm}

See \S\ref{sec:lambda} for the proof when $N(\lambda) = -1$.
The following lemma gives a representation of Bolt sets $\mathcal{S}$, $\mathcal{M}$, and $\mathcal{L}$ in terms of lattices.

\begin{lem}
\label{lem:Beatty}
For fixed $\alpha\in (0, 1)\setminus \Q$ and $\rho\in \R_{\alpha}$, the two sets
\begin{align*}
D_{0}&= \left\{
\round{\frac{m + \rho}{1 - \alpha}}\in \Z + \frac{1}{2}
~\middle|~
m\in \Z
\right\},&
D_{1}&= \left\{
\round{\frac{m - \rho}{\alpha}}\in \Z + \frac{1}{2}
~\middle|~
m\in \Z
\right\}
\end{align*}
give a partition $\Z + 1/2 = D_{0}\sqcup D_{1}$.
Moreover, the bi-infinite word $(w_{n})$ defined by
\[
w_{n} = \begin{cases}
0 & \text{if }n + 1/2\in D_{0},\\
1 & \text{if }n + 1/2\in D_{1}
\end{cases}
\]
coincides with the Sturmian word of slope $\alpha$ and intercept $\rho$.
\end{lem}

\begin{proof}
This is essentially a half-integer reformulation of the classical complementary Beatty sequence theorem (Rayleigh's theorem).
For a particularly transparent geometric proof, see Tamura~\cite{Tamura1996}.
The present statement is obtained by translating that partition to $\Z+1/2$ and reading the corresponding symbolic coding.
\end{proof}

By this lemma, we obtain the following representation.
\begin{equation}
\label{eq:Beatty}
\begin{aligned}
\mathcal{S}&\BDsim
\hat{\mathcal{S}} =
\frac{\Z + \rho_{1}}{1 - \alpha}\times
\frac{\Z + \rho_{2}}{1 - \alpha},\\
\mathcal{M}&\BDsim
\hat{\mathcal{M}} = \left(
\frac{\Z - \rho_{1}}{\alpha}\times
\frac{\Z + \rho_{2}}{1 - \alpha}
\right)\vee \left(
\frac{\Z + \rho_{1}}{1 - \alpha}\times
\frac{\Z - \rho_{2}}{\alpha}
\right),\\
\mathcal{L}&\BDsim
\hat{\mathcal{L}} =
\frac{\Z - \rho_{1}}{\alpha}\times
\frac{\Z - \rho_{2}}{\alpha}.
\end{aligned}
\end{equation}
Here, for $\lambda, \mu > 0$ and $\eta_{1}, \eta_{2}\in \R$, we set
\[
\frac{\Z + \eta_{1}}{\lambda}\times
\frac{\Z + \eta_{2}}{\mu}:=
\left\{
\Triplet{*}
{~\frac{j + \eta_{1}}{\lambda}}
{~\frac{k + \eta_{2}}{\mu}}\in \R^{2}
~\middle|~
(j, k)\in \Z^{2}
\right\}.
\]
\color{red!50!black}
There is a trivial BD bijection $\psi\colon \hat{\mathcal{S}}\vee \hat{\mathcal{M}}\vee \hat{\mathcal{L}}\to \mathcal{S}\vee \mathcal{M}\vee \mathcal{L} = (\Z + 1/2)^{2}$, defined by $\psi(x, y) = (\round{x}, \round{y})$.
Our basic strategy is to apply cross-BD pairs (see \S\ref{sec:cross}) between three lattices $\hat{\mathcal{S}}$, $\hat{\mathcal{M}}$, and $\hat{\mathcal{L}}$.
\color{black}

\begin{rem}
\label{rem:Beatty}
As we saw in \S\ref{sec:BDeq}, the BD-equivalence relation on the class of uniformly spread Delone sets is powerful and has natural density as a complete invariant.
In particular, $\mathcal{M}$ also has the representation
\[
\mathcal{M}\BDsim
\frac{1}{\sqrt{2\alpha(1 - \alpha)}}
(\Z\times \Z),
\]
which is more concise than \eqref{eq:Beatty}.
{\color{HMDcomment}In fact, this representation alone suffices in theory to guarantee the existence of $\SL(\alpha)$-tile sets.}

{\color{HMDcomment}Nevertheless, the representation \eqref{eq:Beatty} is useful for the explicit construction of $\SL(\alpha)$-tile sets: Lemma~\ref{lem:Beatty} provides a way to fully recover the Bolt tiling for a symbolic Sturmian lattice $\SL(\alpha\mid \bm{\rho})$ from the lattices on the right-hand side of \eqref{eq:Beatty}.}
This makes the construction and description of the BD correspondence concise, and the counting of tiles easy.
For the necessary arguments, see the examples in \S\ref{sec:ex1}.
\end{rem}

\subsection{\texorpdfstring
{Example for $\alpha = \sqrt{2} - 1$}
{Example for \textbackslash alpha = \textbackslash sqrt\{2\} - 1}}
\label{sec:sqrt2}

\begin{figure}[htb]\centering
\includegraphics[width = \linewidth]
{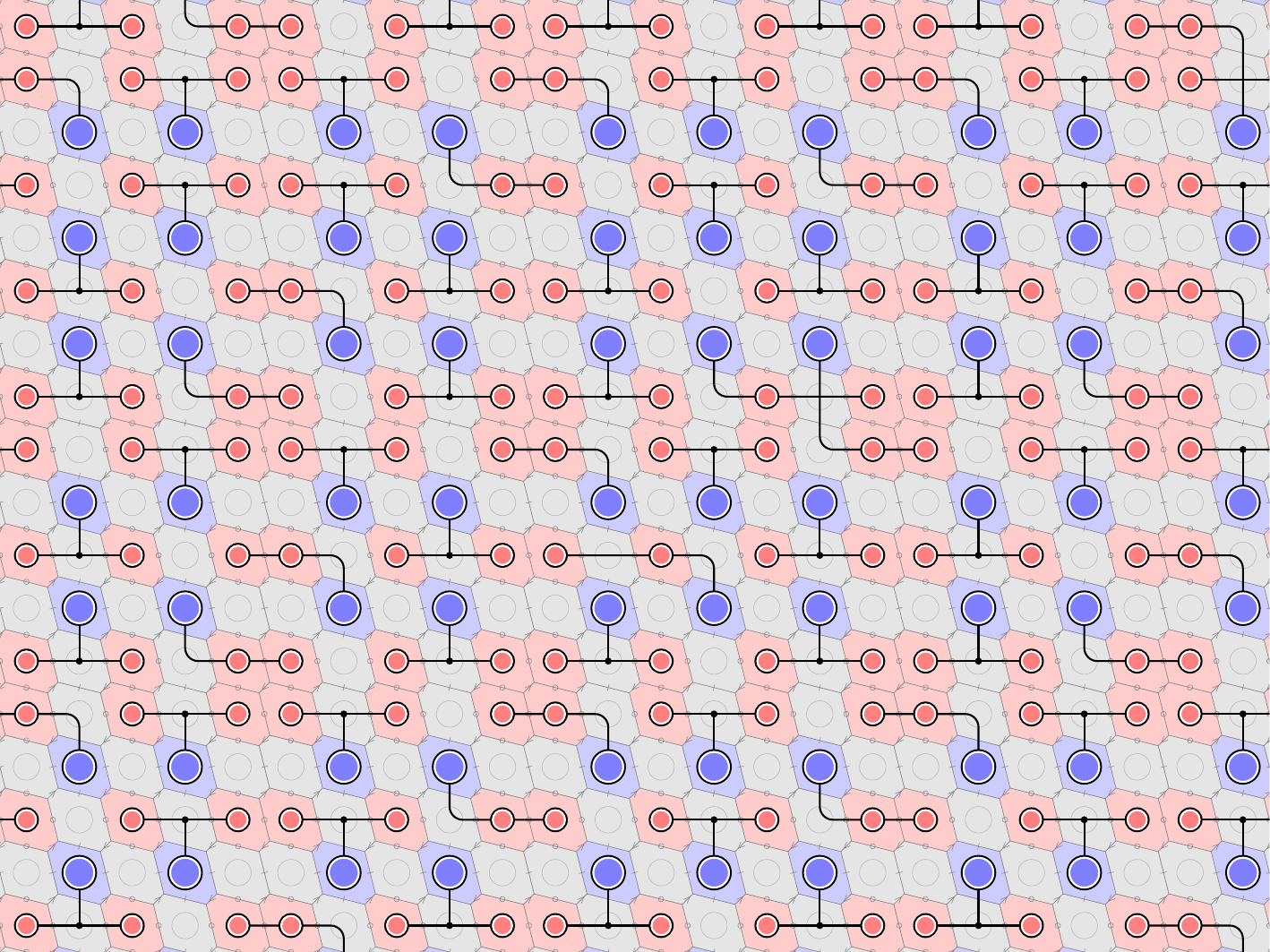}
\caption{Tiling by the $\SL(\sqrt{2} - 1)$-tile set in Figure~\ref{fig:sqrt2-tiles}}
\label{fig:sqrt2-tiling}
\end{figure}

\begin{figure}[htb]\centering
\includegraphics[width = \linewidth,
pagebox = artbox, page = 1]
{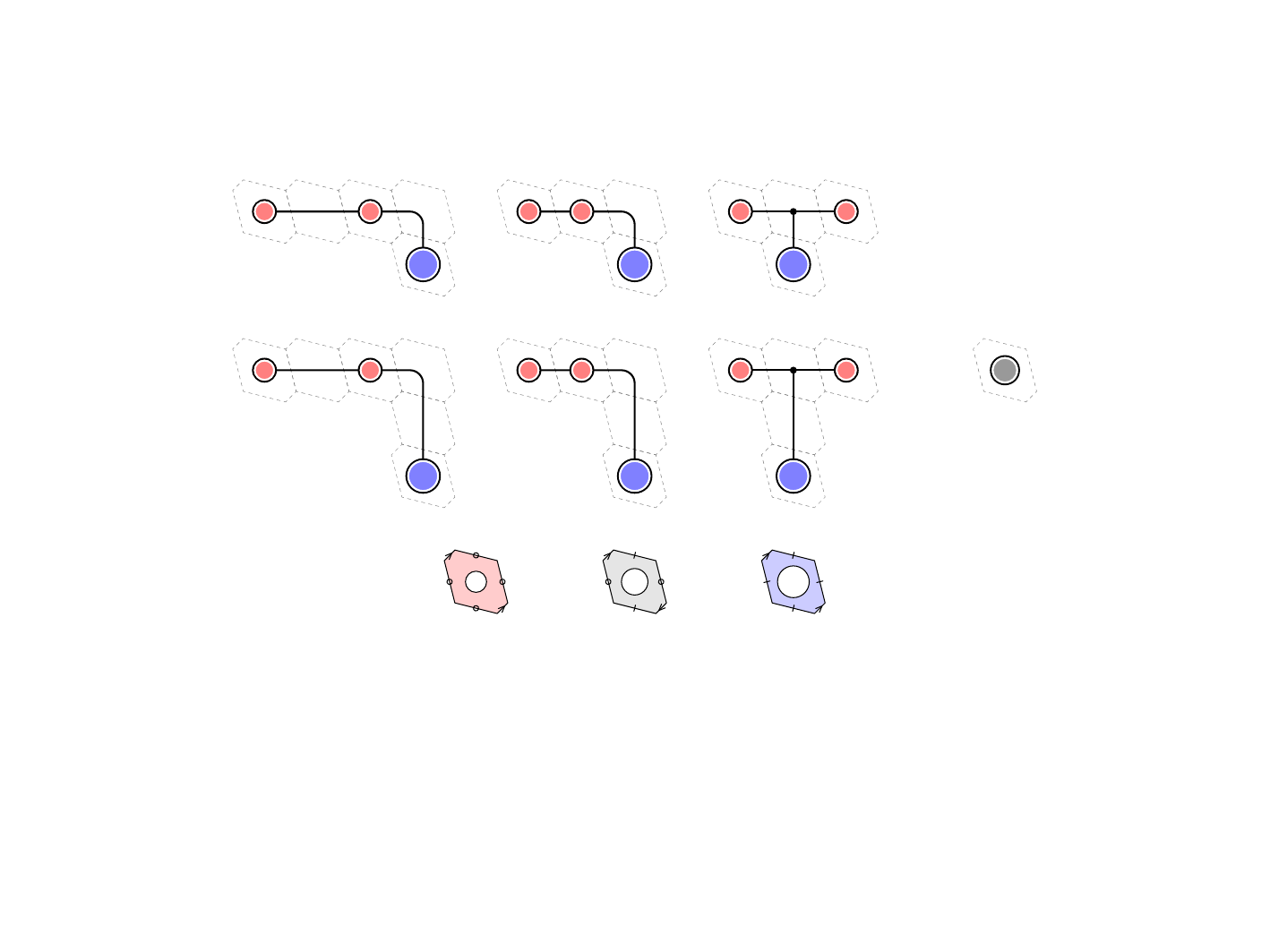}
\caption{$\SL(\sqrt{2} - 1)$-tile set composed of $10$ patch-tiles}
\label{fig:sqrt2-tiles}
\end{figure}

We begin with slope $\alpha = \sqrt{2} - 1$.
Since $\alpha^{2} + 2\alpha - 1 = 0$ ($u = 2$, $v = -1$), we have
\[
\varphi(S) : \varphi(M) : \varphi(L) =
-1 : 0 : 2,
\]
so that we may construct a tile set of type
\[
\AA/{\BDsim} =
\{2S + L,~M,~\Nut(S),~\Nut(M),~\Nut(L)\}.
\]
Figure~\ref{fig:sqrt2-tiling} shows an example of a tiling.
Here, by applying an affine transformation, we have placed each Bolt center $B(j,k)\in \R^{2}$ at the Cartesian coordinates $(j, k)\in \Z^{2}$.
The lines connecting Nuts $S$ and $L$ are meant to make the BD correspondence $\mathcal{S}\stackrel{2}{\to} \mathcal{L}$ explicit; the patch-tile $2S + L$ obtained in this way is a disjoint union of three disks.
Figure~\ref{fig:sqrt2-tiles} depicts a set of $10$ patch-tiles obtained from this construction.
We only give the result and will not go into detail; see also \S\ref{sec:sqrt6}.

\subsection{\texorpdfstring
{Example for $\alpha = \sqrt{6} - 2$}
{Example for \textbackslash alpha = \textbackslash sqrt\{6\} - 2}}
\label{sec:sqrt6}

\begin{figure}[htb]\centering
\includegraphics[width = \linewidth]
{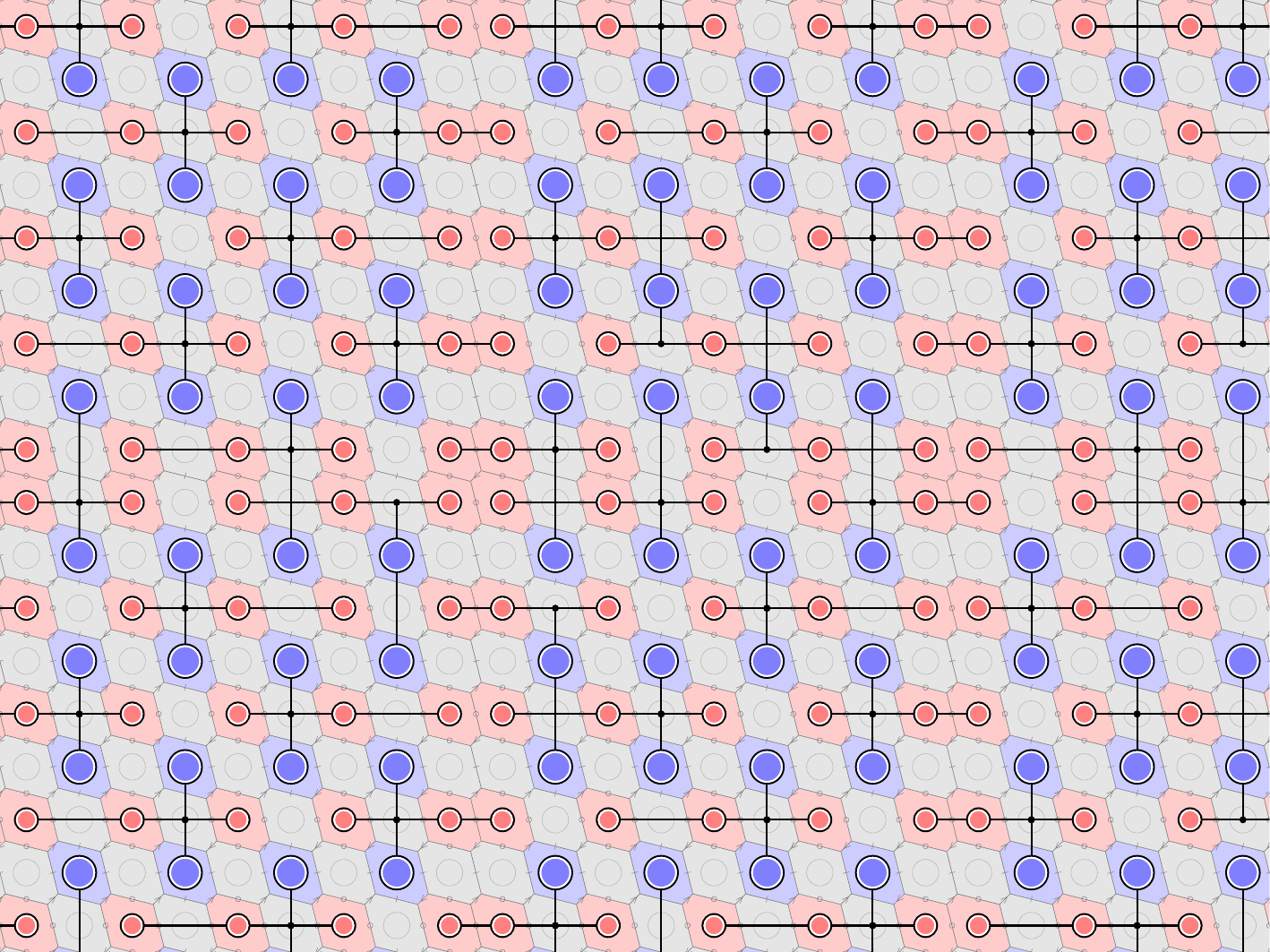}
\caption{Tiling by the $\SL(\sqrt{6} - 2)$-tile set in Figure~\ref{fig:sqrt6-tiles}}
\label{fig:sqrt6-tiling}
\end{figure}

We next apply our method to the case $\alpha=\sqrt6-2$.
Since $\alpha^{2} + 4\alpha - 2 = 0$ ($u = 4$, $v = -2$), we have
\[
\varphi(S) : \varphi(M) : \varphi(L) =
-2 : 0 : 3.
\]
We need a tile set of type
\[
\AA/{\BDsim} =
\{3S + 2L,~M,~\Nut(S),~\Nut(M),~\Nut(L)\}.
\]
Figure~\ref{fig:sqrt6-tiling} shows an example of the resulting tiling.
Accordingly, we construct a cross-BD pair $\mathcal{S}\stackrel{3}{\to}\Lambda\stackrel{2}{\gets}\mathcal{L}$ (see \S\ref{sec:cross}).
Here, applying a cross-BD pair to the Bolt sets means applying it to the lattices on the right-hand side of \eqref{eq:Beatty},
\begin{align*}
\hat{\mathcal{S}}&=
\frac{\Z + \rho_{1}}{1 - \alpha}\times
\frac{\Z + \rho_{2}}{1 - \alpha},&
\hat{\mathcal{L}}&=
\frac{\Z - \rho_{1}}{\alpha}\times
\frac{\Z - \rho_{2}}{\alpha}.
\end{align*}
{\color{HMDcomment}Note that when $\varphi(M) = 0$, the Bolt $M$ is not related to the BD correspondence we need; i.e., it is so versatile that it can be used as a single tile, or be put anywhere to shape other tiles.}
We obtain the desired BD correspondence
\[
\mathcal{S}\stackrel{\psi}{\gets}
\hat{\mathcal{S}}\stackrel{3}{\to}
\Lambda\stackrel{2}{\gets}
\hat{\mathcal{L}}\stackrel{\psi}{\to}
\mathcal{L}.
\]
A comparison of Figure~\ref{fig:sqrt6-BD} and Figure~\ref{fig:sqrt6-tiling} will visually aid understanding.

\begin{figure}[htb]\centering
\includegraphics[width = \linewidth,
pagebox = artbox, page = 1]
{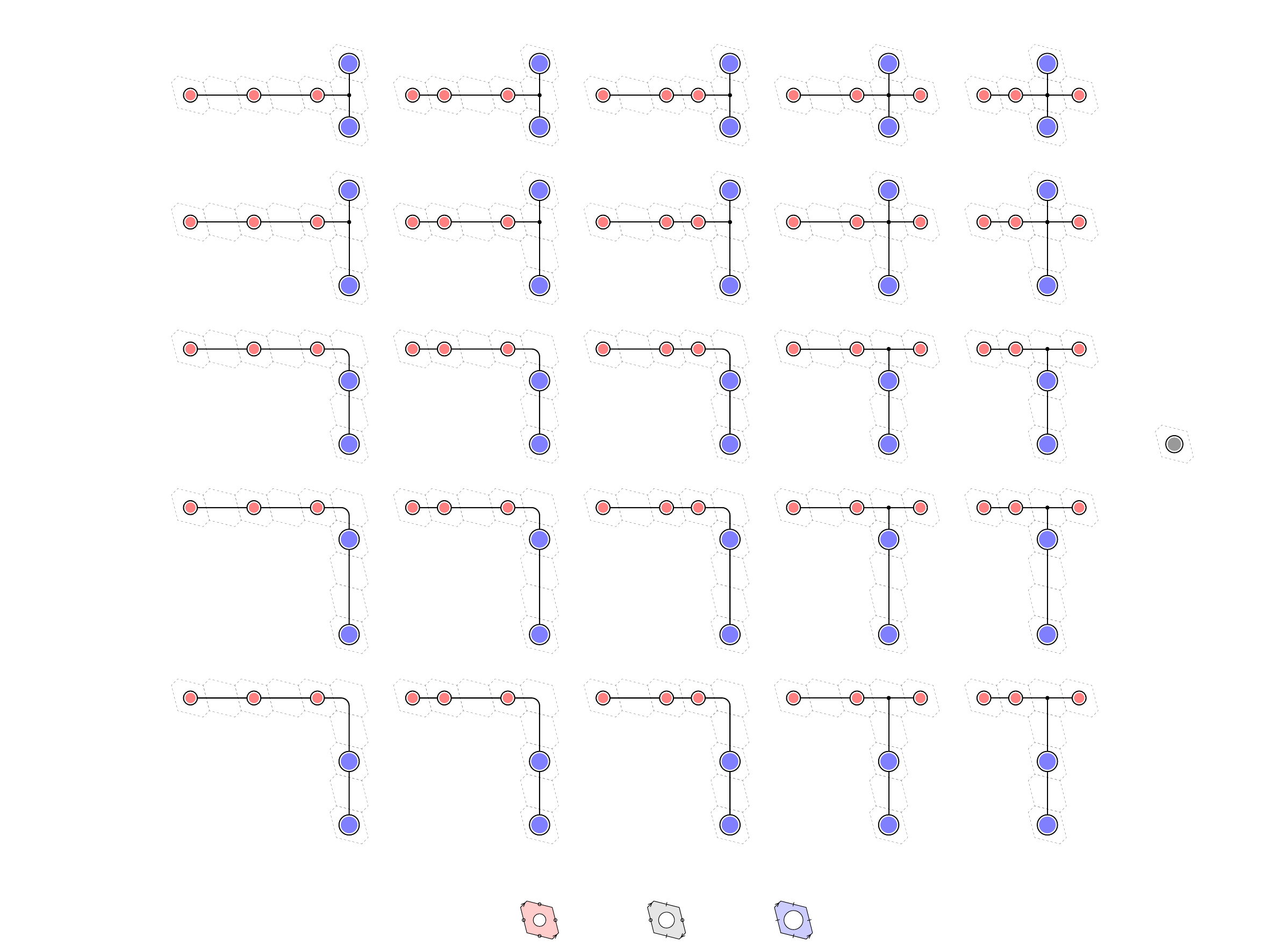}
\caption{$\SL(\sqrt{6} - 2)$-tile set composed of $29$ patch-tiles}
\label{fig:sqrt6-tiles}
\end{figure}

\subsection{\texorpdfstring
{Example for $\alpha = \sqrt{7} - 2$}
{Example for \textbackslash alpha = \textbackslash sqrt\{7\} - 2}}
\label{sec:sqrt7}

We move on the case where $\varphi(M)\ne 0$, specifically when $\alpha = \sqrt{7} - 2$.
In this situation, 
we must also apply a BD correspondence to $\mathcal{M} = \mathcal{M}_{1}\vee \mathcal{M}_{2}$, which it is natural to view as the union of the two lattices
\begin{align*}
\mathcal{M}_{1}&\BDsim
\hat{\mathcal{M}}_{1}:=
\frac{1}{\alpha}\Z\times
\frac{1}{1 - \alpha}\Z,&
\mathcal{M}_{2}&\BDsim
\hat{\mathcal{M}}_{2}:=
\frac{1}{1 - \alpha}\Z\times
\frac{1}{\alpha}\Z,
\end{align*}
recall Remark~\ref{rem:Beatty}.
Hence, instead of the Bolt $B(j,k)\in \mathcal{S}\cup \mathcal{L}$, for a sufficiently small $\eps > 0$ we use two kinds of Bolt centered respectively at
\begin{align*}
B_{1}(j,k)&:= (j + \eps, k - \eps),&
B_{2}(j,k)&:= (j - \eps, k + \eps).
\end{align*}
By this, the set $\mathcal{S}$ is duplicated into the two sets
\begin{align*}
\mathcal{S}_{1}&:=
\{B_{1}(j,k)\in \R^{2}\mid
B(j,k)\in \mathcal{S}\},\\
\mathcal{S}_{2}&:=
\{B_{2}(j,k)\in \R^{2}\mid
B(j,k)\in \mathcal{S}\}
\end{align*}
The set $\mathcal{L}$ is treated in the same way.
From now on we write the elements of $\mathcal{S}_{1}\sqcup \mathcal{S}_{2}$ and $\mathcal{L}_{1}\sqcup \mathcal{L}_{2}$ as $\frac{S}{2}$ and $\frac{L}{2}$ respectively.
Accordingly, we replace both $\Nut(S)$ and $\Nut(L)$ with a two-hole Nut.
We apply BD correspondences to $\mathcal{S}_{j}$, $\mathcal{M}_{j}$, and $\mathcal{L}_{j}$ for each $j$ to form desired patch-tiles.
It suffices to consider only $j = 1$.

\begin{figure}[htb]\centering
\includegraphics[width = \linewidth]
{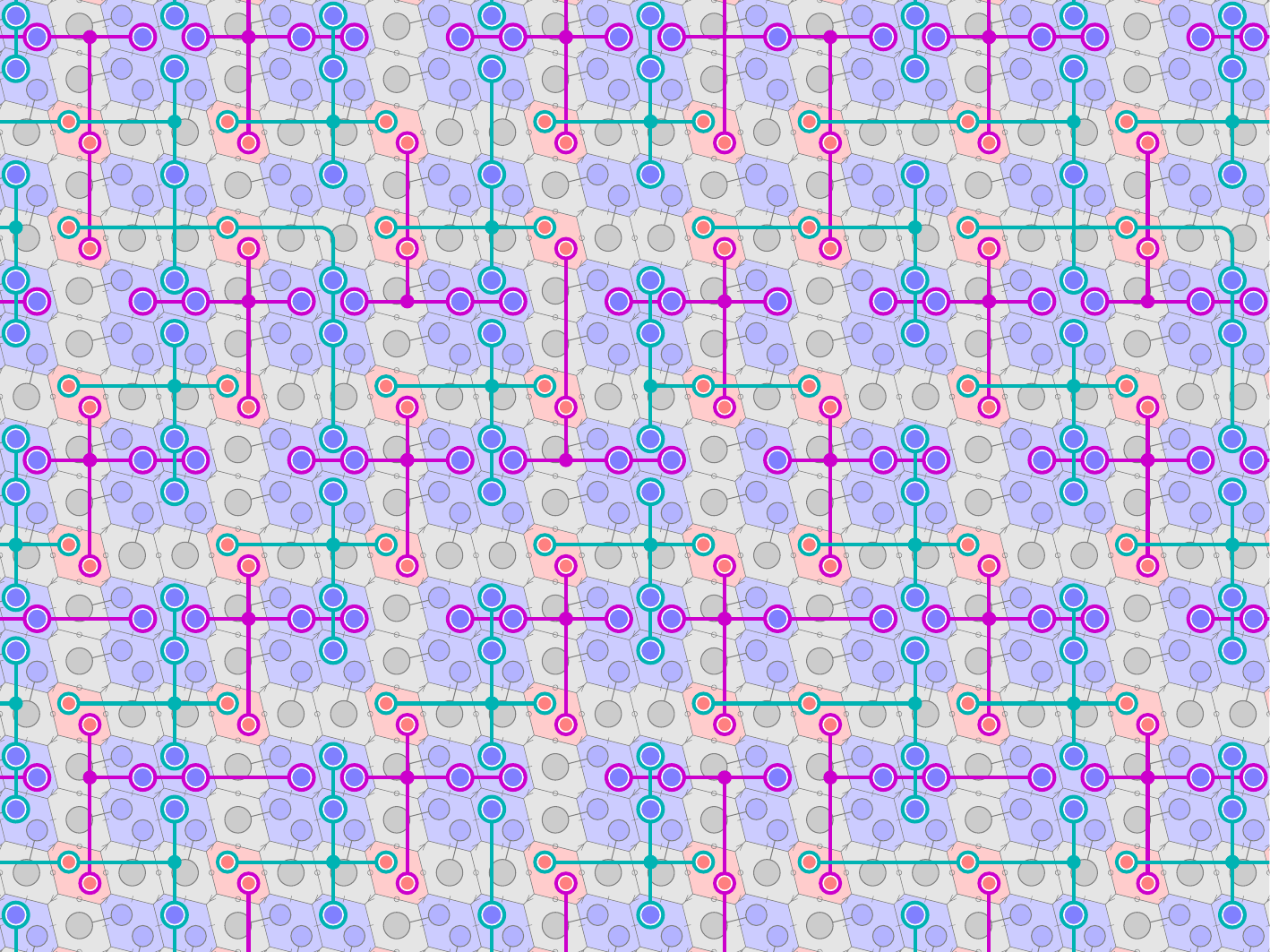}
\caption{Tiling by the $\SL(\sqrt{7} - 2)$-tile set in Figure~\ref{fig:sqrt7-tiles}}
\label{fig:sqrt7-tiling}
\end{figure}

Since $\alpha^{2} + 4\alpha - 3 = 0$ ($u = 4$, $v = -3$), we have
\[
\varphi(S) : \varphi(M) : \varphi(L) =
-3 : -1 : 2
\]
and our target is
\[
\AA/{\BDsim} = \left\{
2\tfrac{S}{2} + 3\tfrac{L}{2},~
M + \tfrac{L}{2},~
\Nut(S),~\Nut(M),~\Nut(L)
\right\}.
\]
Figure~\ref{fig:sqrt7-tiling} shows an example of the resulting tiling.
Applying Proposition~\ref{prop:add-sub}, we divide $\mathcal{L}_{1}$ into two sets $\mathcal{L}'\vee \mathcal{L}''$ with suitable natural densities since we use it for patch-tiles of two different types.
That is, we construct the two BD correspondences
\begin{align*}
\mathcal{S}_{1}&\stackrel{2}{\to} \Lambda'\stackrel{3}{\gets} \mathcal{L}',&
\mathcal{M}_{1}&\stackrel{1}{\to}\Lambda''\stackrel{1}{\gets} \mathcal{L}''.
\end{align*}
The latter reduces to a BD between one-dimensional lattices, and to the former we apply a cross-BD pair.

\begin{figure}[htb]\centering
\includegraphics[width = \linewidth,
pagebox = artbox, page = 1]
{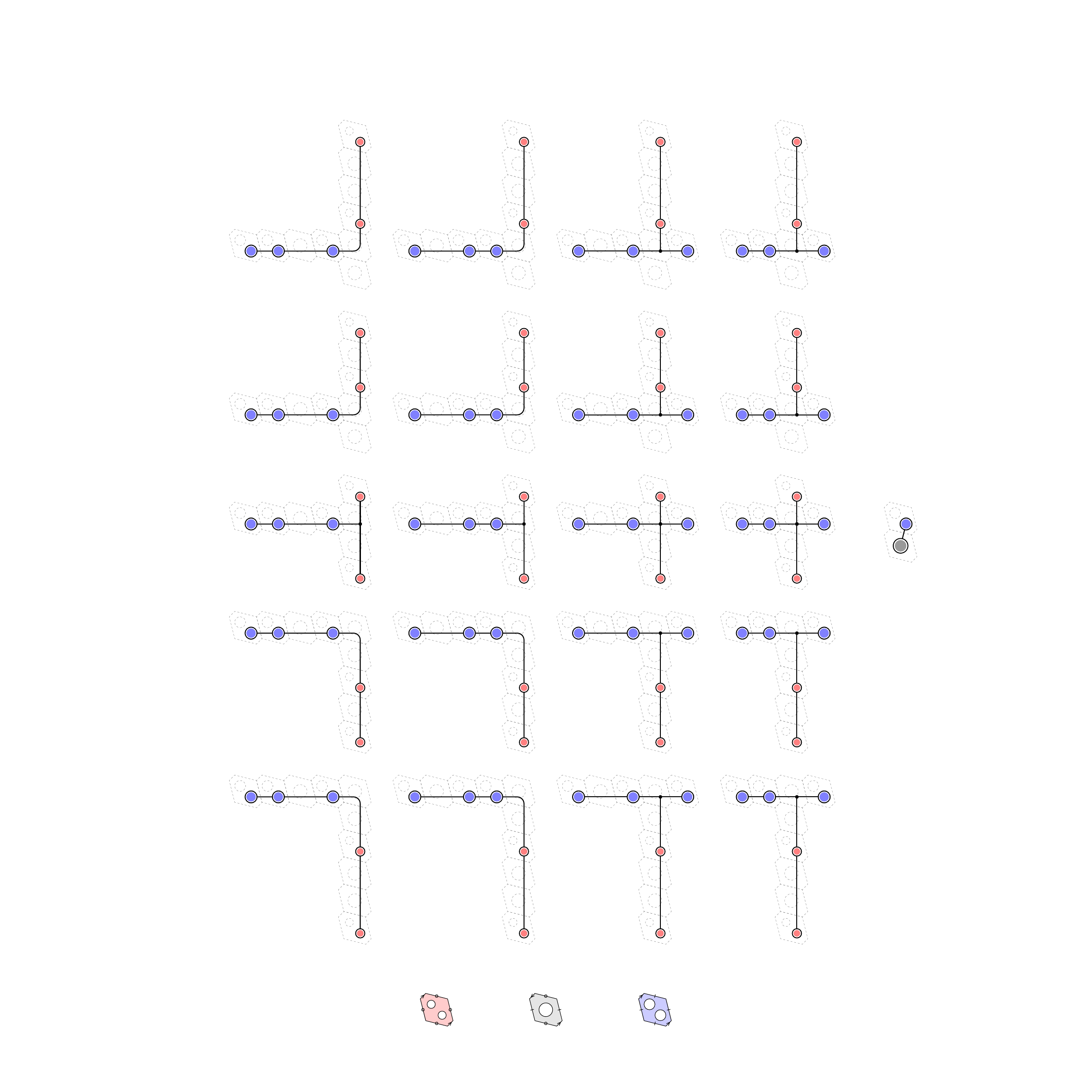}
\caption{$\SL(\sqrt{7} - 2)$-tile set composed of $24$ patch-tiles}
\label{fig:sqrt7-tiles}
\end{figure}

\subsection{Proof of Theorem~\ref{thm:lambda}}
\label{sec:lambda}

We give a proof of Theorem~\ref{thm:lambda} by constructing explicit tile sets.
In this paper, we only discuss the case where $N(\lambda) = -1$.
We fix an integer $h\ge 3$ and let
\[
\alpha:= \frac{1}{\lambda} =
\frac{-h+\sqrt{h^{2}+4}}{2} = [\bar{h}],
\]
which satisfies that $\alpha^{2} + h\alpha - 1 = 0$ and
\[
\frac{1}{h+1} < \alpha < \frac{1}{h}.
\]
Since
\[
\varphi(S) : \varphi(M) : \varphi(L) =
-1 : \frac{h - 2}{2} : h
\]
holds, after duplicating the Bolts as in \S\ref{sec:sqrt7}, one sees that the desired tile set is
\[
\AA/{\BDsim} = \left\{
\begin{aligned}
T' &= (h-2)\tfrac{S}{2} + M,\\
T''&= h\tfrac{S}{2} + \tfrac{L}{2},
\end{aligned}~
\Nut(S),~\Nut(M),~\Nut(L)
\right\}.
\]
Here again we decompose $\mathcal{S}_{1} = \mathcal{S}'\vee \mathcal{S}''$ and make $T'$ by a periodic BD correspondence $\mathcal{S}'\xrightarrow{h-2} \Lambda'\stackrel{1}{\gets} \mathcal{M}_{1}$ and $T''$ by a cross-BD pair $\mathcal{S}''\stackrel{h}{\to} \Lambda''\stackrel{1}{\gets} \mathcal{L}_{1}$.
As a result, a single prototile of type $T'$ and $2h + C$ prototiles of type $T''$ are produced, which gives the claim.
%
%

On the Bolt tiling $\mathscr{B}$ derived by a (symbolic) Sturmian lattice $(a_{i}, b_{j}, c_{k})$ of slope $\alpha$, 
we put
\begin{align*}
\mathcal{S}_{1}&=
\{B(j,k)\mid (b_{j}, c_{k}) = (0, 0)\},\\
\mathcal{M}_{1}&=
\{B(j,k)\mid (b_{j}, c_{k}) = (1, 0)\},\\
\mathcal{L}_{1}&=
\{B(j,k)\mid (b_{j}, c_{k}) = (1, 1)\}.
\end{align*}
To make a prototile of type $T' = (h-2)\frac{S}{2} + M$, we define
\[
\mathcal{S}':=
\bigsqcup_{j'=1}^{h-2}
\{B(j-j',k)\in \mathcal{S}_{1}\mid
B(j,k)\in \mathcal{M}_{1}\}
\]
and take a patch-tile
\[
T'(j,k):=
\bigsqcup_{j'=0}^{h-2} B(j-j',k)
\]
for integers $j$ and $k$ with $b_{j} = 1$, $c_{k} = 0$.
Then $T'(j,k)$ consists of tiles in $\mathcal{S}'\cup \mathcal{M}_{1}$ and is of type $T'$.

\begin{figure}[htb]\centering
\includegraphics[width = .4\linewidth,
pagebox = artbox, page = 1]
{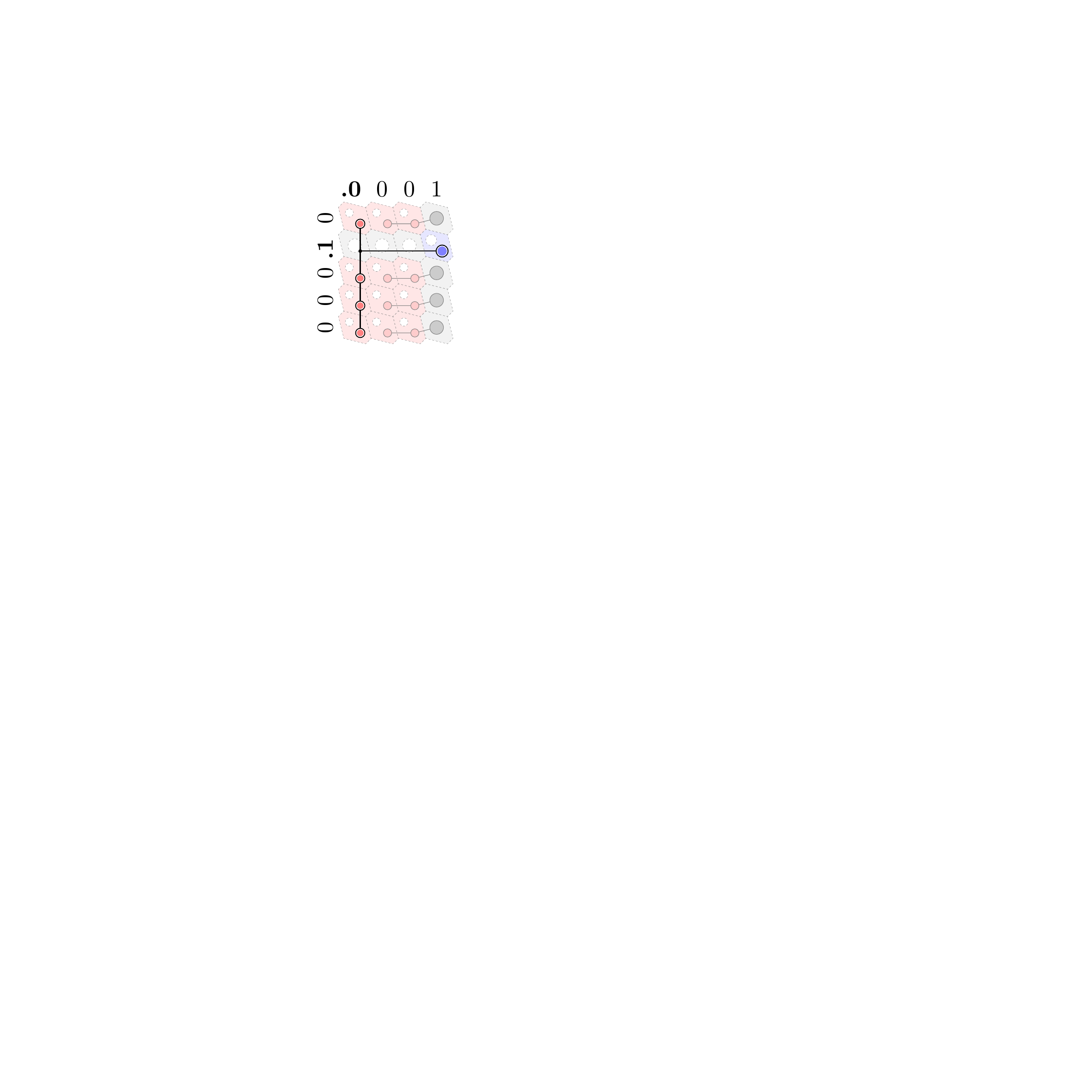}
\caption{A patch-tile of type $T'' = 4\frac{S}{2} + \frac{L}{2}$ characterized by a pair of factors $(v, w) = (b_{[0,3]}, c_{[-3,1]}) = (.0001, 000.10)\in \mathcal{V}\times \mathcal{W}$}
\label{fig:lambda-tiles}
\end{figure}

To make another type of patch-tiles $T'' = h\frac{S}{2} + \frac{L}{2}$, we apply a cross-BD pair to the remaining Bolts, say
\begin{align*}
\mathcal{S}'':=
\mathcal{S}_{1}\setminus \mathcal{S}'\BDsim
\hat{\mathcal{S}}''&=
\frac{1}{\alpha}\left(
\frac{1 - \alpha}{h\alpha}\Z\times
\frac{h\alpha}{1 - \alpha}
\frac{1}{h}\Z
\right),&
\mathcal{L}_{1}\BDsim
\hat{\mathcal{L}}_{1}&=
\frac{1}{\alpha}
(\Z\times \Z).
\end{align*}
We note that $\mu = (1 - \alpha)/h\alpha < 1$ holds since $1/(h+1) < \alpha$.
Let $P$ be a patch-tile of type $T''$, which is obtained by applying the cross-BD pair.
By acting translation, we may assume that $P$ has its nexus at $B(0, 0)$ without loss of generality.
Then $P$ has the form
\[
P = B(j_{1}, 0)\cup \bigcup_{n=1}^{h} B(0, k_{n})
\]
as the disjoint union of Bolts $B(j_{1}, 0)\in \mathcal{L}_{1}$ and $B(0, k_{n})\in \mathcal{S}''$ with $k_{1} < k_{2} < \dots < k_{h}$.
We take the minimal patch
\[
\mathscr{X}:=
\{B(j,0)\mid j\in J = [-j',j'']\}\cup
\{B(0,k)\mid k\in K = [-k',k'']\}
\]
that covers $P$.
We wish to enumerate all possible two factors $v = b_{J}$ and $w = c_{K}$ that uniquely determine the patch $P$.
See Figure~\ref{fig:lambda-tiles} for an example when $h = 4$.

By construction, both $v$ and $w$ are minimal factors of a Sturmian word of slope $\alpha$ and they satisfy
\begin{align*}
b_{[0,h-2]}&= 00^{h-2},&
\abs{v}_{1}&= 1,&
c_{0}&= 1,&
\abs{w}_{0}&= h.
\end{align*}
Equivalently, $(v, w)\in \mathcal{V}\times \mathcal{W}$, where
\begin{align*}
\mathcal{V}&:= \left\{\begin{array}
{r@{\str0}l}
10 & 0^{h-2}\\
 1 & 0^{h-2}\\
   & 0^{h-2}1\\
   & 0^{h-2}01
\end{array}\right\},&
\mathcal{W}&:= \left\{\begin{array}
{r@{\str1}l}
010^{h-1} & \\
 00^{h-1} & \\
  0^{h-1} & 0^{1}\\
  0^{h-2} & 0^{2}\\
\multicolumn{2}{c}{\vdots}\\
  0^{2} & 0^{h-2}\\
  0^{1} & 0^{h-1}\\
        & 0^{h-1}0\\
        & 0^{h-1}10
\end{array}\right\}.
\end{align*}
Since $\Card(\mathcal{V}) = 4$ and $\Card(\mathcal{W}) = h + 3$, the patch-tiles $T''$ can be collected as $\Card(\mathcal{V}\times \mathcal{W}) = 4h + 12$ up to translation.
Moreover, since $\mathcal{W}$ is closed under reflection (with axis $c_{0} = \str{1}$), we find that there are $2h + C$ of them up to isometry, which gives the claim for the case $N(\lambda) = -1$.

\section{Examples II}
\label{sec:ex2}

In this section, we present additional examples where we tried to minimize
the number of patch-tiles.
Indeed, the previous examples used a general method applicable to any quadratic slope, resulting in a large number of tiles and inefficient counting.

\subsection{\texorpdfstring
{Another example for $\alpha = \sqrt{2} - 1$}
{Another example for \textbackslash alpha = \textbackslash sqrt\{2\} - 1}}
\label{sec:sqrt2-2}

\begin{figure}[htb]\centering
\includegraphics[pagebox = artbox, 
width = \linewidth, page = 1]
{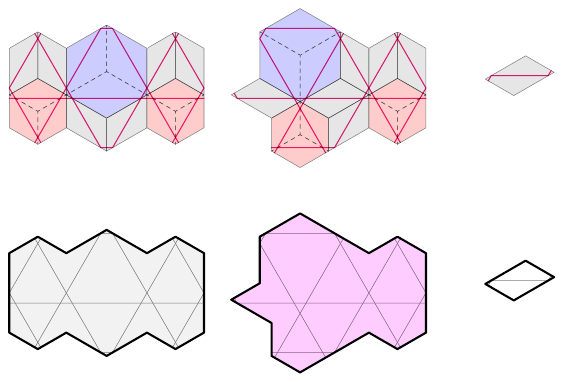}
\caption{Another example of $\SL(\sqrt{2} - 1)$-tile sets}
\label{fig:ex1-tiles}
\end{figure}

\begin{figure}[htb]\centering
\includegraphics[
width = \linewidth, page = 1]
{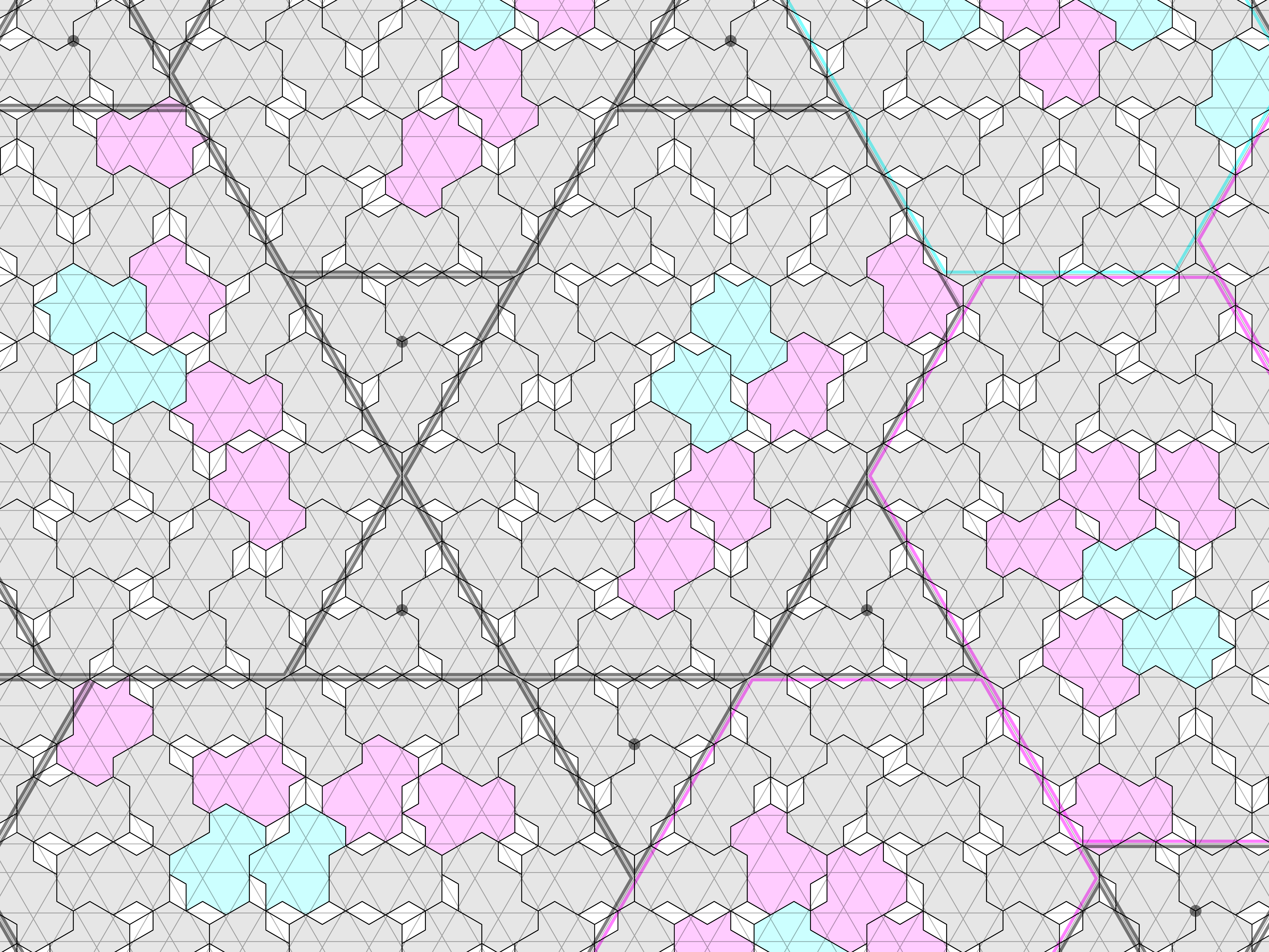}
\caption{Tiling by $\SL(\sqrt{2} - 1)$-tile set in Figure~\ref{fig:ex1-tiles}}
\label{fig:ex1-tiling}
\end{figure}

We revisit the slope $\alpha = \sqrt{2} - 1$.
We produce another tile set and we need only three patch-tiles up to isometry.
The resulting patch-tiles are in Figure~\ref{fig:ex1-tiles}.
{\color{HMDcomment}Figure~\ref{fig:ex1-tiling} shows a tiling; here a cyan tile is the mirror image of a magenta one.}

The following two ideas are central to the construction:
\begin{itemize}
\item
We use another type of cells, called ``Voronoi cells''.
Their tiling has a complicated structure described as the sum of three rotated lattices, but they instead exhibit higher symmetry.
\item
When constructing a necessary BD map, we use the self-similarity of Sturmian lattices (see \S\ref{sec:negative}).
We do not yet know how to generalize this method.
\end{itemize}

\begin{figure}[htb]\centering
\subfigure[Parallelogram cells $H_{a}(j,k)$]{%
\includegraphics[page = 1,
width = .48\linewidth]
{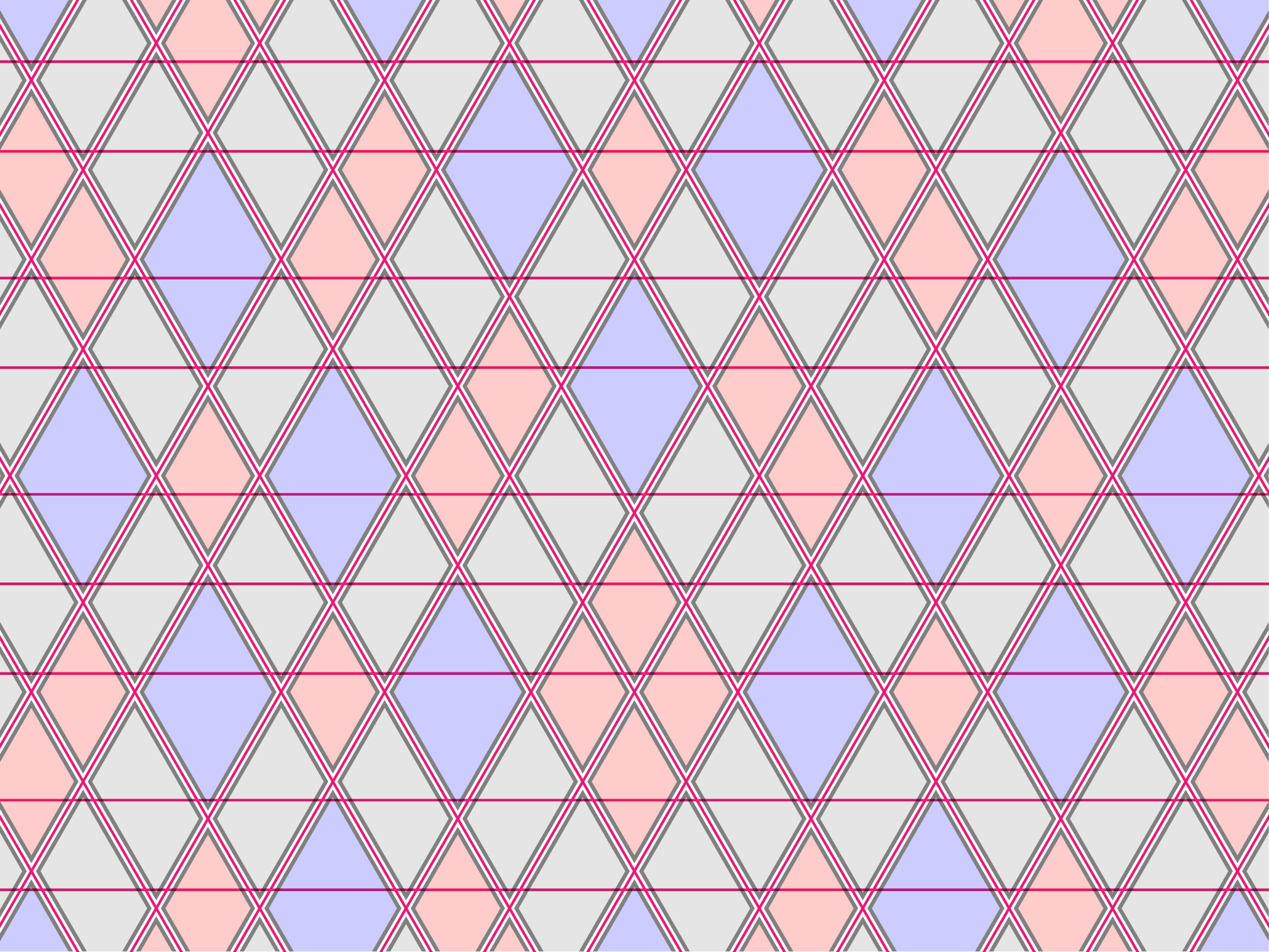}
}\hfill%
\subfigure[Voronoi cells $\VC_{a}(j,k)$]{%
\includegraphics[page = 1,
width = .48\linewidth]
{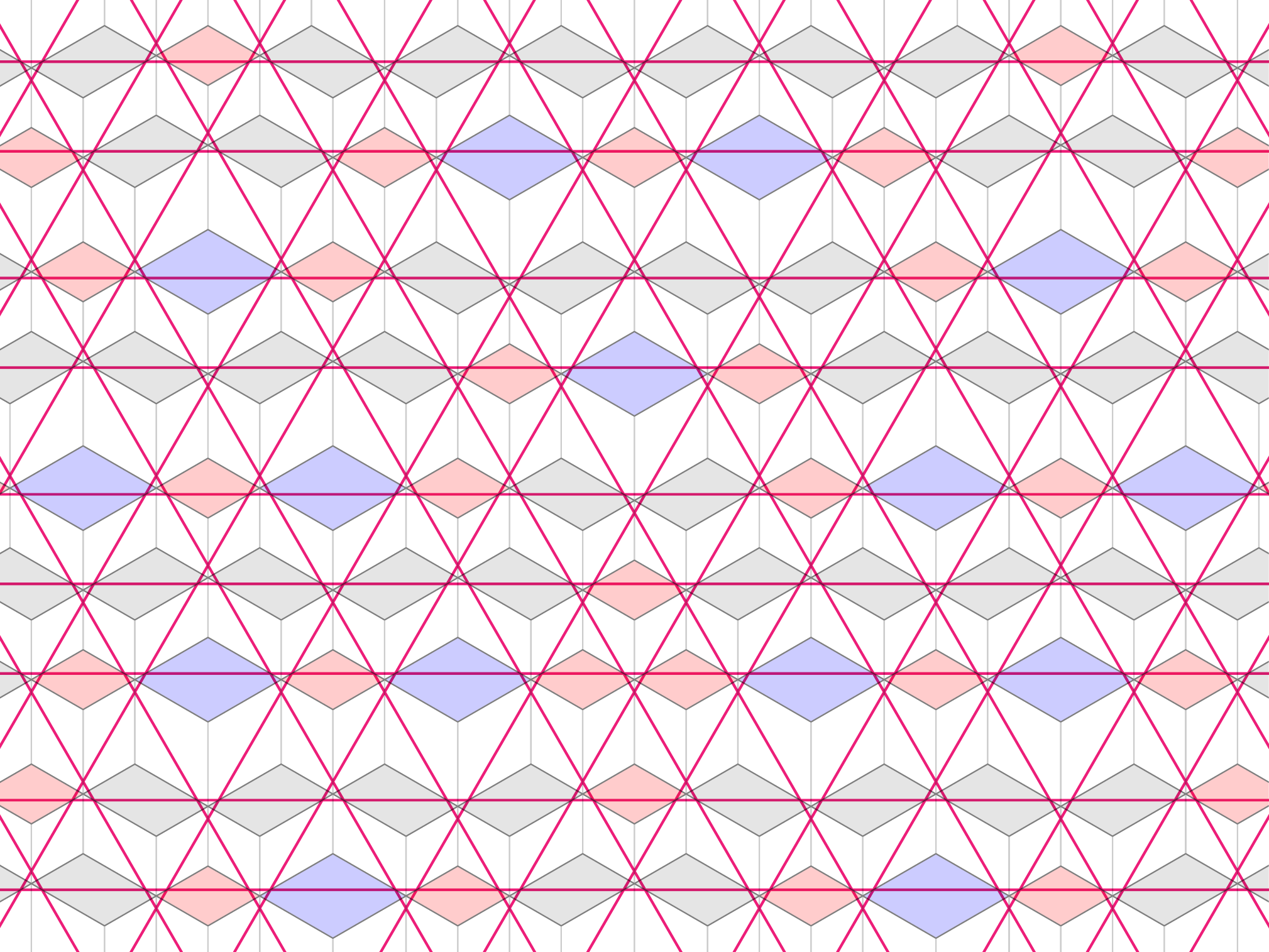}
}
\caption{Comparison of parallelogram cells and Voronoi cells}
\label{fig2-isometric1}
\end{figure}

We define Voronoi cells.
Let $(a(i), b(j), c(k))\in \SL(\kappa, \alpha)$ be a $2$-color Sturmian lattice.
For integers $j$ and $k$, we define (the support of) the Voronoi cell $\VC_{a}(j, k)$ as follows: let $H_{a}(j, k)\subset \R^{2}$ be the parallelogram surrounded by four lines:
\begin{align*}
b&= b(j),&
b&= b(j+1),&
c&= c(k),&
c&= c(k+1)
\end{align*}
as depicted in Figure~\ref{fig2-isometric1}(a).
Then the horizontal line $a = a(i)$ with $i + j + k + 1 = 0$ crosses almost diagonally across $H_{a}(j, k)$.
By acting translation, we set $a(i) = 0$ without loss of generality.
We define
\[
\supp(\VC_{a}(j, k)):= \begin{pmatrix}
1 & 0\\ 0 & 1/3
\end{pmatrix}H_{a}(j, k),
\]
whose vertices are located at the centroids of the four triangles generated by the five lines above, see Figure~\ref{fig2-isometric1}(b).
We also define $\supp(\VC_{b}(k, i))$ and $\supp(\VC_{c}(i, j))$ in a similar way, and we put
\[
\supp(\mathcal{VC}):=
\supp(\mathcal{VC}_{a})\cup
\supp(\mathcal{VC}_{b})\cup
\supp(\mathcal{VC}_{c})
\]
where
\begin{align*}
\supp(\mathcal{VC}_{a})&:=
\{\supp(\VC_{a}(j, k))\subset \R^{2}\mid j, k\in \Z\},\\
\supp(\mathcal{VC}_{b})&:=
\{\supp(\VC_{b}(k, i))\subset \R^{2}\mid k, i\in \Z\},\\
\supp(\mathcal{VC}_{c})&:=
\{\supp(\VC_{c}(i, j))\subset \R^{2}\mid i, j\in \Z\}.
\end{align*}

The Voronoi cells have Ammann bars (see \S\ref{sec:aperiodic}), which describe our matching rule.
Let
\begin{align*}
\mathfrak{A}&:= \bigsqcup_{i\in \Z}
\{a = a(i)\},&
\mathfrak{B}&:= \bigsqcup_{j\in \Z}
\{b = b(j)\},&
\mathfrak{C}&:= \bigsqcup_{k\in \Z}
\{c = c(k)\}
\end{align*}
be the collection of lines in each direction.
We draw the inherited pattern
\[
\ell(V):= (\mathfrak{A}\cup \mathfrak{B}\cup \mathfrak{C})\cap \supp(V)
\]
on each $\supp(V)\in \supp(\mathcal{VC})$.
We define \emph{Voronoi cells} $V:= (\supp(V), \ell(V))$ with Ammann bars and denote by
\[
\mathcal{VC} = \mathcal{VC}_{a}\cup \mathcal{VC}_{b}\cup \mathcal{VC}_{c}
\]
the set of Voronoi cells on the Sturmian lattice.

\begin{rem}
\label{rem:Voronoi}
The name ``Voronoi cells'' actually derives from the fact that $\mathcal{VC}$ is obtained as a Voronoi diagram with respect to open segments on the plane.
For $i + j + k = 0$, we define $\mathfrak{T}(i,j,k)\subset \mathfrak{A}\cup \mathfrak{B}\cup \mathfrak{C}$ as the boundary of the tiny triangle surrounded by three lines $a = a(i)$, $b = b(j)$, and $c = c(k)$.
Then the complement $\mathfrak{S}:= (\mathfrak{A}\cup \mathfrak{B}\cup \mathfrak{C})\setminus \mathfrak{T}$ of
\[
\mathfrak{T}:= \bigcup_{i + j + k = 0}
\mathfrak{T}(i,j,k)
\]
is a collection of open segments.
The Voronoi-cell tiling is equivalent with the Voronoi diagram with respect to the connected components in $\mathfrak{S}$.
\end{rem}

By definition of the Voronoi cells, it is plain to see that there are three types, up to isometry.
We say that $\VC_{a}(j, k)\in \mathcal{VC}_{a}$ is
\begin{itemize}
\item of type $S_{a}$ if $b_{j} = c_{k} = 0$,
\item of type $M_{a}$ if $b_{j}\ne c_{k}$, and
\item of type $L_{a}$ if $b_{j} = c_{k} = 1$.
\end{itemize}
We also define types $S_{b}$, $M_{b}$, and $L_{b}$ for $\mathcal{VC}_{b}$, and types $S_{c}$, $M_{c}$, and $L_{c}$ for $\mathcal{VC}_{c}$.
We denote by $\mathcal{S}_{*}$, $\mathcal{M}_{*}$, and $\mathcal{L}_{*}$ (${*}\in \{a, b, c\}$) the set of Voronoi cells of types $S_{*}$, $M_{*}$, and $L_{*}$.
If we tile the plane by the Voronoi cells $\{S_{*}, M_{*}, L_{*}\}$ following the matching rule, then the Ammann bars must form a Sturmian lattice by Lemma~\ref{lem:SL}.
Since the natural densities can be again calculated as
\[
\delta(\mathcal{S}_{*}) :
\delta(\mathcal{L}_{*}) =
(1 - \alpha)^{2} :
\alpha^{2} = 2 : 1,
\]
we can apply the same calculation as for hexagonal cells, and we may make a $\SL(\alpha)$-tile set of type $\{2S_{*} + L_{*}, M_{*}\}$.
For higher symmetry, we combine three Voronoi cells to prepare hexagonal tiles
\begin{align*}
S&= S_{a} + S_{b} + S_{c},&
L&= L_{a} + L_{b} + L_{c}.
\end{align*}
We denote by $\mathcal{S}$ and $\mathcal{L}$ the sets of the centroids of cells $S$ and $L$, respectively.

\begin{figure}[htb]\centering
\includegraphics[width = \linewidth,
pagebox = artbox, page = 1]
{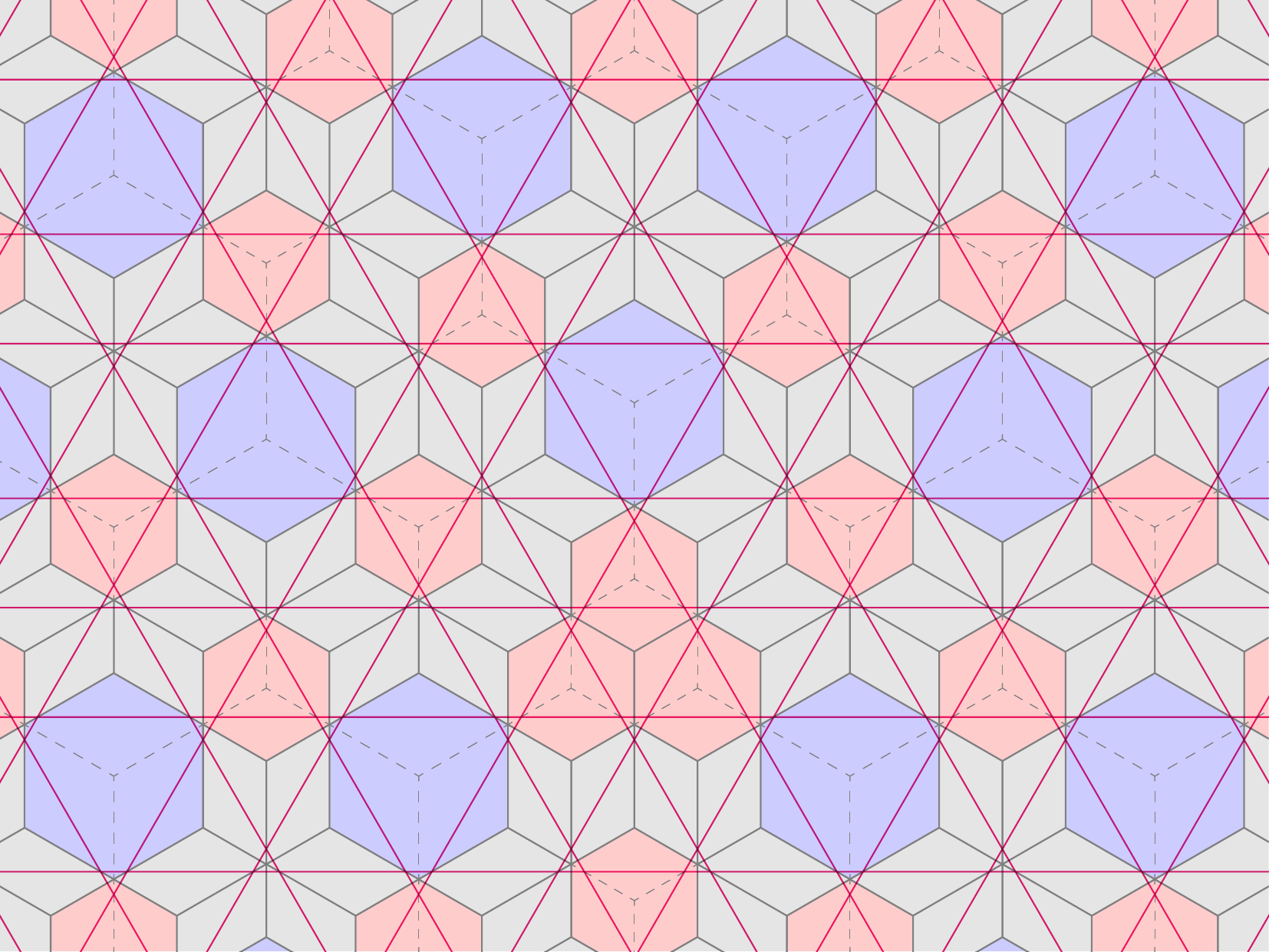}
\caption{Voronoi-cell tiling by $\{S, M_{a}, M_{b}, M_{c}, L\}$}
\label{fig2-isometric2}
\end{figure}

\begin{figure}[htb]\centering
\subfigure[$\SL(\kappa, \alpha)$]{%
\includegraphics[width = .48\linewidth]{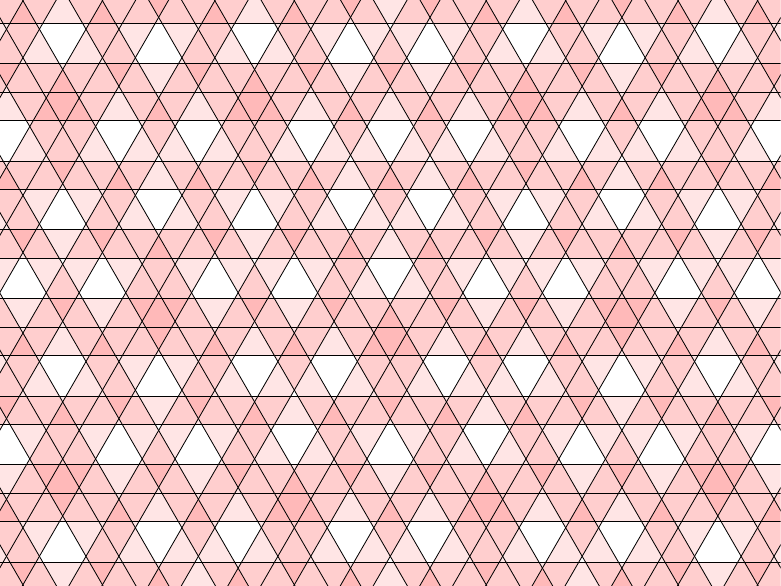}
}\hfill%
\subfigure[$\Psi^{2}\SL(\kappa, \alpha) = \kappa^{2}\SL(\kappa, \alpha)$]{%
\includegraphics[width = .48\linewidth]{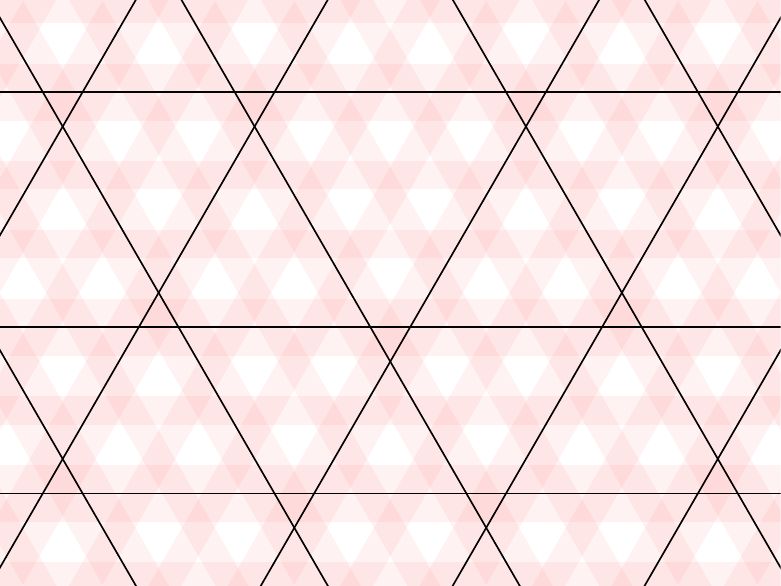}
}
\caption{Super Sturmian lattice on level $2$ for $\alpha = \sqrt{2} - 1 = [\overline{2, 2}]$}
\label{fig:ex1-superSL}
\end{figure}

{\color{HMDcomment}Since the slope $\alpha = \sqrt{2} - 1 = [\bar{2}]$ is purely periodic, taking the passage $\kappa = 1/\alpha = \sqrt{2} + 1$ gives $\Psi\SL(\kappa, \alpha) = \kappa\SL(\kappa, \alpha)$.
Note further that, for a Sturmian lattice $(a(i), b(j), c(k))\in \SL(\kappa, \alpha)$, the lattice $(A(I), B(J), C(K)) = \Psi^{2}(a(i), b(j), c(k))\in \kappa^{2}\SL(\kappa, \alpha)$ is characterized by}
\begin{align*}
\{A(I)\in \R\mid I\in \Z\}&=
\{a(i)\in \R\mid a_{i-1}a_{i} = 00\},\\
\{B(J)\in \R\mid J\in \Z\}&=
\{b(j)\in \R\mid b_{j-1}b_{j} = 00\},\\
\{C(K)\in \R\mid K\in \Z\}&=
\{c(k)\in \R\mid c_{k-1}c_{k} = 00\}.
\end{align*}
{\color{HMDcomment}In particular, $\Psi^{2}(a(i), b(j), c(k))$ is a Sturmian sub-lattice of $(a(i), b(j), c(k))$; see Figure~\ref{fig:ex1-superSL}.}

To get a desired tile set, we may construct a BD map $\mathcal{S}\stackrel{2}{\to} \mathcal{L}$.
{\color{HMDcomment}We first prepare four kinds of polygonal tiles
whose boundaries consist of the lines constituting a Sturmian lattice $(a(i), b(j), c(k))\in \SL(\kappa, \alpha)$, and describe the self-similarity of the Sturmian lattice as a self-similar tiling.}
\begin{itemize}
\item Two hexagons $\SLH_{1}$ and $\SLH_{2}$.
They include two $S$'s and a single $L$.
\item Two triangles $\SLT_{M}$ and $\SLT_{0}$.
\end{itemize}
Figure~\ref{fig:ex1-substitution} depicts them.
Applying $\Psi^{2}$ to these tiles, the images can be written by
\begin{align*}
\Psi^{2}(\SLH_{1})&= 23\SLH_{1} + 6\SLH_{2} + 46\SLT_{M} + 10\SLT_{0},\\
\Psi^{2}(\SLH_{2})&= 22\SLH_{1} + 7\SLH_{2} + 46\SLT_{M} + 10\SLT_{0},\\
\Psi^{2}(\SLT_{M})&= 3\SLH_{1} + 6\SLT_{M} + \SLT_{0},\\
\Psi^{2}(\SLT_{0})&= \SLT_{M},
\end{align*}
which define a self-similar tiling that is compatible with the underlying Sturmian lattice.
By infinitely iterating $\Psi^{2}$, we obtain a BD map $\mathcal{S}\stackrel{2}{\to}\mathcal{L}$.
This may also be viewed as an application of the extension theorem for tilings; see, for example, Grünbaum and Shephard~\cite{GS}.
To get a tiling as in Figure~\ref{fig:ex1-tiling}, we substitute hexagons $\SLH_{1}$ and $\SLH_{2}$ for the patch-tiles including $2S + L$, and we may fill the gap with Voronoi cells $M_{*}$.

\begin{figure}[htb]\centering
\includegraphics[pagebox = artbox, 
width = \linewidth, page = 4]
{Example1.pdf}
\caption{Substitution rules (from left: $\SLH_{1}$, $\SLH_{2}$, $\SLT_{M}$, and $\SLT_{0}$)}
\label{fig:ex1-substitution}
\end{figure}

Note that this tile set has local interchangeable pairs.
Figure~\ref{fig:ex1-substitution-pair} shows examples.
If $\{\P_{1}, \P_{2}\}$ is a local interchangeable pair, 
then so is $\{\Psi^{2n}(\P_{1}), \Psi^{2n}(\P_{2})\}$ for any $n\ge 1$.
In particular,
there are uncountably many tilings that fit a given Sturmian lattice, and the associated tiling space
has positive topological
entropy, see \cite{DGG:17}.

\begin{figure}[htb]\centering
\includegraphics[pagebox = artbox, 
width = \linewidth, page = 5]
{Example1.pdf}
\caption{Local interchangeable pairs}
\label{fig:ex1-substitution-pair}
\end{figure}

\subsection{Smith Turtle revisited}
\label{sec:exTurtle}

\begin{figure}[htb]\centering\hfill
\subfigure[Cell $S$]{%
\includegraphics[page = 10,
width = .30\linewidth]
{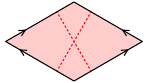}
}\hfill%
\subfigure[Cell $M$]{%
\includegraphics[page = 11,
width = .30\linewidth]
{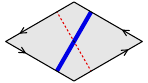}
}\hfill%
\subfigure[Cell $L$]{%
\includegraphics[page = 12,
width = .30\linewidth]
{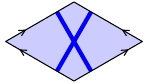}
}\hfill\mbox{}
\caption{Three rhombi with two-colored Ammann bars}
\label{fig:ex2-Voronoi}
\end{figure}

In this subsection we try to interpret Smith Turtle in \cite{SMKGS:23_1} using our method.
{\color{HMDcomment}We prepare the congruent rhombi $S$, $M$, and $L$ shown in Figure~\ref{fig:ex2-Voronoi}; their interior angles are $\pi/3$ and $2\pi/3$, and their four edges are oriented.
The segments of length $1$ on the support pass through the centroid of the rhombus and are perpendicular to the edges; two red segments are drawn on $S$, one red and one blue on $M$, and two blue segments on $L$.
We tile the plane with these three kinds of tiles, glued edge-to-edge respecting the orientations, and use the two-colored segments as Ammann bars.}
Then, the Ammann bars must form a kagome lattice, and focusing on only one of the colors produces a geometric Sturmian lattice, see Figure~\ref{fig:ex3-tiling}(a).
If we set $q\in [0, 1]$ to be the natural density of blue lines, then the natural densities of each type of cell are calculated as
\[
\delta(\mathcal{S}) :
\delta(\mathcal{M}) :
\delta(\mathcal{L}) =
(1 - q)^{2} :
2q(1 - q) :
q^{2}.
\]
{\color{HMDcomment}Therefore, by using the tiles of Figure~\ref{fig:ex2-Voronoi} instead of the Nuts and Bolts of \S\ref{sec:Voronoi}, we can also construct an $\SL(q)$-tile set for any quadratic irrational $q\in (0, 1)$.}

\begin{figure}[htb]\centering
\subfigure[original]{%
\includegraphics[pagebox = artbox, 
width = .48\linewidth, page = 2]
{Example3.pdf}}\hfill%
\subfigure[transformed]{%
\includegraphics[pagebox = artbox, 
width = .48\linewidth, page = 3]
{Example3.pdf}}%
\caption{Transformed rhombic cells}
\label{fig:ex3-tiling}
\end{figure}


{\color{HMDcomment}As long as the Sturmian lattice is still enforced, deforming the three kinds of tiles while preserving their areas and the total length of the Ammann bars does not change the values of the natural densities.
In particular, using the pentagons appearing in Figure~\ref{fig:ex3-tiling}(b), a Turtle and its flipped version are described as the patch-tiles $3S + 2M$ and $3L + 2M$ respectively, as shown in Figure~\ref{fig3-Turtle}.
The two-colored Ammann bars induced from the transformed cells correspond to the Golden Ammann bars introduced in \cite{Akiyama-Araki:23} and their complement, and it is shown there that they enforce a Sturmian lattice.}

\begin{figure}[htb]\centering
\includegraphics[pagebox = artbox, 
width = .5\linewidth, page = 4]
{Example3.pdf}
\caption{Smith Turtle as the union of transformed cells}
\label{fig3-Turtle}
\end{figure}

\begin{figure}[htb]\centering
\includegraphics[pagebox = artbox, 
width = \linewidth, page = 1]
{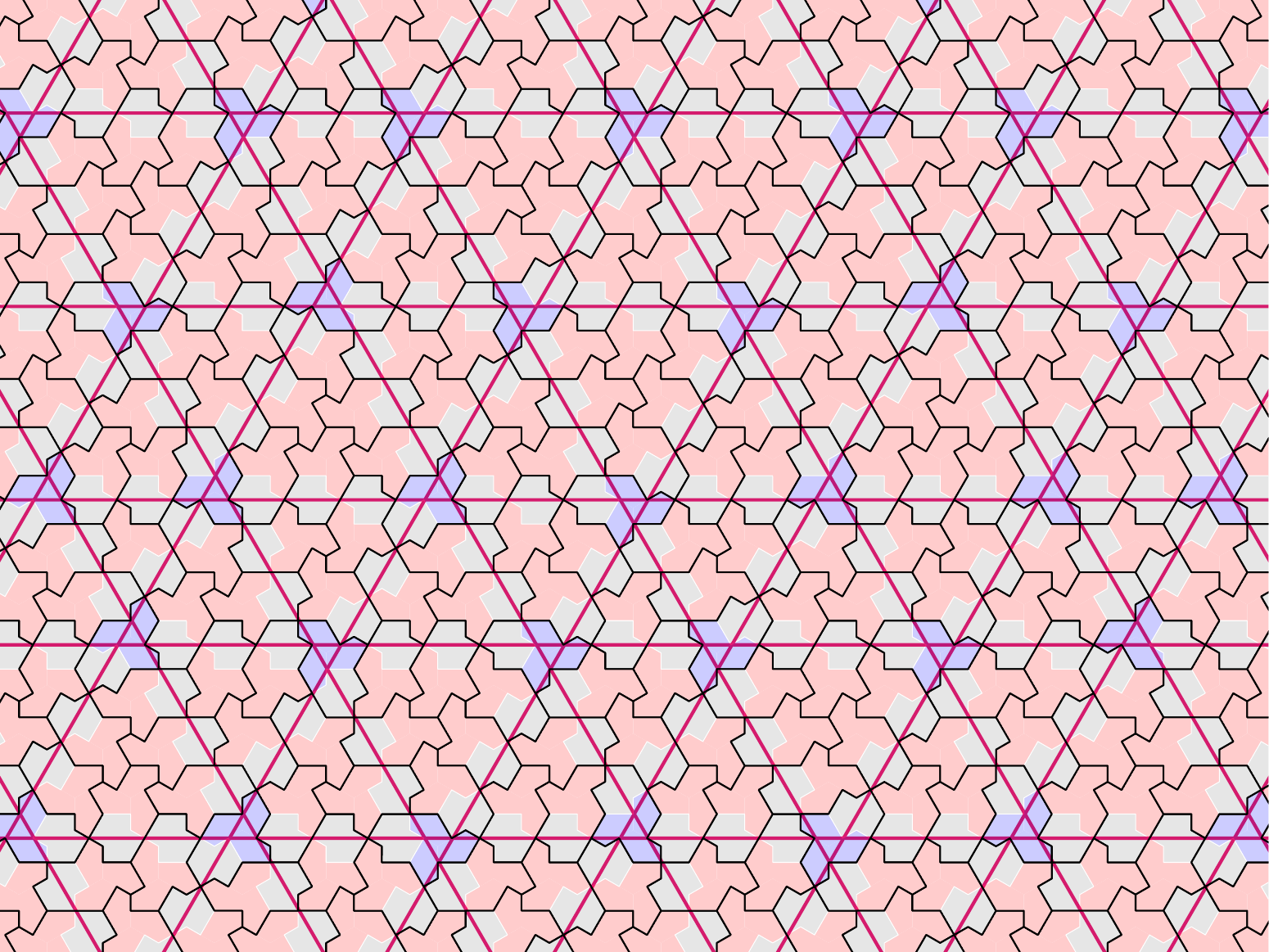}
\caption{Tiling by Turtle as in Figure~\ref{fig3-Turtle}}
\label{fig3-Turtle-tiling}
\end{figure}

{\color{HMDcomment}The natural densities $\delta':= \delta(3S + 2M)$ and $\delta'' = \delta(3L + 2M)$ of the two kinds of Turtle in the tiling satisfy the following equalities.}
\begin{align*}
(1 - q)^{2}&= 3\delta',&
2q(1 - q)&= 2\delta' + 2\delta'',&
q^{2}&= 3\delta''.
\end{align*}
{\color{HMDcomment}Therefore we have}
\[
q = \frac{5\pm \sqrt{5}}{10},
\]
{\color{HMDcomment}which provides another view of the aperiodicity proof in \cite{Akiyama-Araki:23}.}

\section{Open problems}
\label{Open}

We give a list of intriguing problems.

\begin{itemize}

\item
Can we obtain another aperiodic monotile? We did not put much effort in
minimizing the number of tiles in this paper. There should be
a lot of freedom left in the usage of bounded displacement equivalence. 

\item Is there a way to construct an aperiodic tile set by a Sturmian lattice
of non-quadratic slope? What about cubic irrational slopes? 
Note that Sturmian lattices exist for any irrational slope,
and we do not essentially use the self-similar structure in our construction. 

\item
Can we make the tiles homeomorphic to a disk? Extending our method of construction, 
we can make 
all tiles connected but it is not clear if we can make their interiors 
connected.
Can we simplify our matching rule? When can we use super SAB instead of the original SAB? When can we eliminate SAB? 
Bounded displacement equivalence does not take care the topology of resulting tiles and SAB's.
\item
Regarding Theorem~\ref{thm:lambda}, can the upper bound on the number of tiles be made uniformly bounded in $h$?
\end{itemize}

\section*{Acknowledgments}

This research was partially supported by JSPS grants (21H00989, 20K03528,
24K06662). 
We are grateful to Yotam Smilansky who
drew our attention to Duneau--Oguey~\cite{DO:91} for the explicit construction
of the bounded displacement map.
Our original
draft is drastically improved 
by the idea of their construction. We are also deeply indebted to 
S\'ebastien Labb\'e who reminded us of the importance of the 
classification of bi-infinite balanced words of 
Markoff refined by \cite{Reutenauer:06} (see also \cite{Labbe:25}). 


\end{document}